\documentclass[a4paper]{amsart}

\usepackage[mathscr]{eucal}
\usepackage{amsmath,amsfonts,amssymb}
\usepackage{amsthm}
\usepackage{mathrsfs}
\usepackage{bm}
\usepackage{cite}
\usepackage[dvips]{graphicx,color}
\graphicspath{{./graphs/}}
\usepackage{fancyhdr}

\newcommand{\mathsym}[1]{{}}

\theoremstyle{plain}

\newtheorem{theorem}{Theorem}[section]
\newtheorem{proposition}[theorem]{Proposition}
\newtheorem{lemma}[theorem]{Lemma}

\newtheorem{conjecture}{Conjecture}[section]
\theoremstyle{definition}
\newtheorem{definition}{Definition}[section]
\newtheorem{remark}{Remark}[section]

\begin{document}
\title[Construction of spherical designs from finite graph]{A construction of spherical designs from finite graphs \\ with the theory of crystal lattice}
\author{Junichi Shigezumi}

\maketitle \vspace{-0.1in}
\begin{center}
Graduate School of Mathematics Kyushu University\\
Hakozaki 6-10-1 Higashi-ku, Fukuoka, 812-8581 Japan\\
{\it E-mail address} : j.shigezumi@math.kyushu-u.ac.jp \vspace{-0.05in}
\end{center} \quad

\begin{quote}
{\small\bfseries Abstract.}
We want to introduce a construction of spherical designs from finite graphs with the theory of crystal lattice. We start from a finite graph, and we consider standard realization of the crystal lattices as the maximal Abelian covering of the graph. Then, we take the set of vectors which form the crystal lattice. If every vector has the same norm, then we can consider a finite set on Euclidean sphere, and then we get a spherical design.

In this paper, we observe the results by numerical calculations. We tried constructing vectors from various finite graphs, strongly regular graphs, distance regular graphs, and so on. We also introduce some facts and conjectures.\\  \vspace{-0.15in}

\vspace{-0.15in}

 \vspace{0.15in}
\end{quote}

\section{Introduction}

Let $S^{d-1} := \{ (x_1, \ldots, x_d) \in \mathbb{R}^d \, ; \, x_1^2 + \cdots + x_d^2 = 1 \}$ be the Euclidean sphere for $d \geqslant 1$.

\begin{definition}[Spherical design \cite{DGS}]\label{def-design}
Let $X$ be a non-empty finite set on the Euclidean sphere $S^{d-1}$, and let $t$ be a positive integer. $X$ is called a spherical $t$-design if
\begin{equation}
\frac{1}{| S^{d-1} |} \int_{S^{d-1}} f(\xi) \, d \xi \ = \ \frac{1}{| X |} \hspace{0.05in} \sum_{\xi \in X} \hspace{0.05in} f(\xi) \label{eq-design}
\end{equation}
for every polynomial $f(x) = f(x_1, \ldots, x_d)$ of degree at most $t$.
\end{definition}

Here, the left-hand side of the above equation (\ref{eq-design}) means the average on the sphere $S^{d-1}$, and the right-hand side means the average on the finite subset $X$. Thus, if $X$ is a spherical design, then $X$ gives a certain approximation of the sphere $S^{d-1}$. We also have another condition for spherical design:
\begin{proposition}[cf. \cite{DGS}]\label{prop-con-design}
Let $t$ be a positive integer. $X \subset S^{d-1}$ is a spherical $t$-design if and only if
\begin{equation}
\sum_{\xi \in X} \hspace{0.05in} f(\xi) = 0 \label{eq-design1}
\end{equation}
for every polynomial $f \in \text{\upshape Harm}_1(S^{d-1}) \oplus \cdots \oplus \text{\upshape Harm}_t(S^{d-1})$, where $\text{\upshape Harm}_i(S^{d-1})$ is the space of homogeneous harmonic polynomials of degree $i$ on $S^{d-1}$.
\end{proposition}

Let us denote $X = \{ u_i \, ; \, 1 \leqslant i \leqslant n \}$, and $\sqrt{d/n} \ u_k = (w_{1, k}, \ldots, w_{d, k})$ for every $1 \leqslant k \leqslant n$. We have the polynomials $x_i$ ($1 \leqslant i \leqslant d$) as basis of $\text{Harm}_1(S^{d-1})$, and the polynomials $x_i x_j$, $x_i^2 - x_j^2$ ($1 \leqslant i \leqslant d$) as basis of $\text{Harm}_2(S^{d-1})$. Thus, by the above proposition, the finite set $X \subset S^{d-1}$ is a spherical $2$-design if and only if $\{ w_{i, k} \}_{1 \leqslant i \leqslant d, 1 \leqslant k \leqslant n}$ satisfy the following: (cf. \cite{M})
\begin{allowdisplaybreaks}
\begin{align}
\tag{$C0$} &\sum_{i=1}^d w_{i, k}^2 = \frac{d}{n}, \qquad 1 \leqslant k \leqslant n.\\
\tag{$C1$} &\sum_{k=1}^n w_{i, k} = 0, \qquad 1 \leqslant i \leqslant d.\\
\tag{$C2$} &\sum_{k=1}^n w_{i, k} w_{j, k} = 0, \qquad 1 \leqslant i < j \leqslant d.\\
\tag{$C3$} &\sum_{k=1}^n w_{i, k}^2 = 1, \qquad 1 \leqslant i \leqslant d.
\end{align}
\end{allowdisplaybreaks}

Now, we consider a finite graph $X_0 = (V_0, E_0)$, and we consider standard realization of the crystal lattices as the maximal Abelian covering of the graph. Then we take the set of vectors $X^0 := \{P(e_0) \, ; \, e_0 \in E_0\}$ which form the crystal lattice, and then we define $X$ the set of normalized vectors of $X^0$. In details, we can see in Section \ref{sec-sr-lattice}.

H. Kurihara \cite{K} proved the following lemma:

\begin{lemma}[Kurinara \cite{K}]\label{lem-sn-3}
Let $X_0$ be a finite graph, and let $X^0$ be the vectors of the crystal lattice which is the maximal Abelian covering and constructed with the method in Section \ref{sec-sr-lattice}, then we define $X$ the set of normalized vectors of $X^0$. $X$ is a spherical $3$-design if and only if the vectors of $X^0$ have the same norm.
\end{lemma}

For the proof of the above lemma, we use the conditions $({C0}')$, \ldots, $({C3}')$ in the Remark \ref{rem-sn-3} with $w_{i, k} = c_{i, k} / \sqrt{2}$. Then, we have the above conditions $(C0)$, \ldots, $(C3)$ for spherical $2$-design. In addition, we consider the fact that $X$ is antipodal.

Furhtermore, Kurihara proved the following theorems:

\begin{theorem}[Kurinara \cite{K}]\label{th-et}
Let $X_0$ be a finite graph, and let $X^0$ be the vectors constructed from $X_0$ with the method in Section \ref{sec-sr-lattice}. If $X_0$ is edge-transitive, then the vectors of $X^0$ have the same norm.
\end{theorem}

\begin{theorem}[Kurinara \cite{K}]\label{th-am3}
Let $X_0$ be a finite graph, and let $X^0$ be the vectors constructed with the method in Section \ref{sec-sr-lattice}, then we define $X$ the set of normalized vectors of $X^0$. Then, $X$ is spherical $t$-design for $t \leqslant 3$, except for the case of the Hexagonal lattice $($See Section \ref{sec-sr-lattice} and Figure \ref{fig-2dia}$)$.
\end{theorem}

Here, we consider strongly regular graphs. Let $\Gamma_x$ be a set of vertices of the graph $X_0$ which are adjacent to $x$. $k$-regular graph with $v$ vertices is a strongly regular graph with parameters $(v, k, \lambda, \mu)$ if $i)$ $|\Gamma_x \cap \Gamma_y|$ is equal to the constant $\lambda$ for every pair of vertices $(x, y)$ which are adjacent to each other and $ii)$ $|\Gamma_x \cap \Gamma_y|$ is equal to the constant $\mu$ for every pair of vertices $(x, y)$ which are not adjacent. E. Bannai \cite{B} proposed the following conjecture:

\begin{conjecture}[Bannai \cite{B}]\label{conj-srg}
Let $X_0$ be a finite graph, and let $X^0$ be the vectors constructed from $X_0$ with the method in Section \ref{sec-sr-lattice}. If $X_0$ is strongly regular, then the vectors of $X^0$ have the same norm.
\end{conjecture}
By Lemma \ref{lem-sn-3}, the normalized vectors $X$ is also a spherical $3$-design. In Section \ref{sec-sr}, we calculate every strongly regular graphs with at most $30$ vertices and many other strongly regular graphs, where the vectors from these graphs verify the above conjecture, that is, the vectors of every constructed crystal lattice have the same norm. In addition, the sets of vectors from these graphs are at most $10$-distance sets.\\

Furthermore, there are some strongly regular graphs which are not edge-transitive (See Section \ref{sec-dif-sr}). Then, we want to consider an sufficient condition of the graphs from which we construct vectors of crystal lattices of the same norm. Here, we have the following fact:

\begin{proposition}[Dauber and Harary \cite{DH}]
If a connected edge-transitive graph is not vertex-transitive, then it is bipartite.
\end{proposition}

Now, we want to suggest the graphs from association schemes as a candidate for such graphs, because they include both vertex and edge-transitive graphs and strongly regular graphs.

Let $X$ be a finite set, and we call the subset of $X$ {\it relation on $X$}. We denote $s^{*} := \{ (y, x) \in X \times X \, ; \, (x, y) \in s \}$ for $s$ the relation on $X$. Furthermore, we define $x s := \{ y \in X \, ; \, (x, y) \in s \}$ and $s x := \{ y \in X \, ; \, (y, x) \in s \}$ for every $x \in X$ and every relation $s$. Now, we have the following definition:

\begin{definition}[Association scheme \cite{BI}, \cite{Ha}]\label{def-as}
Let $X$ be a finite set, and $s_0$, $s_1$, \ldots, $s_d$ be non-empty relations on $X$. We call $(X, S := \{ s_i \}_{0 \leqslant i \leqslant d})$ a {\it association scheme of class $d$} if it satisfies the following conditions:
\def\labelenumi{(\arabic{enumi})}
\begin{enumerate}
\item $S$ is a partition of $X \times X$, that is $X \times X = \bigcup_{0 \leqslant i \leqslant d}$ and $s_i \cap s_j = \phi$ when $i \ne j$.
\item $1 := \{ (x, x) \, ; \, x \in X \} \in S$. (We define $s_0 := 1$ without loss of generality.)
\item If $s \in S$, then $s^{*} \in S$.
\item For every $s_i, s_j, s_k \in S$, $p_{i, j}^{k} := | x s_i \cap s_j y |$ is constant if $(x, y) \in s_k$.
\end{enumerate}
In addition, if an association scheme $(X, S)$ satisfies the following condition
\begin{enumerate}
\item[(5)] $p_{i, j}^{k} = p_{j, i}^{k}$ for every $i, j, k \in \{ 0, \ldots , d \}$,
\end{enumerate}
then we call $(X, S)$ is commutative. Also, we call an association scheme $(X, S)$ symmetric if it satisfies the following condition:
\begin{enumerate}
\item[(6)] $s^{*} = s$ for every relation $s \in S$.
\end{enumerate}
\end{definition}

We define the {\it adjacency matrix $A_i$} which corresponds to $s_i$ the relation on $X$, where every $(x, y)$ component of $A_i$ is defined by $A_i(x, y) = 1$ and $0$ according as $(x, y) \in s$ and otherwise, respectively. If $A_i$ is symmetric, then we regard $A_i$ as an adjacency matrix of a graph.

Note that strongly regular graphs corresponds to symmetric association schemes of class $2$. On the other hand, for vertex and edge-transitive graphs, we have the following definition:

\begin{definition}[Schurian scheme \cite{BI} (Example 2.1)]\label{def-ss}
Let $X$ be a finite set, and $G$ be a transitive permutation group on $X$. We can consider the action of $G$ on $X \times X$, then let $s_0$, \ldots, $s_d$ be the orbits of $G$ on $X \times X$. We call $(X, S := \{ s_i \}_{0 \leqslant i \leqslant d})$ a {\it Schurian scheme}.
\end{definition}
Note that every Schurian scheme $(X, \{ s_i \}_{0 \leqslant i \leqslant d})$ is an association scheme.

We also have the fact that vertex and edge-transitive graphs have correspondences to Schurian schemes, in exact word, we have the following fact:

\begin{proposition}\label{prop-ss-et}
Every graph from symmetric adjacency matrix of every association scheme is vertex and edge-transitive. On the other hand, every vertex and edge-transitive graph has an adjacency matrix which is also an adjacency matrix of an association scheme.
\end{proposition}

We can prove the above proposition easily, in the proof, we can consider the correspondence between $G$ the transitive permutation group on $X$ and the automorphism group on the graph.

Now, we have the following conjecture:

\begin{conjecture}\label{conj-as}
Let $X_0$ be a finite graph, and let $X^0$ be the vectors constructed from $X_0$ with the method in Section \ref{sec-sr-lattice}. If $X_0$ is a graph from an adjacency matrix of an association scheme, then the vectors of $X^0$ have the same norm.
\end{conjecture}
In Section \ref{sec-oth-as}, we calculate every graph from adjacency matrix of non-Schurian association scheme with at most $30$ vertices, where the vectors from these graph verify the above conjecture. For Schurian schemes, we can prove the above conjecture by Theorem \ref{th-et} and Proposition \ref{prop-ss-et}. In conclusion, we can prove the above conjecture if $X_0$ has at most $30$ vertices.\\

On the other hand, it is difficult to consider necessary condition of the graphs from which we construct vectors of crystal lattices of the same norm. There are some such graphs which are not regular and not edge-transitive (See Section \ref{sec-net}).\\

\section{Standard realization of crystal lattices}\label{sec-sr-lattice}
In this section, we introduce the method of {\it standard realization} of crystal lattices as Abelian covering of finite graphs. For the details, we refer \cite{St1} and \cite{St2}.

Let $X_0 = (V_0, E_0)$ be a finite graph with $V_0$ the set of vertices and $E_0$ be the set of edges. Here, we mean by a graph loops and multiple edges are allowed and edges have orientations. Furthermore, $E_0$ contains edges $e_0$ and $\overline{e_0}$ simultaneously, where $\overline{e_0} = - e_0$ denote the edge with the opposite direction as $e_0$.

We denote by $C_0(X_0, \mathbb{R})$ the group of $0$-chain on $X_0$, that is the set of $\mathbb{R}$-linear combination of $V_0$, and by $C_1(X_0, \mathbb{R})$ the group of $1$-chain on $X_0$, that is the set of $\mathbb{R}$-linear combination of $E_0$. Consider the boundary operator $\partial$ which is the linear map from $C_1(X_0, \mathbb{R})$ to $C_0(X_0, \mathbb{R})$ defined by $\partial (e_0) = t(e_0) - o(e_0)$ for $e_0 \in E_0$, where $t(e_0)$ and $o(e_0)$ are the terminal and the origin of every edge $e_0$. We define $H_1(X_0, \mathbb{R})$ the homology group of degree $1$ by the kernel of $\partial$.

Now, we assume that there exists some surjective homomorphism $\mu$ from $H_1(X_0, \mathbb{Z})$ to $L$, where $L$ is an Abelian group which is isomorphic to $\mathbb{Z}^d$. Let $W$ be the subspace of $H_1(X_0, \mathbb{R})$ which is spanned by $\text{Ker} \, \mu$, and let $H$ be the intersection of the orthogonal complement $W^{\perp}$ and $H_1(X_0, \mathbb{R})$. Note that $H = H_1(X_0, \mathbb{R})$ if $L = H_1(X_0, \mathbb{Z})$.

Let $P : C_1(X_0, \mathbb{R}) \rightarrow H (\subset C_1(X_0, \mathbb{R}))$ be the orthogonal projection with respect to the inner product on $C_1(X_0, \mathbb{R})$ defined by $\langle e_1, e_2 \rangle = 1, - 1, 0$ according as $e_1 = e_2, e_1 = \overline{e_2}$, otherwise (respectively). The set $\{P(e_0) \, ; \, e_0 \in E_0\}$ is a set of vectors in $\mathbb{R}^d(= H_1(X_0, \mathbb{\mathbb{R}}))$. Then, we get a realization of $d$-dimensional crystal lattice as the Abelian covering of the graph $X_0$. Furthermore, this realization is a {\it standard realization} (See \cite{St2}).

In particular, if $L = H_1(X_0, \mathbb{Z})$, then the crystal lattice is the maximal Abelian covering of $X_0$. In this paper, we consider realization of the crystal lattices as the maximal Abelian covering of the graph. Then, we have $d = \dim H_1(X_0, \mathbb{R}) = 1 - |V_0| + |E_0| / 2$, if the graph $X_0$ is connected. Recall that $E_0$ contains edges $e_0$ and $\overline{e_0}$ simultaneously, thus $|E_0| / 2$ means the number of edges up to orientation.\\

\begin{figure}[htbp]
\begin{center}
{Graph}\includegraphics[height=1.3in]{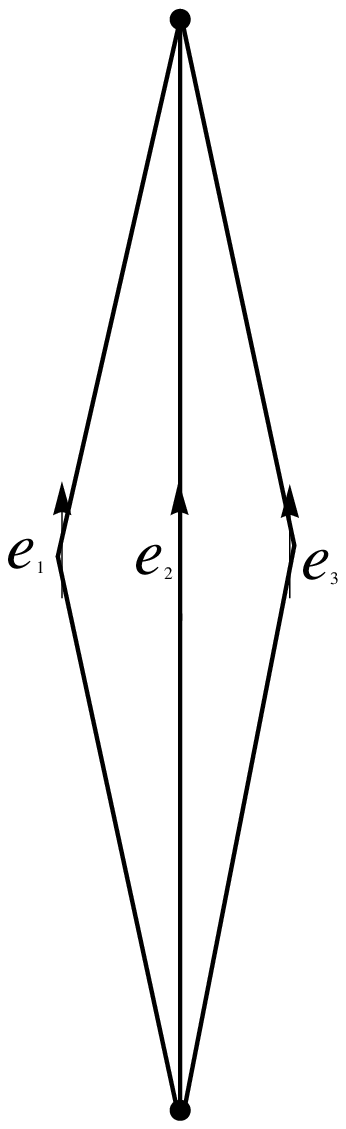}
\quad $\Rightarrow$ \quad
{Lattice}\includegraphics[height=1.3in]{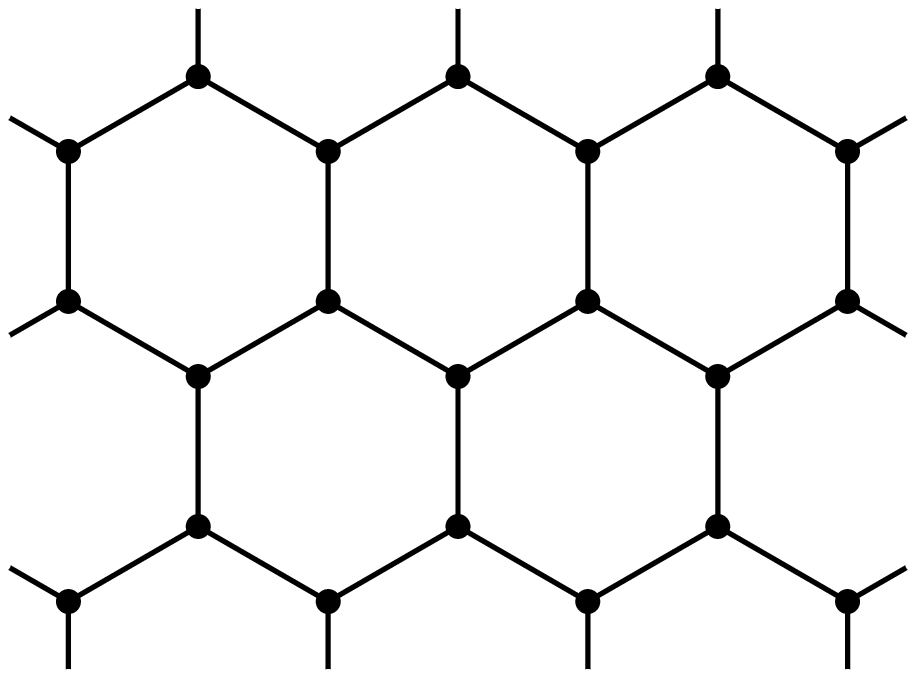}
\quad $\Rightarrow$ \quad
{Vectors}\includegraphics[height=1.3in]{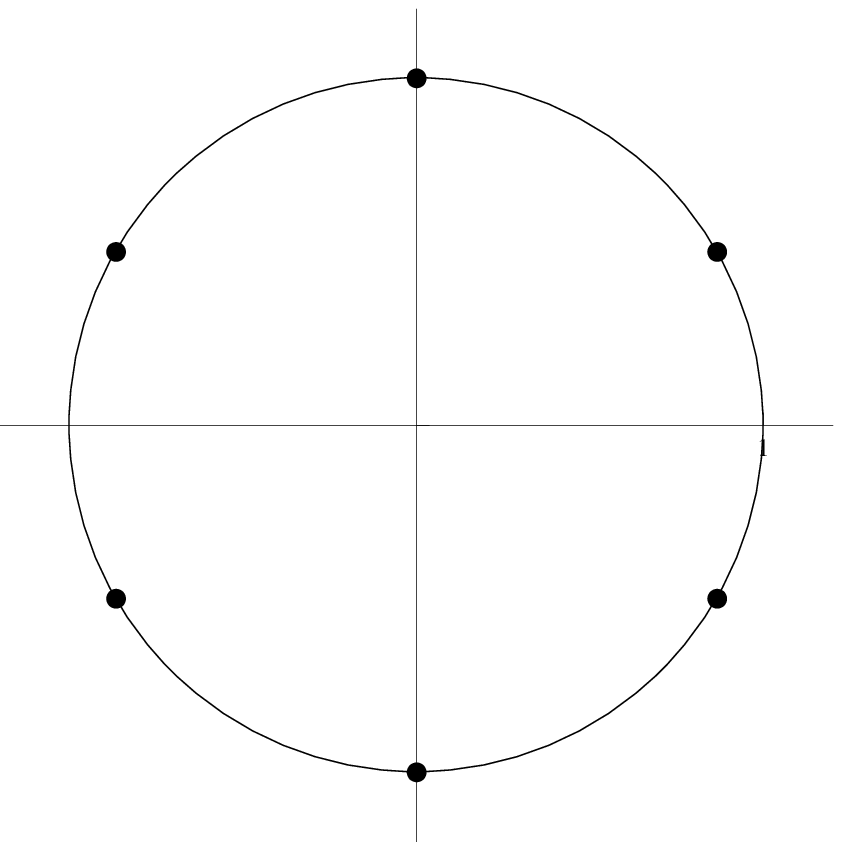}
\end{center}
\caption{Hexagonal lattice}\label{fig-2dia}
\end{figure}

We can see the concrete example, the construction of $K_4$ crystal from the complete graph $K_4$, in \cite{St1} and \cite{St2}. In this paper, we consider another example, the {\it Hexagonal lattice} from the graph in Figure \ref{fig-2dia}. Let $X_0 = (V_0, E_0)$ be the graph in Figure \ref{fig-2dia}, we have two vertices in $V_0$ and we have $E_0 = \{ \pm e_1, \pm e_2, \pm e_3 \}$, where the orientations of $e_1$, $e_2$, and $e_3$ are fixed as in Figure \ref{fig-2dia}.

Now, we have two closed path $c_1 = (e_1, - e_2)$ and $c_2 = (e_2, - e_3)$ as the basis of $H_1(X_0, \mathbb{Z})$, which satisfy $\| c_1 \|^2 = \| c_2 \|^2 = 2$ and $\langle c_1, c_2 \rangle = -1$. Since we have $\langle P(e_1), c_1 \rangle = \langle e_1, c_1 \rangle = 1$ and $\langle P(e_1), c_2 \rangle = \langle e_1, c_2 \rangle = 0$, we have $P(e_1) = \frac{2}{3} c_1 + \frac{1}{3} c_2$. Similarly, we have $P(e_2) = - \frac{1}{3} c_1 + \frac{1}{3} c_2$ and $P(e_3) = - \frac{1}{3} c_1 - \frac{2}{3} c_2$. Furthermore, to regard $P(e_i)$ as the vectors in $\mathbb{R}^2$, we consider orthonormalization of $\{ c_1, c_2 \}$. Then, we have $\{ {c_1}', {c_2}' \}$ a orthonormal basis of $H_1(X_0, \mathbb{Z})$, where ${c_1}' = c_1 / \sqrt{2}$ and ${c_2}' = (c_1 + 2 c_2) / \sqrt{6}$. And then, we have $P(e_1) = \frac{1}{\sqrt{2}} {c_1}' + \frac{1}{\sqrt{6}} {c_2}'$, $P(e_2) = - \frac{1}{\sqrt{2}} {c_1}' + \frac{1}{\sqrt{6}} {c_2}'$, and $P(e_3) = - \sqrt{\frac{2}{3}} {c_2}'$. Here, these vectors have the same norm $\| P(e_1) \| = \| P(e_2) \| = \| P(e_3) \| = \frac{2}{3}$. Finally, after normalization, we have
the six vectors
\begin{equation}
\left\{ \pm \left( \frac{\sqrt{3}}{2}, \ \frac{1}{2} \right), \quad \pm \left( \frac{\sqrt{3}}{2}, \ - \frac{1}{2} \right), \quad \pm \left( 0, \ 1 \right) \right\}.
\end{equation}
Here, these vectors form a regular hexagon (Figure \ref{fig-2dia}). Furthermore, note that set of these vectors is a spherical $5$-design.

In another word, since $\{ e_1, e_2, e_3 \}$ are orthonormal basis of $C_1(X_0, \mathbb{R})$, we denote $e_1 = (1, 0, 0)$, $e_2 = (0, 1, 0)$, and $e_3 = (0, 0, 1)$. Then we can write $c_1 = (1, -1, 0)$ and $c_2 = (0, 1, -1)$, and then we have ${c_1}' = \frac{1}{\sqrt{2}} (1, -1, 0)$ and ${c_2}' = \frac{1}{\sqrt{6}} (1, 1, -2)$. Here, we denote ${c_i}' = (c_{i, 1}, c_{i, 2}, c_{i, 3})$ for $i = 1, 2$. Then we have $\langle P(e_k), {c_i}' \rangle = \langle e_k, {c_i}' \rangle = c_{i, k}$ for $1 \leqslant i \leqslant 3$ and $k = 1, 2$, and then we can write $P(e_k) = (c_{1, k}, c_{2, k})$ for $1 \leqslant k \leqslant 3$ since $\{ {c_1}', {c_2}' \}$ are orthonormal basis.

\begin{remark}\label{rem-sn-3}
In general, let $d := H_1(X_0, \mathbb{Z})$ and $e := |E_0| / 2$, and let $\{ e_1, e_2, \ldots, e_{2 e} \} \subset E_0$ be basis of $C_1(X_0, \mathbb{R})$. Then, we can write $e_1 = (1, 0, \ldots, 0)$, $e_2 = (0, 1, \ldots, 0)$, and $e_e = (0, 0, \ldots, 1)$. Here, let $\{ {c_1}', {c_2}', \ldots, {c_d}' \}$ be orthonormal basis of $H_1(X_0, \mathbb{Z})$, then we can write ${c_k}' = (c_{i, 1}, c_{i, 2}, \ldots, c_{i, e})$ for $1 \leqslant k \leqslant d$. Furthermore, we have $P(e_k) = (c_{1, k}, c_{2, k}, \ldots, c_{d, k})$ for $1 \leqslant k \leqslant e$, and denote $P(e_{k + e}) := P(-e_k) = -(c_{1, k}, c_{2, k}, \ldots, c_{d, k})$. Thus, we have the following relations:
\begin{align}
\tag{${C1}'$} &\sum_{k=1}^{2 e} c_{i, k} = 0, \qquad 1 \leqslant i \leqslant d.\\
\tag{${C2}'$} &\sum_{k=1}^{2 e} c_{i, k} c_{j, k} = 0, \qquad 1 \leqslant i < j \leqslant d.\\
\tag{${C3}'$} &\sum_{k=1}^{2 e} c_{i, k}^2 = 2, \qquad 1 \leqslant i \leqslant d.
\end{align}
Furthermore, if every $P(e_0)$ has the same norm for $e_0 \in E_0$, then we have
\begin{equation}
\tag{${C0}'$} \sum_{i=1}^d c_{i, k}^2 = \frac{d}{e}, \qquad 1 \leqslant k \leqslant 2 e.
\end{equation}
\end{remark}\quad

\section{History}

This problems were researched in the seminar for graduate students in master course of Laboratory of Professor Eiichi Bannai in Kyushu University, and it was first suggested by E. Bannai. Then, H. Kurihara \cite{K} considered the vectors of crystal lattices from the graph from the regular polyhedrons (See Section \ref{sec-3_reg_poly}) and some other graphs at once, and then we knew that we can construct some good spherical designs from some graphs.

After that, by the result of Kurihara, Bannai suggested the fact of Theorem \ref{th-et}. Furthermore, M. Shinohara \cite{Sm} suggested the complete bipartite graphs, which are edge-transitive but not regular. Then, I confirm that the vectors of crystal lattices from many complete multipartite graphs satisfy Theorem \ref{th-et} (See Section \ref{sec-bip} and \ref{sec-ukn_cm}) And then, Kurihara succeeded to prove Lemma \ref {lem-sn-3} and Theorem \ref{th-et}. Since then, we started to search for the graphs of not edge-transitive, from which we can construct vectors of crystal lattices with the same norms.

D. Fujino \cite{FIM} confirmed that the vectors from the complement of Shrikhande graph have the same norm, where it is not edge-transitive (See Section \ref{sec-dif-l24}). Then, K. Mitsumaru \cite{FIM} considered the vectors from the Chang graphs, and showed that they have the same norm and consider the difference among the vectors of crystal lattices from them (See Section \ref{sec-dif-t8}). And then, Fujino consider the every strongly regular graph with at most $30$ vertices, and confirm that they satisfy Conjecture \ref{conj-srg}. Furthermore, T. Ishii \cite{FIM} found the graph $netg(6, 10)$ in Figure \ref{fig-net}, where it is not regular and not edge-transitive, but the vectors of crystal lattice from it have the same norm. Finally, Kurihara proved Theorem \ref{th-am3}.

I had supported them in numerical calculations and I also did independently. Now, in the paper, I supplement the result by numerical calculation from the above graphs and also calculate from the graphs from association schemes. \newpage

\section{Notations}
Let $X$ be a non-empty finite set on the Euclidean sphere $S^{d-1}$ ($\subset \mathbb{R}^{d}$). We denote the distance set of $X$ by $A(X) := \{ (x, y) ; x, y \in X, x \ne y \}$ where $(x, y)$ denotes the standard inner product of the vectors $x, y \in \mathbb{R}^{d}$, then we call $X$ a $s$-distance set if $| A(X) | = s$. Furthermore, $X$ is said to be a $(d, n, s, t)$-configuration if $X \subset S^{d-1}$ is of order $n (:= | X |)$, a $s$-distance set, and a spherical $t$-design. We also denote $A_x(X, a) := \{ y \in X  ; (x, y) = a \}$ for $x \in X$ and $a \in A(X)$.

In this research, we start from a finite graph with $v$ vertices and $e$ edges. Then, we consider standard realization of the crystal lattices as the maximal Abelian covering of the graph. And then, we take the set of vectors which form the crystal lattice. If every vector has the same norm, then we can consider a finite set on Euclidean sphere, and then we can consider spherical design.\\

\section{Trivial examples}

\subsection{$d$-dimensional standard crystal lattice}\label{sec-sta}
For the graph with one vertex and $d$ loops, we have the $d$-dimensional standard crystal lattice as the maximal Abelian covering of the graph. For the $(d, n, s, t)$-configurations of the vectors of the crystal lattices, we have the following table:

\begin{table}[h]
\begin{center}
\begin{tabular}{cccc}
\hline
$d$ & $n$ & $s$ & $t$\\
\hline
$2$ & $4$ & $2$ & $3$\\
$3$ & $6$ & $2$ & $3$\\
$4$ & $8$ & $2$ & $3$\\
\hline
\end{tabular}
\qquad
\begin{tabular}{cccc}
\hline
$d$ & $n$ & $s$ & $t$\\
\hline
$5$ & $10$ & $2$ & $3$\\
$6$ & $12$ & $2$ & $3$\\
$7$ & $14$ & $2$ & $3$\\
\hline
\end{tabular}
\end{center}\quad
\caption{Standard crystal lattice}\vspace{-0.2in}
\end{table}

\noindent
For every integer $d \geqslant 2$, we can easily prove that the configuration of the vectors is $(d, n, s, t) = (d, 2 d, 2, 3)$.\\

\quad\\

\subsection{$d$-dimensional diamond crystal lattice}\label{sec-dia}
For the graph with two vertices and $d+1$ edges ($(d+1)$-fold multiple edge) connecting the two vertices, we have the $d$-dimensional diamond crystal lattice as the maximal Abelian covering of the graph. For the $(d, n, s, t)$-configurations of the vectors of the crystal lattices, we have the following table:

\begin{table}[h]
\begin{center}
\begin{tabular}{cccc}
\hline
$d$ & $n$ & $s$ & $t$\\
\hline
$2$ & $6$ & $3$ & \underline{$5$}\\
$3$ & $8$ & $3$ & $3$\\
$4$ & $10$ & $3$ & $3$\\
$5$ & $12$ & $3$ & $3$\\
\hline
\end{tabular}
\qquad
\begin{tabular}{cccc}
\hline
$d$ & $n$ & $s$ & $t$\\
\hline
$6$ & $14$ & $3$ & $3$\\
$7$ & $16$ & $3$ & $3$\\
$8$ & $18$ & $3$ & $3$\\
$9$ & $20$ & $3$ & $3$\\
\hline
\end{tabular}
\end{center}\quad
\caption{Diamond crystal lattice}\vspace{-0.2in}
\end{table}

\noindent
For every integer $d \geqslant 2$, we can easily prove that we have the distance set $A(X) = \{ -1, \: \pm 1/d \}$ and that the configuration of the vectors is $(d, n, s, t) = (d, 2 (d+1), 3, 3)$ (If $d=2$, then $t=5$). Furthermore, we have $X = \mathcal{S}_d \cup (-\mathcal{S}_d)$, where $\mathcal{S}_d$ is the $d$-dimensional simplex.

\newpage

\section{Graphs from regular polytopes}\label{sec-reg_poly}
In this section, we consider the graphs from regular and semi-regular polyhedrons and regular polytopes, where we can consider edges as the pair of vertices between which we take minimum distance. Note that the graphs from regular polyhedrons and regular polytopes are vertex and edge-transitive.\\

\subsection{Regular polyhedrons (in $\mathbb{R}^3$)}\label{sec-3_reg_poly}
For the $(d, n, s, t)$-configurations of the regular polyhedrons and those of the vectors of the constructed crystal lattices, we have the following table:

\begin{table}[h]
\begin{center}
\begin{tabular}{ccccccc}
\multicolumn{6}{c}{Regular polyhedron}\\
\hline
& $d$ & $v$ & $e$ & $s$ & $t$\\
\hline
Tetrahedron & $3$ & $4$ & $6$ & $1$ & $2$\\
Hexahedron & $3$ & $8$ & $12$ & $3$ & $3$\\
Octahedron & $3$ & $6$ & $12$ & $2$ & $3$\\
Dodecahedron & $3$ & $20$ & $30$ & $5$ & $5$\\
Icosahedron & $3$ & $12$ & $30$ & $3$ & $5$\\
\hline
\end{tabular}
\qquad $\Rightarrow$ \qquad
\begin{tabular}{cccc}
\multicolumn{4}{c}{Vectors of lattice}\\
\hline
$d$ & $n$ & $s$ & $t$\\
\hline
$3$ & $12$ & $4$ & $3$\\
$5$ & $24$ & $9$ & $3$\\
$7$ & $24$ & $9$ & $3$\\
$11$ & $60$ & $12$ & $3$\\
$19$ & $60$ & $12$ & $3$\\
\hline
\end{tabular}
\end{center}\quad
\caption{Regular polyhedron}\vspace{-0.2in}
\end{table}

\noindent
The vectors of the crystal lattice from tetrahedron is cuboctahedron, which is a semi-regular polyhedron.\\

\quad\\

\subsection{Semi-regular polyhedrons (in $\mathbb{R}^3$)}\label{sec-3_semireg_poly}
The graphs from semi-regular polyhedrons are vertex-transitive, and only the graphs from cuboctahedron and icosidodecahedron are edge-transitive. For the $(d, n, s, t)$-configurations of the semi-regular polyhedrons (Archimedean solids) and those of the vectors of the constructed crystal lattices, we have the following table:

\begin{table}[h]
\begin{center}
\begin{tabular}{ccccccccccc}
\multicolumn{6}{c}{Semi-regular polyhedron} & & \multicolumn{4}{c}{Vectors of lattice}\\
\cline{1-6} \cline{8-11}
& $d$ & $v$ & $e$ & $s$ & $t$ && $d$ & $n$ & $s$ & $t$\\
\cline{1-6} \cline{8-11}
Truncated tetrahedron & $3$ & $12$ & $18$ & $4$ & $2$ && $7$ & $36$ & $20 (2)$\\
Cuboctahedron & $3$ & $12$ & $24$ & $4$ & $3$ && $13$ & $48$ & $16$ & $3$\\
Truncated octahedron & $3$ & $24$ & $36$ & $10$ & $3$ && $13$ & $72$ & $43(2)$\\
Truncated hexahedron & $3$ & $24$ & $36$ & $7$ & $3$ && $13$ & $72$ & $32(2)$\\
Rhombicuboctahedron & $3$ & $24$ & $48$ & $7$ & $3$ && $25$ & $96$ & $65(2)$\\
Great rhombicuboctahedron & $3$ & $48$ & $72$ & $25$ & $3$ && $25$ & $144$ & $145(3)$\\
Icosidodecahedron & $3$ & $30$ & $60$ & $8$ & $5$ & \qquad {\Large $\Rightarrow$} \quad \ & $31$ & $120$ & $33$ & $3$\\
Truncated icosahedron & $3$ & $60$ & $90$ & $21$ & $5$ && $31$ & $180$ & $88(2)$\\
Truncated dodecahedron & $3$ & $60$ & $90$ & $19$ & $5$ && $31$ & $180$ & $66(2)$\\
Snub cube & $3$ & $24$ & $60$ & $11$ & $3$ && $37$ & $120$ & $171(3)$\\
Rhombicosidodecahedron & $3$ & $60$ & $120$ & $17$ & $5$ && $61$ & $240$ & $151(2)$\\
Great rhombicosidodecahedron & $3$ & $120$ & $180$ & $61$ & $5$ && $61$ & $360$ & $331(3)$\\
Snub dodecahedron & $3$ & $60$ & $180$ & $31$ & $5$ && $91$ & $300$ & $411(3)$\\
\cline{1-6} \cline{8-11}
\end{tabular}
\end{center}\quad
\caption{Semi-regular polyhedron}\vspace{-0.2in}
\end{table}

\noindent
Here, we denote the number in `$( \ )$' by the number of the norms of the vectors of lattices.

\newpage

\subsection{$4$-dimensional regular polytopes}\label{sec-4_reg_poly}
For the $(d, n, s, t)$-configurations of the $4$-dimensional regular polytopes (regular polychorons) and those of the vectors of the constructed crystal lattices, we have the following table:

\begin{table}[h]
\begin{center}
\begin{tabular}{cccccc}
\multicolumn{6}{c}{$4$-dimensional regular polytope}\\
\hline
& $d$ & $v$ & $e$ & $s$ & $t$\\
\hline
$5$-cell polytope & $4$ & $5$ & $10$ & $1$ & $2$\\
$8$-cell polytope & $4$ & $16$ & $32$ & $3$ & $3$\\
$16$-cell polytope & $4$ & $8$ & $24$ & $2$ & $3$\\
$24$-cell polytope & $4$ & $24$ & $96$ & $4$ & $5$\\
$120$-cell polytope & $4$ & $600$ & $1200$ & $30$ & $11$\\
$600$-cell polytope & $4$ & $120$ & $720$ & $8$ & $11$\\
\hline
\end{tabular}
\quad $\Rightarrow$ \quad
\begin{tabular}{cccc}
\multicolumn{4}{c}{Vectors of lattice}\\
\hline
$d$ & $n$ & $s$ & $t$\\
\hline
$6$ & $20$ & $4$ & $3$\\
$17$ & $64$ & $11$ & $3$\\
$17$ & $48$ & $10$ & $3$\\
$73$ & $192$ & $24$ & $3$\\
$601$ & $2400$ & $154$ & $3$\\
$601$ & $1440$ & $64$ & $3$\\
\hline
\end{tabular}
\end{center}\quad
\caption{$4$-dimensional regular polytope}
\end{table}\quad

\subsection{$k$-dimensional regular $(k+1)$-cell polytopes (Complete graph $K_{k+1}$)}\label{sec-kn}
The graph from $k$-dimensional regular $(k+1)$-cell polytope is the complete graph $K_{k+1}$. For the $(d, n, s, t)$-configurations of the $k$-dimensional regular $(k+1)$-cell polytopes and those of the vectors of the constructed crystal lattices, we have the following table:

\begin{table}[h]
\begin{center}
\begin{tabular}{ccccc}
\multicolumn{5}{c}{\small $(k+1)$-cell polytope}\\
\hline
$k$ & $v$ & $e$ & $s$ & $t$\\
\hline
$3$ & $4$ & $6$ & $1$ & $2$\\
$4$ & $5$ & $10$ & $1$ & $2$\\
$5$ & $6$ & $15$ & $1$ & $2$\\
$6$ & $7$ & $21$ & $1$ & $2$\\
$7$ & $8$ & $28$ & $1$ & $2$\\
$8$ & $9$ & $36$ & $1$ & $2$\\
$9$ & $10$ & $45$ & $1$ & $2$\\
$10$ & $11$ & $55$ & $1$ & $2$\\
$11$ & $12$ & $66$ & $1$ & $2$\\
$12$ & $13$ & $78$ & $1$ & $2$\\
$13$ & $14$ & $91$ & $1$ & $2$\\
$14$ & $15$ & $105$ & $1$ & $2$\\
$15$ & $16$ & $120$ & $1$ & $2$\\
$16$ & $17$ & $136$ & $1$ & $2$\\
\hline
\end{tabular}
\ $\Rightarrow$ \
\begin{tabular}{cccc}
\multicolumn{4}{c}{\small Vectors of lattice}\\
\hline
$d$ & $n$ & $s$ & $t$\\
\hline
$3$ & $12$ & $4$ & $3$\\
$6$ & $20$ & $4$ & $3$\\
$10$ & $30$ & $4$ & $3$\\
$15$ & $42$ & $4$ & $3$\\
$21$ & $56$ & $4$ & $3$\\
$28$ & $72$ & $4$ & $3$\\
$36$ & $90$ & $4$ & $3$\\
$45$ & $110$ & $4$ & $3$\\
$55$ & $132$ & $4$ & $3$\\
$66$ & $156$ & $4$ & $3$\\
$78$ & $182$ & $4$ & $3$\\
$91$ & $310$ & $4$ & $3$\\
$105$ & $240$ & $4$ & $3$\\
$120$ & $272$ & $4$ & $3$\\
\hline
\end{tabular}
\qquad \quad
\begin{tabular}{ccccc}
\multicolumn{5}{c}{\small $(k+1)$-cell polytope}\\
\hline
$k$ & $v$ & $e$ & $s$ & $t$\\
\hline
$17$ & $18$ & $153$ & $1$ & $2$\\
$18$ & $19$ & $171$ & $1$ & $2$\\
$19$ & $20$ & $190$ & $1$ & $2$\\
$20$ & $21$ & $210$ & $1$ & $2$\\
$21$ & $22$ & $231$ & $1$ & $2$\\
$22$ & $23$ & $253$ & $1$ & $2$\\
$23$ & $24$ & $276$ & $1$ & $2$\\
$24$ & $25$ & $300$ & $1$ & $2$\\
$25$ & $26$ & $325$ & $1$ & $2$\\
$26$ & $27$ & $351$ & $1$ & $2$\\
$27$ & $28$ & $378$ & $1$ & $2$\\
$28$ & $29$ & $406$ & $1$ & $2$\\
$29$ & $30$ & $435$ & $1$ & $2$\\
$30$ & $31$ & $465$ & $1$ & $2$\\
\hline
\end{tabular}
\ $\Rightarrow$ \
\begin{tabular}{cccc}
\multicolumn{4}{c}{\small Vectors of lattice}\\
\hline
$d$ & $n$ & $s$ & $t$\\
\hline
$136$ & $306$ & $4$ & $3$\\
$153$ & $342$ & $4$ & $3$\\
$171$ & $380$ & $4$ & $3$\\
$190$ & $420$ & $4$ & $3$\\
$210$ & $462$ & $4$ & $3$\\
$231$ & $506$ & $4$ & $3$\\
$253$ & $552$ & $4$ & $3$\\
$276$ & $600$ & $4$ & $3$\\
$300$ & $650$ & $4$ & $3$\\
$325$ & $702$ & $4$ & $3$\\
$351$ & $756$ & $4$ & $3$\\
$378$ & $812$ & $4$ & $3$\\
$406$ & $870$ & $4$ & $3$\\
$435$ & $930$ & $4$ & $3$\\
\hline
\end{tabular}
\end{center}\quad
\caption{$k$-dimensional regular $(k+1)$-cell polytope}\vspace{-0.2in}
\end{table}

\noindent
For every integer $k \geqslant 3$, we can easily prove that we have the distance set $A(X) = \{ -1, \: \pm 1/(k-1), \: 0 \}$ and that the configuration of the vectors is $(d, n, s, t) = (k(k-1)/2, k(k+1), 4, 3)$. Furthermore, for each $x \in X$, we have the following table:
\begin{center}
\begin{tabular}{c|ccc}
$a$ & $-1$ & $\pm 1/(k-1)$ & $0$\\
\hline
$|A_x(X, a)|$ & $1$ & $2 (k-1)$ & $(k-1) (k-2)$\\
\end{tabular}
\end{center}

Incidentally, the configuration of the $k$-dimensional regular $(k+1)$-cell polytope is $(d, n, s, t) = (k, k+1, 1, 2)$.

\newpage

\subsection{$k$-dimensional regular $2 k$-cell polytopes (Hamming graph $H(k, 2)$)}\label{sec-cube}
The graph from $k$-dimensional regular $2 k$-cell polytope is the Hamming graph $H(k, 2)$. For the $(d, n, s, t)$-configurations of the $k$-dimensional regular $2 k$-cell polytopes and those of the vectors of the constructed crystal lattices, we have the following table:

\begin{table}[h]
\begin{center}
\begin{tabular}{ccccc}
\multicolumn{5}{c}{$2 k$-cell polytope}\\
\hline
$k$ & $v$ & $e$ & $s$ & $t$\\
\hline
$3$ & $8$ & $12$ & $3$ & $3$\\
$4$ & $16$ & $32$ & $4$ & $3$\\
$5$ & $32$ & $80$ & $5$ & $3$\\
\hline
\end{tabular}
 $\Rightarrow$ 
\begin{tabular}{cccc}
\multicolumn{4}{c}{\small Vectors of lattice}\\
\hline
$d$ & $n$ & $s$ & $t$\\
\hline
$5$ & $24$ & $9$ & $3$\\
$17$ & $64$ & $11$ & $3$\\
$49$ & $160$ & $17$ & $3$\\
\hline
\end{tabular}
\qquad
\begin{tabular}{ccccc}
\multicolumn{5}{c}{$2 k$-cell polytope}\\
\hline
$k$ & $v$ & $e$ & $s$ & $t$\\
\hline
$6$ & $64$ & $192$ & $6$ & $3$\\
$7$ & $128$ & $448$ & $7$ & $3$\\
$8$ & $256$ & $1024$ & $8$ & $3$\\
\hline
\end{tabular}
 $\Rightarrow$ 
\begin{tabular}{cccc}
\multicolumn{4}{c}{\small Vectors of lattice}\\
\hline
$d$ & $n$ & $s$ & $t$\\
\hline
$129$ & $384$ & $19$ & $3$\\
$321$ & $896$ & $25$ & $3$\\
$769$ & $2048$ & $27$ & $3$\\
\hline
\end{tabular}
\end{center}\quad
\caption{$k$-dimensional regular $2 k$-cell polytope}\vspace{-0.2in}
\end{table}

\noindent
For every integer $3 \leqslant k \leqslant 8$, by numerical calculation, we can prove that the configuration of the vectors is $(d, n, s, t) = (2^{k-1} (k-2)+1, 2^k \cdot k, s(k), 3)$, where we have $s(k) = 4k-3$ for every odd integer $k$ and we have $s(k) = 4k-5$ for every even integer $k$. Incidentally, the configuration of the $k$-dimensional regular $2 k$-cell polytope is $(d, n, s, t) = (k, 2^k, k, 3)$.\\

\quad\\

\subsection{$k$-dimensional regular $2^k$-cell polytopes (Cocktail party graph $CP(k)$)}\label{sec-cp}
The graph from $k$-dimensional regular $2^k$-cell polytope is the cocktail party graph $CP(k)$. For the $(d, n, s, t)$-configurations of the $k$-dimensional regular $2^k$-cell polytopes and those of the vectors of the constructed crystal lattices, we have the following table:\vspace{-0.07in}

\begin{table}[h]
\begin{center}
\begin{tabular}{ccccc}
\multicolumn{5}{c}{$2^k$-cell polytope}\\
\hline
$k$ & $v$ & $e$ & $s$ & $t$\\
\hline
$3$ & $6$ & $12$ & $2$ & $3$\\
$4$ & $8$ & $24$ & $2$ & $3$\\
$5$ & $10$ & $40$ & $2$ & $3$\\
$6$ & $12$ & $60$ & $2$ & $3$\\
$7$ & $14$ & $84$ & $2$ & $3$\\
$8$ & $16$ & $112$ & $2$ & $3$\\
$9$ & $18$ & $144$ & $2$ & $3$\\
\hline
\end{tabular}
\ $\Rightarrow$ \
\begin{tabular}{cccc}
\multicolumn{4}{c}{\small Vectors of lattice}\\
\hline
$d$ & $n$ & $s$ & $t$\\
\hline
$7$ & $24$ & $\underline{9}$ & $3$\\
$17$ & $48$ & $10$ & $3$\\
$31$ & $80$ & $10$ & $3$\\
$49$ & $120$ & $10$ & $3$\\
$71$ & $168$ & $10$ & $3$\\
$97$ & $224$ & $10$ & $3$\\
$127$ & $288$ & $10$ & $3$\\
\hline
\end{tabular}
\qquad
\begin{tabular}{ccccc}
\multicolumn{5}{c}{$2^k$-cell polytope}\\
\hline
$k$ & $v$ & $e$ & $s$ & $t$\\
\hline
$10$ & $20$ & $180$ & $2$ & $3$\\
$11$ & $22$ & $220$ & $2$ & $3$\\
$12$ & $23$ & $264$ & $2$ & $3$\\
$13$ & $26$ & $312$ & $2$ & $3$\\
$14$ & $28$ & $364$ & $2$ & $3$\\
$15$ & $30$ & $420$ & $2$ & $3$\\
$16$ & $32$ & $480$ & $2$ & $3$\\
\hline
\end{tabular}
\ $\Rightarrow$ \
\begin{tabular}{cccc}
\multicolumn{4}{c}{\small Vectors of lattice}\\
\hline
$d$ & $n$ & $s$ & $t$\\
\hline
$161$ & $360$ & $10$ & $3$\\
$199$ & $440$ & $10$ & $3$\\
$241$ & $528$ & $10$ & $3$\\
$287$ & $624$ & $10$ & $3$\\
$337$ & $728$ & $10$ & $3$\\
$391$ & $840$ & $10$ & $3$\\
$449$ & $960$ & $10$ & $3$\\
\hline
\end{tabular}
\end{center}\quad
\caption{$k$-dimensional regular $2^k$-cell polytope}\vspace{-0.27in}
\end{table}

\noindent
For every integer $3 \leqslant k \leqslant 16$, by numerical calculation, we can prove that we have the distance set $A(X) = \{ -1, \: \pm (2k-1) / (2 d(k)), \: \pm (k-1) / d(k), \: \pm 1 / d(k), \: \pm 1 / (2 d(k)), \: 0 \}$ (If $k=3$, $A(X)$ does not include $0$.) and that the configuration of the vectors is $(d, n, s, t) = (d(k), 4 k (k-1), 10, 3)$ (If $k=3$, then $s=9$.), where we denote $d(k) := 2 k^2 - 4 k + 1$. Furthermore, for each $x \in X$, we have the following table:
\begin{center}
\begin{tabular}{c|cccccc}
$a$ & $-1$ & $\pm (2k-1) / (2 d(k))$ & $\pm (k-1) / d(k)$ & $\pm 1 / d(k)$ & $\pm 1 / (2 d(k))$ & $0$\\
\hline
$|A_x(X, a)|$ & $1$ & $4 (k-2)$ & $2$ & $1$ & $4 (k-2)$ & $4 (k-2) (k-3)$\\
\end{tabular}
\end{center}

Incidentally, the configuration of the $k$-dimensional regular $2^k$-cell polytope is $(d, n, s, t) = (k, 2 k, 2, 3)$.

\newpage

\section{Complete bipartite graphs}\label{sec-bip}
In this section, we consider complete bipartite graph $K_{m_1, m_2}$ with two parts of $m_1$ vertices and $m_2$ vertices. We may assume that $m_1 \leqslant m_2$. Here, bipartite graph is a graph with two parts, where the pair of vertices from the same part is not connected with any edges. Note that complete bipartite graphs are edge-transitive, and they are not vertex-transitive if $m_1 < m_2$. Moreover, if $m_1 = 2$, then the constructed crystal lattice degenerate into the $m_2 - 1$-dimensional diamond crystal lattice (See Section \ref{sec-dia}). Thus, we assume that $m_1 \geqslant 3$.

Now, for the $(d, n, s, t)$-configurations of the vectors of the constructed crystal lattices, we have the following table:

{\small
\begin{center}
\begin{tabular}{ccccccccc}
\hline
$m_1$ & $m_2$ & $v$ & $e$ & \ & $d$ & $n$ & $s$ & $t$\\
\hline
$3$ & $3$ & $6$ & $9$ && $4$ & $18$ & $5$ & $3$\\
& $4$ & $7$ & $12$ && $6$ & $24$ & $7$ & $3$\\
& $5$ & $8$ & $15$ && $8$ & $30$ & $7$ & $3$\\
& $6$ & $9$ & $18$ && $10$ & $36$ & $7$ & $3$\\
& $7$ & $10$ & $21$ && $12$ & $42$ & $7$ & $3$\\
& $8$ & $11$ & $24$ && $14$ & $48$ & $7$ & $3$\\
& $9$ & $12$ & $27$ && $16$ & $54$ & $7$ & $3$\\
& $10$ & $13$ & $30$ && $18$ & $60$ & $7$ & $3$\\
& $11$ & $14$ & $33$ && $20$ & $66$ & $7$ & $3$\\
& $12$ & $15$ & $36$ && $22$ & $72$ & $7$ & $3$\\
& $13$ & $16$ & $39$ && $24$ & $78$ & $7$ & $3$\\
& $14$ & $17$ & $42$ && $26$ & $84$ & $7$ & $3$\\
& $15$ & $18$ & $45$ && $28$ & $90$ & $7$ & $3$\\
& $16$ & $19$ & $48$ && $30$ & $96$ & $7$ & $3$\\
& $17$ & $20$ & $51$ && $32$ & $102$ & $7$ & $3$\\
& $18$ & $21$ & $54$ && $34$ & $108$ & $7$ & $3$\\
& $19$ & $22$ & $57$ && $36$ & $114$ & $7$ & $3$\\
& $20$ & $23$ & $60$ && $38$ & $120$ & $7$ & $3$\\
& $21$ & $24$ & $63$ && $40$ & $126$ & $7$ & $3$\\
& $22$ & $25$ & $66$ && $42$ & $132$ & $7$ & $3$\\
& $23$ & $26$ & $69$ && $44$ & $138$ & $7$ & $3$\\
& $24$ & $27$ & $72$ && $46$ & $144$ & $7$ & $3$\\
& $25$ & $28$ & $75$ && $48$ & $150$ & $7$ & $3$\\
& $26$ & $29$ & $78$ && $50$ & $156$ & $7$ & $3$\\
& $27$ & $30$ & $81$ && $52$ & $162$ & $7$ & $3$\\
& $28$ & $31$ & $84$ && $54$ & $168$ & $7$ & $3$\\
& $29$ & $32$ & $87$ && $56$ & $174$ & $7$ & $3$\\
& $30$ & $33$ & $90$ && $58$ & $180$ & $7$ & $3$\\
& $31$ & $34$ & $93$ && $60$ & $186$ & $7$ & $3$\\
& $32$ & $35$ & $96$ && $62$ & $192$ & $7$ & $3$\\
& $33$ & $36$ & $99$ && $64$ & $198$ & $7$ & $3$\\
& $34$ & $37$ & $102$ && $66$ & $204$ & $7$ & $3$\\
& $35$ & $38$ & $105$ && $68$ & $210$ & $7$ & $3$\\
& $36$ & $39$ & $108$ && $70$ & $216$ & $7$ & $3$\\
& $37$ & $40$ & $111$ && $72$ & $222$ & $7$ & $3$\\
& $38$ & $41$ & $114$ && $74$ & $228$ & $7$ & $3$\\
& $39$ & $42$ & $117$ && $76$ & $234$ & $7$ & $3$\\
& $40$ & $43$ & $120$ && $78$ & $240$ & $7$ & $3$\\
& $41$ & $44$ & $123$ && $80$ & $246$ & $7$ & $3$\\
& $42$ & $45$ & $126$ && $82$ & $252$ & $7$ & $3$\\
& $43$ & $46$ & $129$ && $84$ & $258$ & $7$ & $3$\\
& $44$ & $47$ & $132$ && $86$ & $264$ & $7$ & $3$\\
& $45$ & $48$ & $135$ && $88$ & $270$ & $7$ & $3$\\
& $46$ & $49$ & $138$ && $90$ & $276$ & $7$ & $3$\\
& $47$ & $50$ & $141$ && $92$ & $282$ & $7$ & $3$\\
& $48$ & $51$ & $144$ && $94$ & $288$ & $7$ & $3$\\
& $49$ & $52$ & $147$ && $96$ & $294$ & $7$ & $3$\\
& $50$ & $53$ & $150$ && $98$ & $300$ & $7$ & $3$\\
& $51$ & $54$ & $153$ && $100$ & $306$ & $7$ & $3$\\
& $52$ & $55$ & $156$ && $102$ & $312$ & $7$ & $3$\\
\hline
\end{tabular}
\qquad
\begin{tabular}{ccccccccc}
\hline
$m_1$ & $m_2$ & $v$ & $e$ & \ & $d$ & $n$ & $s$ & $t$\\
\hline
$3$ & $53$ & $56$ & $159$ && $104$ & $318$ & $7$ & $3$\\
& $54$ & $57$ & $162$ && $106$ & $324$ & $7$ & $3$\\
& $55$ & $58$ & $165$ && $108$ & $330$ & $7$ & $3$\\
& $56$ & $59$ & $168$ && $110$ & $336$ & $7$ & $3$\\
& $57$ & $60$ & $171$ && $112$ & $342$ & $7$ & $3$\\
& $58$ & $61$ & $174$ && $114$ & $348$ & $7$ & $3$\\
& $59$ & $62$ & $177$ && $116$ & $354$ & $7$ & $3$\\
& $60$ & $63$ & $180$ && $118$ & $360$ & $7$ & $3$\\
& $61$ & $64$ & $183$ && $120$ & $366$ & $7$ & $3$\\
& $62$ & $65$ & $186$ && $122$ & $372$ & $7$ & $3$\\
& $63$ & $66$ & $189$ && $124$ & $378$ & $7$ & $3$\\
& $64$ & $67$ & $192$ && $126$ & $384$ & $7$ & $3$\\
& $65$ & $68$ & $195$ && $128$ & $390$ & $7$ & $3$\\
& $66$ & $69$ & $198$ && $130$ & $396$ & $7$ & $3$\\
& $67$ & $70$ & $201$ && $132$ & $402$ & $7$ & $3$\\
& $68$ & $71$ & $204$ && $134$ & $408$ & $7$ & $3$\\
& $69$ & $72$ & $207$ && $136$ & $414$ & $7$ & $3$\\
& $70$ & $73$ & $210$ && $138$ & $420$ & $7$ & $3$\\
& $71$ & $74$ & $213$ && $140$ & $426$ & $7$ & $3$\\
& $72$ & $75$ & $216$ && $142$ & $432$ & $7$ & $3$\\
& $73$ & $76$ & $219$ && $144$ & $438$ & $7$ & $3$\\
& $74$ & $77$ & $222$ && $146$ & $444$ & $7$ & $3$\\
& $75$ & $78$ & $225$ && $148$ & $450$ & $7$ & $3$\\
& $76$ & $79$ & $228$ && $150$ & $456$ & $7$ & $3$\\
& $77$ & $80$ & $231$ && $152$ & $462$ & $7$ & $3$\\
& $78$ & $81$ & $234$ && $154$ & $468$ & $7$ & $3$\\
& $79$ & $82$ & $237$ && $156$ & $474$ & $7$ & $3$\\
& $80$ & $83$ & $240$ && $158$ & $480$ & $7$ & $3$\\
& $81$ & $84$ & $243$ && $160$ & $486$ & $7$ & $3$\\
& $82$ & $85$ & $246$ && $162$ & $492$ & $7$ & $3$\\
& $83$ & $86$ & $249$ && $164$ & $498$ & $7$ & $3$\\
& $84$ & $87$ & $252$ && $166$ & $504$ & $7$ & $3$\\
& $85$ & $88$ & $255$ && $168$ & $510$ & $7$ & $3$\\
& $86$ & $89$ & $258$ && $170$ & $516$ & $7$ & $3$\\
& $87$ & $90$ & $261$ && $172$ & $522$ & $7$ & $3$\\
& $88$ & $91$ & $264$ && $174$ & $528$ & $7$ & $3$\\
& $89$ & $92$ & $267$ && $176$ & $534$ & $7$ & $3$\\
& $90$ & $93$ & $270$ && $178$ & $540$ & $7$ & $3$\\
& $91$ & $94$ & $273$ && $180$ & $546$ & $7$ & $3$\\
& $92$ & $95$ & $276$ && $182$ & $552$ & $7$ & $3$\\
& $93$ & $96$ & $279$ && $184$ & $558$ & $7$ & $3$\\
& $94$ & $97$ & $282$ && $186$ & $564$ & $7$ & $3$\\
& $95$ & $98$ & $285$ && $188$ & $570$ & $7$ & $3$\\
& $96$ & $99$ & $288$ && $190$ & $576$ & $7$ & $3$\\
& $97$ & $100$ & $291$ && $192$ & $582$ & $7$ & $3$\\
& $98$ & $101$ & $294$ && $194$ & $588$ & $7$ & $3$\\
& $99$ & $102$ & $297$ && $196$ & $594$ & $7$ & $3$\\
& $100$ & $103$ & $300$ && $198$ & $600$ & $7$ & $3$\\
& $101$ & $104$ & $303$ && $200$ & $606$ & $7$ & $3$\\
\hline\\
\end{tabular}
\end{center}}

\newpage

{\small
\begin{center}
\begin{tabular}{ccccccccc}
\hline
$m_1$ & $m_2$ & $v$ & $e$ && $d$ & $n$ & $s$ & $t$\\
\hline
$4$ & $4$ & $8$ & $16$ && $9$ & $32$ & $5$ & $3$\\
& $5$ & $9$ & $20$ && $12$ & $40$ & $7$ & $3$\\
& $6$ & $10$ & $24$ && $15$ & $48$ & $7$ & $3$\\
& $7$ & $11$ & $28$ && $18$ & $56$ & $7$ & $3$\\
& $8$ & $12$ & $32$ && $21$ & $64$ & $7$ & $3$\\
& $9$ & $13$ & $36$ && $24$ & $72$ & $7$ & $3$\\
& $10$ & $14$ & $40$ && $27$ & $80$ & $7$ & $3$\\
& $11$ & $15$ & $44$ && $30$ & $88$ & $7$ & $3$\\
& $12$ & $16$ & $48$ && $33$ & $96$ & $7$ & $3$\\
& $13$ & $17$ & $52$ && $36$ & $104$ & $7$ & $3$\\
& $14$ & $18$ & $56$ && $39$ & $112$ & $7$ & $3$\\
& $15$ & $19$ & $60$ && $42$ & $120$ & $7$ & $3$\\
& $16$ & $20$ & $64$ && $45$ & $128$ & $7$ & $3$\\
& $17$ & $21$ & $68$ && $48$ & $136$ & $7$ & $3$\\
& $18$ & $22$ & $72$ && $51$ & $144$ & $7$ & $3$\\
& $19$ & $23$ & $76$ && $54$ & $152$ & $7$ & $3$\\
& $20$ & $24$ & $80$ && $57$ & $160$ & $7$ & $3$\\
& $21$ & $25$ & $84$ && $60$ & $168$ & $7$ & $3$\\
& $22$ & $26$ & $88$ && $63$ & $176$ & $7$ & $3$\\
& $23$ & $27$ & $92$ && $66$ & $184$ & $7$ & $3$\\
& $24$ & $28$ & $96$ && $69$ & $192$ & $7$ & $3$\\
& $25$ & $29$ & $100$ && $72$ & $200$ & $7$ & $3$\\
& $26$ & $30$ & $104$ && $75$ & $208$ & $7$ & $3$\\
& $27$ & $31$ & $108$ && $78$ & $216$ & $7$ & $3$\\
& $28$ & $32$ & $112$ && $81$ & $224$ & $7$ & $3$\\
& $29$ & $33$ & $116$ && $84$ & $232$ & $7$ & $3$\\
& $30$ & $34$ & $120$ && $87$ & $240$ & $7$ & $3$\\
& $31$ & $35$ & $124$ && $90$ & $248$ & $7$ & $3$\\
& $32$ & $36$ & $128$ && $93$ & $256$ & $7$ & $3$\\
& $33$ & $37$ & $132$ && $96$ & $264$ & $7$ & $3$\\
& $34$ & $38$ & $136$ && $99$ & $272$ & $7$ & $3$\\
& $35$ & $39$ & $140$ && $102$ & $280$ & $7$ & $3$\\
\hline\\
\hline
$m_1$ & $m_2$ & $v$ & $e$ && $d$ & $n$ & $s$ & $t$\\
\hline
$5$ & $5$ & $10$ & $25$ && $16$ & $50$ & $5$ & $3$\\
& $6$ & $11$ & $30$ && $20$ & $60$ & $7$ & $3$\\
& $7$ & $12$ & $35$ && $24$ & $70$ & $7$ & $3$\\
& $8$ & $13$ & $40$ && $28$ & $80$ & $7$ & $3$\\
& $9$ & $14$ & $45$ && $32$ & $90$ & $7$ & $3$\\
& $10$ & $15$ & $50$ && $36$ & $100$ & $7$ & $3$\\
& $11$ & $16$ & $55$ && $40$ & $110$ & $7$ & $3$\\
& $12$ & $17$ & $60$ && $44$ & $120$ & $7$ & $3$\\
& $13$ & $18$ & $65$ && $48$ & $130$ & $7$ & $3$\\
& $14$ & $19$ & $70$ && $52$ & $140$ & $7$ & $3$\\
& $15$ & $20$ & $75$ && $56$ & $150$ & $7$ & $3$\\
& $16$ & $21$ & $80$ && $60$ & $160$ & $7$ & $3$\\
& $17$ & $22$ & $85$ && $64$ & $170$ & $7$ & $3$\\
& $18$ & $23$ & $90$ && $68$ & $180$ & $7$ & $3$\\
& $19$ & $24$ & $95$ && $72$ & $190$ & $7$ & $3$\\
& $20$ & $25$ & $100$ && $76$ & $200$ & $7$ & $3$\\
& $21$ & $26$ & $105$ && $80$ & $210$ & $7$ & $3$\\
& $22$ & $27$ & $110$ && $84$ & $220$ & $7$ & $3$\\
& $23$ & $28$ & $115$ && $88$ & $230$ & $7$ & $3$\\
& $24$ & $29$ & $120$ && $92$ & $240$ & $7$ & $3$\\
& $25$ & $30$ & $125$ && $96$ & $250$ & $7$ & $3$\\
& $26$ & $31$ & $130$ && $100$ & $260$ & $7$ & $3$\\
& $27$ & $32$ & $135$ && $104$ & $270$ & $7$ & $3$\\
& $28$ & $33$ & $140$ && $108$ & $280$ & $7$ & $3$\\
\hline
\end{tabular}
\qquad
\begin{tabular}{ccccccccc}
\hline
$m_1$ & $m_2$ & $v$ & $e$ && $d$ & $n$ & $s$ & $t$\\
\hline
$4$ & $36$ & $40$ & $144$ && $105$ & $288$ & $7$ & $3$\\
& $37$ & $41$ & $148$ && $108$ & $296$ & $7$ & $3$\\
& $38$ & $42$ & $152$ && $111$ & $304$ & $7$ & $3$\\
& $39$ & $43$ & $156$ && $114$ & $312$ & $7$ & $3$\\
& $40$ & $44$ & $160$ && $117$ & $320$ & $7$ & $3$\\
& $41$ & $45$ & $164$ && $120$ & $328$ & $7$ & $3$\\
& $42$ & $46$ & $168$ && $123$ & $336$ & $7$ & $3$\\
& $43$ & $47$ & $172$ && $126$ & $344$ & $7$ & $3$\\
& $44$ & $48$ & $176$ && $129$ & $352$ & $7$ & $3$\\
& $45$ & $49$ & $180$ && $132$ & $360$ & $7$ & $3$\\
& $46$ & $50$ & $184$ && $135$ & $368$ & $7$ & $3$\\
& $47$ & $51$ & $188$ && $138$ & $376$ & $7$ & $3$\\
& $48$ & $52$ & $192$ && $141$ & $384$ & $7$ & $3$\\
& $49$ & $53$ & $196$ && $144$ & $392$ & $7$ & $3$\\
& $50$ & $54$ & $200$ && $147$ & $400$ & $7$ & $3$\\
& $51$ & $55$ & $204$ && $150$ & $408$ & $7$ & $3$\\
& $52$ & $56$ & $208$ && $153$ & $416$ & $7$ & $3$\\
& $53$ & $57$ & $212$ && $156$ & $424$ & $7$ & $3$\\
& $54$ & $58$ & $216$ && $159$ & $432$ & $7$ & $3$\\
& $55$ & $59$ & $220$ && $162$ & $440$ & $7$ & $3$\\
& $56$ & $60$ & $224$ && $165$ & $448$ & $7$ & $3$\\
& $57$ & $61$ & $228$ && $168$ & $456$ & $7$ & $3$\\
& $58$ & $62$ & $232$ && $171$ & $464$ & $7$ & $3$\\
& $59$ & $63$ & $236$ && $174$ & $472$ & $7$ & $3$\\
& $60$ & $64$ & $240$ && $177$ & $480$ & $7$ & $3$\\
& $61$ & $65$ & $244$ && $180$ & $488$ & $7$ & $3$\\
& $62$ & $66$ & $248$ && $183$ & $496$ & $7$ & $3$\\
& $63$ & $67$ & $252$ && $186$ & $504$ & $7$ & $3$\\
& $64$ & $68$ & $256$ && $189$ & $512$ & $7$ & $3$\\
& $65$ & $69$ & $260$ && $192$ & $520$ & $7$ & $3$\\
& $66$ & $70$ & $264$ && $195$ & $528$ & $7$ & $3$\\
& $67$ & $71$ & $268$ && $198$ & $536$ & $7$ & $3$\\
\hline\\
\hline
$m_1$ & $m_2$ & $v$ & $e$ && $d$ & $n$ & $s$ & $t$\\
\hline
$5$ & $29$ & $34$ & $145$ && $112$ & $290$ & $7$ & $3$\\
& $30$ & $35$ & $150$ && $116$ & $300$ & $7$ & $3$\\
& $31$ & $36$ & $155$ && $120$ & $310$ & $7$ & $3$\\
& $32$ & $37$ & $160$ && $124$ & $320$ & $7$ & $3$\\
& $33$ & $38$ & $165$ && $128$ & $330$ & $7$ & $3$\\
& $34$ & $39$ & $170$ && $132$ & $340$ & $7$ & $3$\\
& $35$ & $40$ & $175$ && $136$ & $350$ & $7$ & $3$\\
& $36$ & $41$ & $180$ && $140$ & $360$ & $7$ & $3$\\
& $37$ & $42$ & $185$ && $144$ & $370$ & $7$ & $3$\\
& $38$ & $43$ & $190$ && $148$ & $380$ & $7$ & $3$\\
& $39$ & $44$ & $195$ && $152$ & $390$ & $7$ & $3$\\
& $40$ & $45$ & $200$ && $156$ & $400$ & $7$ & $3$\\
& $41$ & $46$ & $205$ && $160$ & $410$ & $7$ & $3$\\
& $42$ & $47$ & $210$ && $164$ & $420$ & $7$ & $3$\\
& $43$ & $48$ & $215$ && $168$ & $430$ & $7$ & $3$\\
& $44$ & $49$ & $220$ && $172$ & $440$ & $7$ & $3$\\
& $45$ & $50$ & $225$ && $176$ & $450$ & $7$ & $3$\\
& $46$ & $51$ & $230$ && $180$ & $460$ & $7$ & $3$\\
& $47$ & $52$ & $235$ && $184$ & $470$ & $7$ & $3$\\
& $48$ & $53$ & $240$ && $188$ & $480$ & $7$ & $3$\\
& $49$ & $54$ & $245$ && $192$ & $490$ & $7$ & $3$\\
& $50$ & $55$ & $250$ && $196$ & $500$ & $7$ & $3$\\
& $51$ & $56$ & $255$ && $200$ & $510$ & $7$ & $3$\\
\hline\\
\end{tabular}
\end{center}}

\newpage

\begin{center}
\begin{tabular}{ccccccccc}
\hline
$m_1$ & $m_2$ & $v$ & $e$ && $d$ & $n$ & $s$ & $t$\\
\hline
$6$ & $6$ & $12$ & $36$ && $25$ & $72$ & $5$ & $3$\\
& $7$ & $13$ & $42$ && $30$ & $84$ & $7$ & $3$\\
& $8$ & $14$ & $48$ && $35$ & $96$ & $7$ & $3$\\
& $9$ & $15$ & $54$ && $40$ & $108$ & $7$ & $3$\\
& $10$ & $16$ & $60$ && $45$ & $120$ & $7$ & $3$\\
& $11$ & $17$ & $66$ && $50$ & $132$ & $7$ & $3$\\
& $12$ & $18$ & $72$ && $55$ & $144$ & $7$ & $3$\\
& $13$ & $19$ & $78$ && $60$ & $156$ & $7$ & $3$\\
& $14$ & $20$ & $84$ && $65$ & $168$ & $7$ & $3$\\
& $15$ & $21$ & $90$ && $70$ & $180$ & $7$ & $3$\\
& $16$ & $22$ & $96$ && $75$ & $192$ & $7$ & $3$\\
& $17$ & $23$ & $102$ && $80$ & $204$ & $7$ & $3$\\
& $18$ & $24$ & $108$ && $85$ & $216$ & $7$ & $3$\\
& $19$ & $25$ & $114$ && $90$ & $228$ & $7$ & $3$\\
& $20$ & $26$ & $120$ && $95$ & $240$ & $7$ & $3$\\
& $21$ & $27$ & $126$ && $100$ & $252$ & $7$ & $3$\\
& $22$ & $28$ & $132$ && $105$ & $264$ & $7$ & $3$\\
& $23$ & $29$ & $138$ && $110$ & $276$ & $7$ & $3$\\
\hline\\
\hline
$m_1$ & $m_2$ & $v$ & $e$ && $d$ & $n$ & $s$ & $t$\\
\hline
$7$ & $7$ & $14$ & $49$ && $36$ & $98$ & $5$ & $3$\\
& $8$ & $15$ & $56$ && $42$ & $112$ & $7$ & $3$\\
& $9$ & $16$ & $63$ && $48$ & $126$ & $7$ & $3$\\
& $10$ & $17$ & $70$ && $54$ & $140$ & $7$ & $3$\\
& $11$ & $18$ & $77$ && $60$ & $154$ & $7$ & $3$\\
& $12$ & $19$ & $84$ && $66$ & $168$ & $7$ & $3$\\
& $13$ & $20$ & $91$ && $72$ & $182$ & $7$ & $3$\\
& $14$ & $21$ & $98$ && $78$ & $196$ & $7$ & $3$\\
& $15$ & $22$ & $105$ && $84$ & $210$ & $7$ & $3$\\
& $16$ & $23$ & $112$ && $90$ & $224$ & $7$ & $3$\\
& $17$ & $24$ & $119$ && $96$ & $238$ & $7$ & $3$\\
& $18$ & $25$ & $126$ && $102$ & $252$ & $7$ & $3$\\
& $19$ & $26$ & $133$ && $108$ & $266$ & $7$ & $3$\\
& $20$ & $27$ & $140$ && $114$ & $280$ & $7$ & $3$\\
\hline\\
\hline
$m_1$ & $m_2$ & $v$ & $e$ && $d$ & $n$ & $s$ & $t$\\
\hline
$8$ & $8$ & $16$ & $64$ && $49$ & $128$ & $5$ & $3$\\
& $9$ & $17$ & $72$ && $56$ & $144$ & $7$ & $3$\\
& $10$ & $18$ & $80$ && $63$ & $160$ & $7$ & $3$\\
& $11$ & $19$ & $88$ && $70$ & $176$ & $7$ & $3$\\
& $12$ & $20$ & $96$ && $77$ & $192$ & $7$ & $3$\\
& $13$ & $21$ & $104$ && $84$ & $208$ & $7$ & $3$\\
& $14$ & $22$ & $112$ && $91$ & $224$ & $7$ & $3$\\
& $15$ & $23$ & $120$ && $98$ & $240$ & $7$ & $3$\\
& $16$ & $24$ & $128$ && $105$ & $256$ & $7$ & $3$\\
& $17$ & $25$ & $136$ && $112$ & $272$ & $7$ & $3$\\
& $18$ & $26$ & $144$ && $119$ & $288$ & $7$ & $3$\\
\hline
\end{tabular}
\qquad
\begin{tabular}{ccccccccc}
\hline
$m_1$ & $m_2$ & $v$ & $e$ && $d$ & $n$ & $s$ & $t$\\
\hline
$6$ & $24$ & $30$ & $144$ && $115$ & $288$ & $7$ & $3$\\
& $25$ & $31$ & $150$ && $120$ & $300$ & $7$ & $3$\\
& $26$ & $32$ & $156$ && $125$ & $312$ & $7$ & $3$\\
& $27$ & $33$ & $162$ && $130$ & $324$ & $7$ & $3$\\
& $28$ & $34$ & $168$ && $135$ & $336$ & $7$ & $3$\\
& $29$ & $35$ & $174$ && $140$ & $348$ & $7$ & $3$\\
& $30$ & $36$ & $180$ && $145$ & $360$ & $7$ & $3$\\
& $31$ & $37$ & $186$ && $150$ & $372$ & $7$ & $3$\\
& $32$ & $38$ & $192$ && $155$ & $384$ & $7$ & $3$\\
& $33$ & $39$ & $198$ && $160$ & $396$ & $7$ & $3$\\
& $34$ & $40$ & $204$ && $165$ & $408$ & $7$ & $3$\\
& $35$ & $41$ & $210$ && $170$ & $420$ & $7$ & $3$\\
& $36$ & $42$ & $216$ && $175$ & $432$ & $7$ & $3$\\
& $37$ & $43$ & $222$ && $180$ & $444$ & $7$ & $3$\\
& $38$ & $44$ & $228$ && $185$ & $456$ & $7$ & $3$\\
& $39$ & $45$ & $234$ && $190$ & $468$ & $7$ & $3$\\
& $40$ & $46$ & $240$ && $195$ & $480$ & $7$ & $3$\\
& $41$ & $47$ & $246$ && $200$ & $492$ & $7$ & $3$\\
\hline\\
\hline
$m_1$ & $m_2$ & $v$ & $e$ && $d$ & $n$ & $s$ & $t$\\
\hline
$7$ & $21$ & $28$ & $147$ && $120$ & $294$ & $7$ & $3$\\
& $22$ & $29$ & $154$ && $126$ & $308$ & $7$ & $3$\\
& $23$ & $30$ & $161$ && $132$ & $322$ & $7$ & $3$\\
& $24$ & $31$ & $168$ && $138$ & $336$ & $7$ & $3$\\
& $25$ & $32$ & $175$ && $144$ & $350$ & $7$ & $3$\\
& $26$ & $33$ & $182$ && $150$ & $364$ & $7$ & $3$\\
& $27$ & $34$ & $189$ && $156$ & $378$ & $7$ & $3$\\
& $28$ & $35$ & $196$ && $162$ & $392$ & $7$ & $3$\\
& $29$ & $36$ & $203$ && $168$ & $406$ & $7$ & $3$\\
& $30$ & $37$ & $210$ && $174$ & $420$ & $7$ & $3$\\
& $31$ & $38$ & $217$ && $180$ & $434$ & $7$ & $3$\\
& $32$ & $39$ & $224$ && $186$ & $448$ & $7$ & $3$\\
& $33$ & $40$ & $231$ && $192$ & $462$ & $7$ & $3$\\
& $34$ & $41$ & $238$ && $198$ & $476$ & $7$ & $3$\\
\hline\\
\hline
$m_1$ & $m_2$ & $v$ & $e$ && $d$ & $n$ & $s$ & $t$\\
\hline
$8$ & $19$ & $27$ & $152$ && $126$ & $304$ & $7$ & $3$\\
& $20$ & $28$ & $160$ && $133$ & $320$ & $7$ & $3$\\
& $21$ & $29$ & $168$ && $140$ & $336$ & $7$ & $3$\\
& $22$ & $30$ & $176$ && $147$ & $352$ & $7$ & $3$\\
& $23$ & $31$ & $184$ && $154$ & $368$ & $7$ & $3$\\
& $24$ & $32$ & $192$ && $161$ & $384$ & $7$ & $3$\\
& $25$ & $33$ & $200$ && $168$ & $400$ & $7$ & $3$\\
& $26$ & $34$ & $208$ && $175$ & $416$ & $7$ & $3$\\
& $27$ & $35$ & $216$ && $182$ & $432$ & $7$ & $3$\\
& $28$ & $36$ & $224$ && $189$ & $448$ & $7$ & $3$\\
& $29$ & $37$ & $232$ && $196$ & $464$ & $7$ & $3$\\
\hline
\end{tabular}
\end{center}

\newpage

\begin{table}[h]
\begin{center}
\begin{tabular}{ccccccccc}
\hline
$m_1$ & $m_2$ & $v$ & $e$ && $d$ & $n$ & $s$ & $t$\\
\hline
$9$ & $9$ & $18$ & $81$ && $64$ & $162$ & $5$ & $3$\\
& $10$ & $19$ & $90$ && $72$ & $180$ & $7$ & $3$\\
& $11$ & $20$ & $99$ && $80$ & $198$ & $7$ & $3$\\
& $12$ & $21$ & $108$ && $88$ & $216$ & $7$ & $3$\\
& $13$ & $22$ & $117$ && $96$ & $234$ & $7$ & $3$\\
& $14$ & $23$ & $126$ && $104$ & $252$ & $7$ & $3$\\
& $15$ & $24$ & $135$ && $112$ & $270$ & $7$ & $3$\\
& $16$ & $25$ & $144$ && $120$ & $288$ & $7$ & $3$\\
& $17$ & $26$ & $153$ && $128$ & $306$ & $7$ & $3$\\
& $18$ & $27$ & $162$ && $136$ & $324$ & $7$ & $3$\\
& $19$ & $28$ & $171$ && $144$ & $342$ & $7$ & $3$\\
& $20$ & $29$ & $180$ && $152$ & $360$ & $7$ & $3$\\
& $21$ & $30$ & $189$ && $160$ & $378$ & $7$ & $3$\\
& $22$ & $31$ & $198$ && $168$ & $396$ & $7$ & $3$\\
& $23$ & $32$ & $207$ && $176$ & $414$ & $7$ & $3$\\
& $24$ & $33$ & $216$ && $184$ & $432$ & $7$ & $3$\\
& $25$ & $34$ & $225$ && $192$ & $450$ & $7$ & $3$\\
& $26$ & $35$ & $234$ && $200$ & $468$ & $7$ & $3$\\
\hline\\
\hline
$m_1$ & $m_2$ & $v$ & $e$ && $d$ & $n$ & $s$ & $t$\\
\hline
$10$ & $10$ & $20$ & $100$ && $81$ & $200$ & $5$ & $3$\\
& $11$ & $21$ & $110$ && $90$ & $220$ & $7$ & $3$\\
& $12$ & $22$ & $120$ && $99$ & $240$ & $7$ & $3$\\
& $13$ & $23$ & $130$ && $108$ & $260$ & $7$ & $3$\\
& $14$ & $24$ & $140$ && $117$ & $280$ & $7$ & $3$\\
& $15$ & $25$ & $150$ && $126$ & $300$ & $7$ & $3$\\
& $16$ & $26$ & $160$ && $135$ & $320$ & $7$ & $3$\\
& $17$ & $27$ & $170$ && $144$ & $340$ & $7$ & $3$\\
& $18$ & $28$ & $180$ && $153$ & $360$ & $7$ & $3$\\
& $19$ & $29$ & $190$ && $162$ & $380$ & $7$ & $3$\\
& $20$ & $30$ & $200$ && $171$ & $400$ & $7$ & $3$\\
& $21$ & $31$ & $210$ && $180$ & $420$ & $7$ & $3$\\
& $22$ & $32$ & $220$ && $189$ & $440$ & $7$ & $3$\\
& $23$ & $33$ & $230$ && $198$ & $460$ & $7$ & $3$\\
\hline
\end{tabular}
\qquad
\begin{tabular}{ccccccccc}
\hline
$m_1$ & $m_2$ & $v$ & $e$ && $d$ & $n$ & $s$ & $t$\\
\hline
$11$ & $11$ & $22$ & $121$ && $100$ & $242$ & $5$ & $3$\\
& $12$ & $23$ & $132$ && $110$ & $264$ & $7$ & $3$\\
& $13$ & $24$ & $143$ && $120$ & $286$ & $7$ & $3$\\
& $14$ & $25$ & $154$ && $130$ & $308$ & $7$ & $3$\\
& $15$ & $26$ & $165$ && $140$ & $330$ & $7$ & $3$\\
& $16$ & $27$ & $176$ && $150$ & $352$ & $7$ & $3$\\
& $17$ & $28$ & $187$ && $160$ & $374$ & $7$ & $3$\\
& $18$ & $29$ & $198$ && $170$ & $396$ & $7$ & $3$\\
& $19$ & $30$ & $209$ && $180$ & $418$ & $7$ & $3$\\
& $20$ & $31$ & $220$ && $190$ & $440$ & $7$ & $3$\\
& $21$ & $32$ & $231$ && $200$ & $462$ & $7$ & $3$\\
\hline\vspace{-0.1in}\\
\hline
$m_1$ & $m_2$ & $v$ & $e$ && $d$ & $n$ & $s$ & $t$\\
\hline
$12$ & $12$ & $24$ & $144$ && $121$ & $288$ & $5$ & $3$\\
& $13$ & $25$ & $156$ && $132$ & $312$ & $7$ & $3$\\
& $14$ & $26$ & $168$ && $143$ & $336$ & $7$ & $3$\\
& $15$ & $27$ & $180$ && $154$ & $360$ & $7$ & $3$\\
& $16$ & $28$ & $192$ && $165$ & $384$ & $7$ & $3$\\
& $17$ & $29$ & $204$ && $176$ & $408$ & $7$ & $3$\\
& $18$ & $30$ & $216$ && $187$ & $432$ & $7$ & $3$\\
& $19$ & $31$ & $228$ && $198$ & $456$ & $7$ & $3$\\
\hline\vspace{-0.1in}\\
\hline
$m_1$ & $m_2$ & $v$ & $e$ && $d$ & $n$ & $s$ & $t$\\
\hline
$13$ & $13$ & $26$ & $169$ && $144$ & $338$ & $5$ & $3$\\
& $14$ & $27$ & $182$ && $156$ & $364$ & $7$ & $3$\\
& $15$ & $28$ & $195$ && $168$ & $390$ & $7$ & $3$\\
& $16$ & $29$ & $208$ && $180$ & $416$ & $7$ & $3$\\
& $17$ & $30$ & $221$ && $192$ & $442$ & $7$ & $3$\\
\hline\vspace{-0.1in}\\
\hline
$m_1$ & $m_2$ & $v$ & $e$ && $d$ & $n$ & $s$ & $t$\\
\hline
$14$ & $14$ & $28$ & $196$ && $169$ & $392$ & $5$ & $3$\\
& $15$ & $29$ & $210$ && $182$ & $420$ & $7$ & $3$\\
& $16$ & $30$ & $224$ && $195$ & $448$ & $7$ & $3$\\
\hline\vspace{-0.1in}\\
\hline
$m_1$ & $m_2$ & $v$ & $e$ && $d$ & $n$ & $s$ & $t$\\
\hline
$15$ & $15$ & $30$ & $225$ && $196$ & $450$ & $5$ & $3$\\
\hline
\end{tabular}
\end{center}\quad
\caption{Complete bipartite graph}\vspace{-0.27in}
\end{table}

\noindent
For every integer $(m_1 - 1) (m_2 - 1) \leqslant 200$, by numerical calculation, we can prove that we have the distance set $A(X) = \{ -1, \: \pm 1/(m_1 - 1), \: \pm 1/(m_2 - 1), \: \pm 1/((m_1 - 1) (m_2 - 1)) \}$ and that the configuration of the vectors is $(d, n, s, t) = ((m_1 - 1) (m_2 - 1), 2 m_1 m_2, 7, 3)$ (If $m_1 = m_2$, then $s=5$). Furthermore, for each $x \in X$, we have the following table:
\begin{center}
\begin{tabular}{c|cccc}
$a$ & $-1$ & $\pm 1/(m_1 - 1)$ & $\pm 1/(m_1 - 1)$ & $\pm 1/((m_1 - 1) (m_2 - 1))$\\
\hline
$|A_x(X, a)|$ & $1$ & $m_1 - 1$ & $m_2 - 1$ & $(m_1 - 1) (m_2 - 1)$\\
\end{tabular}
\end{center}

\newpage

\section{Strongly regular graphs}\label{sec-sr}
In this section, we consider strongly regular graph, we refer to \cite[Chapter2]{CL}. Let $X_0$ be a finite graph, and let $\Gamma_x$ be a set of vertices of $X_0$ which are adjacent to $x$. $k$-regular graph with $v$ vertices is a strongly regular graph with parameters $(v, k, \lambda, \mu)$ if $i)$ $|\Gamma_x \cap \Gamma_y|$ is equal to the constant $\lambda$ for every pair of vertices $(x, y)$ which are adjacent to each other and $ii)$ $|\Gamma_x \cap \Gamma_y|$ is equal to the constant $\mu$ for every pair of vertices $(x, y)$ which are not adjacent. Note that the complement graph of a strongly regular graph with parameters $(v, k, \lambda, \mu)$ is a strongly regular graph with parameters $(v, v - k - 1, v - 2k + \mu - 2, v - 2k + \lambda)$.

\subsection{Triangular graph $T(m)$ (Johnson graph $J(m, 2)$)}\label{sec-tri}
For every positive integer $m \geqslant 3$, the triangular graph $T(m)$ is a strongly regular graph with parameters $(v, k, \lambda, \mu) = (m(m-1)/2, 2(m-2), m-2, 4)$, which is also called the Johnson graph $J(m, 2)$ and is vertex and edge-transitive. Let $\Omega$ be a $m$-point set, and let $V$ be a family of all the $2$-point subsets of $\Omega$. Then, we define $E$ as a set of pairs $(x, y)$ of elements of $V$ such that $|x \cap y| = 1$. Now, we can define the triangular graph $T(m)$ as the graph with $V$ and $E$ as a set of vertices and edges, respectively. Moreover, if $m = 3$, then the constructed crystal lattice degenerate into the $1$-dimensional standard crystal lattice. Thus, we assume that $m \geqslant 4$. For the $(d, n, s, t)$-configurations of the vectors of the constructed crystal lattices, we have the following table:\vspace{-0.05in}

\begin{table}[h]
\begin{center}
\begin{tabular}{ccccccc}
\multicolumn{7}{c}{Triangular graph $T(m)$}\\
\hline
$m$ & & $v$ & $k$ & $\lambda$ & $\mu$ & $e$\\
\hline
$4$ && $6$ & $4$ & $2$ & $4$ & $12$\\
$5$ && $10$ & $6$ & $3$ & $4$ & $30$\\
$6$ && $15$ & $8$ & $4$ & $4$ & $60$\\
$7$ && $21$ & $10$ & $5$ & $4$ & $105$\\
$8$ && $28$ & $12$ & $6$ & $4$ & $168$\\
$9$ && $36$ & $14$ & $7$ & $4$ & $252$\\
$10$ && $45$ & $16$ & $8$ & $4$ & $360$\\
$11$ && $55$ & $18$ & $9$ & $4$ & $495$\\
$12$ && $66$ & $20$ & $10$ & $4$ & $660$\\
$13$ && $78$ & $22$ & $11$ & $4$ & $858$\\
\hline
\end{tabular}
\quad $\Rightarrow$ \quad
\begin{tabular}{cccc}
\multicolumn{4}{c}{Vectors of lattice}\\
\hline
$d$ & $n$ & $s$ & $t$\\
\hline
$7$ & $24$ & $\underline{9}$ & $3$\\
$21$ & $60$ & $10$ & $3$\\
$46$ & $120$ & $10$ & $3$\\
$85$ & $210$ & $10$ & $3$\\
$141$ & $336$ & $10$ & $3$\\
$217$ & $504$ & $10$ & $3$\\
$316$ & $720$ & $10$ & $3$\\
$441$ & $990$ & $10$ & $3$\\
$595$ & $1320$ & $10$ & $3$\\
$781$ & $1716$ & $10$ & $3$\\
\hline
\end{tabular}
\end{center}\quad
\caption{Triangular graph $T(m)$}\vspace{-0.27in}
\end{table}

\noindent
For every integer $4 \leqslant m \leqslant 13$, by numerical calculation, we can prove that we have the distance set $A(X) = \{ -1, \: \pm (m+1) / (2 f(m)), \: \pm m / (2 f(m)), \: \pm 2 / (2 f(m)), \: \pm 1/(2 f(m)), \: 0 \}$ (If $m=4$, then $A(X)$ does not include $0$) and that the configuration of the vectors is $(d, n, s, t) = ((m-2) f(m) / 2, m (m-1) (m-2), 10, 3)$ (If $m=4$, then $s=9$), where we denote $f(m) := m^2 - 2m - 1$. Furthermore, for each $x \in X$, we have the following table:
\begin{center}
\begin{tabular}{c|cccccc}
$a$ & $-1$ & $\pm (m+1) / (2 f(m))$ & $\pm m / (2 f(m))$ & $\pm 2 / (2 f(m))$ & $\pm 1 / (2 f(m))$ & $0$\\
\hline
$|A_x(X, a)|$ & $1$ & $2 (m-2)$ & $2 (m-3)$ & $m-3$ & $2 (m-2) (m-3)$ & $m (m-3) (m-4)$
\end{tabular}
\end{center}

\quad

On the other hand, the complement graph of the triangular graph $T(m)$ is a strongly regular graph with parameters $(v, k, \lambda, \mu) = (m(m-1)/2, (m-2)(m-3)/2, (m-2)(m-3)/2, (m-2)(m-3)/2)$.  Moreover, the graph is empty if $m=3$, and the constructed crystal lattice degenerate into the $0$-dimensional crystal lattice if $m = 4$. Thus, we assume that $m \geqslant 5$. For the $(d, n, s, t)$-configurations of the vectors of the constructed crystal lattices, we have the following table:\vspace{-0.05in}

\begin{table}[h]
\begin{center}
\begin{tabular}{ccccccc}
\multicolumn{7}{c}{$\overline{T(m)}$}\\
\hline
$m$ & & $v$ & $k$ & $\lambda$ & $\mu$ & $e$\\
\hline
$5$ && $10$ & $3$ & $0$ & $1$ & $30$\\
$6$ && $15$ & $6$ & $1$ & $3$ & $60$\\
$7$ && $21$ & $10$ & $3$ & $6$ & $105$\\
$8$ && $28$ & $15$ & $6$ & $10$ & $168$\\
$9$ && $36$ & $21$ & $10$ & $15$ & $252$\\
$10$ && $45$ & $28$ & $15$ & $21$ & $630$\\
\hline
\end{tabular}
\quad $\Rightarrow$ \quad
\begin{tabular}{cccc}
\multicolumn{4}{c}{Vectors of lattice}\\
\hline
$d$ & $n$ & $s$ & $t$\\
\hline
$6$ & $30$ & $6$ & $3$\\
$31$ & $90$ & $10$ & $3$\\
$85$ & $210$ & $10$ & $3$\\
$183$ & $420$ & $10$ & $3$\\
$343$ & $756$ & $10$ & $3$\\
$586$ & $1260$ & $10$ & $3$\\
\hline
\end{tabular}
\end{center}\quad
\caption{Complement of the triangular graph $\overline{T(m)}$}\vspace{-0.27in}
\end{table}

\noindent
In particular, when $m=5$, we call the graph $\overline{T(5)}$ as {\bf Petersen graph}, and we have the distance set  $A(X) = \{ -1, \: \pm 1/2, \: \pm 1/4, \: 0 \}$.

\newpage

\subsection{Square lattice graph $L_2(m)$ (Hamming graph $H(2, m)$)}\label{sec-l2}
For every positive integer $m \geqslant 2$, the square lattice graph $L_2(m)$ is a strongly regular graph with parameters $(v, k, \lambda, \mu) = (m^2, 2(m-1), m-2, 2)$, which is also called the Hamming graph $H(2, m)$ and is vertex and edge-transitive. Let $\Omega$ be a $m$-point set, and let $V = \Omega \times \Omega$. Then, we define $E$ as a set of pairs $(x, y)$ of elements of $V$ such that $x_1 = y_1$ or $x_2 = y_2$ for $x = (x_1, x_2)$ and $y = (y_1, y_2)$. Now, we can define the square lattice graph $L_2(m)$ as the graph with $V$ and $E$ as a set of vertices and edges, respectively. Moreover, if $m = 2$, then the constructed crystal lattice degenerate into the $1$-dimensional standard crystal lattice. Thus, we assume that $m \geqslant 3$.

For the $(d, n, s, t)$-configurations of the vectors of the constructed crystal lattices, we have the following table:

\begin{table}[h]
\begin{center}
\begin{tabular}{ccccccc}
\multicolumn{7}{c}{Square lattice graph $L_2(m)$}\\
\hline
$m$ & & $v$ & $k$ & $\lambda$ & $\mu$ & $e$\\
\hline
$3$ && $9$ & $4$ & $1$ & $2$ & $18$\\
$4$ && $16$ & $6$ & $2$ & $2$ & $48$\\
$5$ && $25$ & $8$ & $3$ & $2$ & $100$\\
$6$ && $36$ & $10$ & $4$ & $2$ & $180$\\
$7$ && $49$ & $12$ & $5$ & $2$ & $294$\\
$8$ && $64$ & $14$ & $6$ & $2$ & $448$\\
$9$ && $81$ & $16$ & $7$ & $2$ & $648$\\
$10$ && $100$ & $18$ & $8$ & $2$ & $900$\\
\hline
\end{tabular}
\quad $\Rightarrow$ \quad
\begin{tabular}{cccc}
\multicolumn{4}{c}{Vectors of lattice}\\
\hline
$d$ & $n$ & $s$ & $t$\\
\hline
$10$ & $36$ & $10$ & $3$\\
$33$ & $96$ & $10$ & $3$\\
$76$ & $200$ & $10$ & $3$\\
$145$ & $360$ & $10$ & $3$\\
$246$ & $588$ & $10$ & $3$\\
$385$ & $896$ & $10$ & $3$\\
$568$ & $1296$ & $10$ & $3$\\
$801$ & $1800$ & $10$ & $3$\\
\hline
\end{tabular}
\end{center}\quad
\caption{Square lattice graph $L_2(m)$}\vspace{-0.27in}
\end{table}

\noindent
For every integer $3 \leqslant m \leqslant 10$, by numerical calculation, we can prove that we have the distance set $A(X) = \{ -1, \: \pm (m+1) / (2 f(m)), \: \pm m / (2 f(m)), \: \pm 2 / (2 f(m)), \: \pm 1/(2 f(m)), \: 0 \}$ and that the configuration of the vectors is $(d, n, s, t) = ((m-1) f(m), 2 m^2 (m-1), 10, 3)$, where we denote $f(m) := m^2 - m - 1$. Furthermore, for each $x \in X$, we have the following table:
\begin{center}
\begin{tabular}{c|cccccc}
$a$ & $-1$ & $\pm (m+1) / (2 f(m))$ & $\pm m / (2 f(m))$ & $\pm 2 / (2 f(m))$ & $\pm 1 / (2 f(m))$ & $0$\\
\hline
$|A_x(X, a)|$ & $1$ & $2 (m-2)$ & $2 (m-1)$ & $m-1$ & $2 (m-1) (m-2)$ & $2 (m-2) f(m)$\\
\end{tabular}
\end{center}

\quad\\

On the other hand, the complement graph of the square lattice graph $L_2(m)$ is a strongly regular graph with parameters $(v, k, \lambda, \mu) = (m^2, (m-1)^2, (m-2)^2, (m-1)(m-2))$. Moreover, if $m = 2$, then the constructed crystal lattice degenerate into the $0$-dimensional crystal lattice. Thus, we assume that $m \geqslant 3$.

For the $(d, n, s, t)$-configurations of the vectors of the constructed crystal lattices, we have the following table:

\begin{table}[h]
\begin{center}
\begin{tabular}{ccccccc}
\multicolumn{7}{c}{$\overline{L_2(m)}$}\\
\hline
$m$ & & $v$ & $k$ & $\lambda$ & $\mu$ & $e$\\
\hline
$3$ && $9$ & $4$ & $1$ & $2$ & $18$\\
$4$ && $16$ & $9$ & $4$ & $6$ & $72$\\
$5$ && $25$ & $16$ & $9$ & $12$ & $200$\\
$6$ && $36$ & $25$ & $16$ & $20$ & $300$\\
$7$ && $49$ & $36$ & $25$ & $30$ & $882$\\
\hline
\end{tabular}
\quad $\Rightarrow$ \quad
\begin{tabular}{cccc}
\multicolumn{4}{c}{Vectors of lattice}\\
\hline
$d$ & $n$ & $s$ & $t$\\
\hline
$10$ & $36$ & $10$ & $3$\\
$57$ & $144$ & $10$ & $3$\\
$176$ & $400$ & $10$ & $3$\\
$415$ & $900$ & $10$ & $3$\\
$834$ & $1764$ & $10$ & $3$\\
\hline
\end{tabular}
\end{center}\quad
\caption{Complement of the square lattice graph $\overline{L_2(m)}$}\vspace{-0.27in}
\end{table}

\noindent
Here, the graph $L_2(3)$ and its complement $\overline{L_2(3)}$ are isomorphic to each other.

\newpage

\subsection{Disjoint unions of the complete graphs $r \cdot K_m$ and complete multipartite graphs $K_{m, \ldots , m}$}\label{sec-ukn_cm}
The disjoint union of the complete graphs $r \cdot K_m$ is a strongly regular graph with parameters $(v, k, \lambda, \mu) = (m r, m-1, m-2, 0)$, which is vertex and edge-transitive. When the vectors of crystal lattice constructed from the complete graph $K_m$ has $(d, n, s, t)$-configuration, then the vectors of crystal lattice constructed from the disjoint union of the complete graph $r \cdot K_m$ has $(r d, r n, s, t)$-configuration. We can refer the result in Section \ref{sec-kn}.\\

On the other hand, the complement graph of the disjoint union of the complete graphs $r \cdot K_m$ is the complete (regular) multipartite graph $K_{m, \ldots , m}$, and it is a strongly regular graph with parameters $(v, k, \lambda, \mu) = (m r, m (r-1), m (r-2), m (r-1))$. In particular, the complement graphs of $r \cdot K_2$ are called the cocktail party graphs (See Section \ref{sec-cp}). Furthermore, we can refer Section \ref{sec-bip} for bipartite graphs.

For the $(d, n, s, t)$-configurations of the vectors of the constructed crystal lattices, we have the following table:

\begin{center}
\begin{tabular}{cccccccc}
\multicolumn{8}{c}{Complete $r$-partite graph}\\
\hline
$r$ & $m$ & \ & $v$ & $k$ & $\lambda$ & $\mu$ & $e$\\
\hline
$3$ & $2$ && $6$ & $4$ & $2$ & $4$ & $12$\\
& $3$ && $9$ & $6$ & $3$ & $6$ & $27$\\
& $4$ && $12$ & $8$ & $4$ & $8$ & $48$\\
& $5$ && $15$ & $10$ & $5$ & $10$ & $75$\\
& $6$ && $18$ & $12$ & $6$ & $12$ & $108$\\
& $7$ && $21$ & $14$ & $7$ & $14$ & $147$\\
& $8$ && $24$ & $16$ & $8$ & $16$ & $192$\\
& $9$ && $27$ & $18$ & $9$ & $18$ & $243$\\
& $10$ && $30$ & $20$ & $10$ & $20$ & $300$\\
& $11$ && $33$ & $22$ & $11$ & $22$ & $363$\\
& $12$ && $36$ & $24$ & $12$ & $24$ & $438$\\
& $13$ && $39$ & $26$ & $13$ & $26$ & $507$\\
\hline
\end{tabular}
\quad $\Rightarrow$ \quad
\begin{tabular}{cccc}
\multicolumn{4}{c}{Vectors of lattice}\\
\hline
$d$ & $n$ & $s$ & $t$\\
\hline
$7$ & $24$ & $9$ & $3$\\
$19$ & $54$ & $9$ & $3$\\
$37$ & $96$ & $9$ & $3$\\
$61$ & $150$ & $9$ & $3$\\
$91$ & $216$ & $9$ & $3$\\
$127$ & $294$ & $9$ & $3$\\
$169$ & $384$ & $9$ & $3$\\
$217$ & $486$ & $9$ & $3$\\
$271$ & $600$ & $9$ & $3$\\
$331$ & $726$ & $9$ & $3$\\
$397$ & $864$ & $9$ & $3$\\
$469$ & $1014$ & $9$ & $3$\\
\hline
\end{tabular}\vspace{0.1in}

\begin{tabular}{cccccccc}
\hline
$r$ & $m$ & \ & $v$ & $k$ & $\lambda$ & $\mu$ & $e$\\
\hline
$4$ & $2$ && $8$ & $6$ & $4$ & $6$ & $24$\\
& $3$ && $12$ & $9$ & $6$ & $9$ & $54$\\
& $4$ && $16$ & $12$ & $8$ & $12$ & $96$\\
& $5$ && $20$ & $15$ & $10$ & $15$ & $150$\\
& $6$ && $24$ & $18$ & $12$ & $18$ & $216$\\
& $7$ && $28$ & $21$ & $14$ & $21$ & $294$\\
& $8$ && $32$ & $24$ & $16$ & $24$ & $384$\\
& $9$ && $36$ & $27$ & $18$ & $27$ & $486$\\
\hline
\end{tabular}
\quad $\Rightarrow$ \quad
\begin{tabular}{cccc}
\hline
$d$ & $n$ & $s$ & $t$\\
\hline
$17$ & $48$ & $10$ & $3$\\
$43$ & $108$ & $10$ & $3$\\
$81$ & $192$ & $10$ & $3$\\
$131$ & $300$ & $10$ & $3$\\
$193$ & $432$ & $10$ & $3$\\
$267$ & $588$ & $10$ & $3$\\
$353$ & $768$ & $10$ & $3$\\
$451$ & $972$ & $10$ & $3$\\
\hline
\end{tabular}\vspace{0.1in}

\begin{tabular}{cccccccc}
\hline
$r$ & $m$ & \ & $v$ & $k$ & $\lambda$ & $\mu$ & $e$\\
\hline
$5$ & $2$ && $10$ & $8$ & $6$ & $8$ & $40$\\
& $3$ && $15$ & $12$ & $9$ & $12$ & $90$\\
& $4$ && $20$ & $16$ & $12$ & $16$ & $160$\\
& $5$ && $25$ & $20$ & $15$ & $20$ & $250$\\
& $6$ && $30$ & $24$ & $18$ & $24$ & $360$\\
& $7$ && $35$ & $28$ & $21$ & $28$ & $490$\\
\hline
\end{tabular}
\quad $\Rightarrow$ \quad
\begin{tabular}{cccc}
\hline
$d$ & $n$ & $s$ & $t$\\
\hline
$31$ & $80$ & $10$ & $3$\\
$76$ & $180$ & $10$ & $3$\\
$141$ & $320$ & $10$ & $3$\\
$226$ & $500$ & $10$ & $3$\\
$331$ & $720$ & $10$ & $3$\\
$456$ & $980$ & $10$ & $3$\\
\hline
\end{tabular}
\end{center}

\newpage

\begin{table}[h]
\begin{center}
\begin{tabular}{cccccccc}
\multicolumn{8}{c}{Complete $r$-partite graph}\\
\hline
$r$ & $m$ & \ & $v$ & $k$ & $\lambda$ & $\mu$ & $e$\\
\hline
$6$ & $2$ && $12$ & $10$ & $8$ & $10$ & $60$\\
& $3$ && $18$ & $15$ & $12$ & $15$ & $135$\\
& $4$ && $24$ & $20$ & $16$ & $20$ & $240$\\
& $5$ && $30$ & $25$ & $20$ & $25$ & $375$\\
\hline
\end{tabular}
\quad $\Rightarrow$ \quad
\begin{tabular}{cccc}
\multicolumn{4}{c}{Vectors of lattice}\\
\hline
$d$ & $n$ & $s$ & $t$\\
\hline
$49$ & $120$ & $10$ & $3$\\
$118$ & $270$ & $10$ & $3$\\
$217$ & $480$ & $10$ & $3$\\
$346$ & $750$ & $10$ & $3$\\
\hline
\end{tabular}\vspace{0.1in}

\begin{tabular}{cccccccc}
\hline
$r$ & $m$ & \ & $v$ & $k$ & $\lambda$ & $\mu$ & $e$\\
\hline
$7$ & $2$ && $14$ & $12$ & $10$ & $12$ & $84$\\
& $3$ && $21$ & $18$ & $15$ & $18$ & $189$\\
& $4$ && $28$ & $24$ & $20$ & $24$ & $336$\\
& $5$ && $35$ & $30$ & $25$ & $30$ & $525$\\
\hline
\end{tabular}
\quad $\Rightarrow$ \quad
\begin{tabular}{cccc}
\hline
$d$ & $n$ & $s$ & $t$\\
\hline
$71$ & $168$ & $10$ & $3$\\
$169$ & $378$ & $10$ & $3$\\
$309$ & $672$ & $10$ & $3$\\
$491$ & $1050$ & $10$ & $3$\\
\hline
\end{tabular}\vspace{0.1in}

\begin{tabular}{cccccccc}
\hline
$r$ & $m$ & \ & $v$ & $k$ & $\lambda$ & $\mu$ & $e$\\
\hline
$8$ & $2$ && $16$ & $14$ & $12$ & $14$ & $112$\\
& $3$ && $24$ & $21$ & $18$ & $21$ & $252$\\
& $4$ && $32$ & $28$ & $24$ & $28$ & $448$\\
\hline
\end{tabular}
\quad $\Rightarrow$ \quad
\begin{tabular}{cccc}
\hline
$d$ & $n$ & $s$ & $t$\\
\hline
$97$ & $224$ & $10$ & $3$\\
$229$ & $504$ & $10$ & $3$\\
$417$ & $896$ & $10$ & $3$\\
\hline
\end{tabular}\vspace{0.1in}

\begin{tabular}{cccccccc}
\hline
$r$ & $m$ & \ & $v$ & $k$ & $\lambda$ & $\mu$ & $e$\\
\hline
$9$ & $2$ && $18$ & $16$ & $14$ & $16$ & $144$\\
& $3$ && $27$ & $24$ & $21$ & $24$ & $324$\\
\hline
\end{tabular}
\quad $\Rightarrow$ \quad
\begin{tabular}{cccc}
\hline
$d$ & $n$ & $s$ & $t$\\
\hline
$127$ & $288$ & $10$ & $3$\\
$298$ & $648$ & $10$ & $3$\\
\hline
\end{tabular}\vspace{0.1in}

\begin{tabular}{cccccccc}
\hline
$r$ & $m$ & \ & $v$ & $k$ & $\lambda$ & $\mu$ & $e$\\
\hline
$10$ & $2$ && $20$ & $18$ & $16$ & $18$ & $180$\\
& $3$ && $30$ & $27$ & $24$ & $27$ & $405$\\
\hline
\end{tabular}
\quad $\Rightarrow$ \quad
\begin{tabular}{cccc}
\hline
$d$ & $n$ & $s$ & $t$\\
\hline
$161$ & $360$ & $10$ & $3$\\
$376$ & $810$ & $10$ & $3$\\
\hline
\end{tabular}\vspace{0.1in}

\begin{tabular}{cccccccc}
\hline
$r$ & $m$ & \ & $v$ & $k$ & $\lambda$ & $\mu$ & $e$\\
\hline
$11$ & $2$ && $22$ & $20$ & $18$ & $20$ & $220$\\
& $3$ && $33$ & $30$ & $27$ & $30$ & $495$\\
\hline
\end{tabular}
\quad $\Rightarrow$ \quad
\begin{tabular}{cccc}
\hline
$d$ & $n$ & $s$ & $t$\\
\hline
$199$ & $440$ & $10$ & $3$\\
$463$ & $990$ & $10$ & $3$\\
\hline
\end{tabular}\vspace{0.1in}

\begin{tabular}{cccccccc}
\hline
$r$ & $m$ & \ & $v$ & $k$ & $\lambda$ & $\mu$ & $e$\\
\hline
$12$ & $2$ && $24$ & $22$ & $20$ & $22$ & $264$\\
\hline
\end{tabular}
\quad $\Rightarrow$ \quad
\begin{tabular}{cccc}
\hline
$d$ & $n$ & $s$ & $t$\\
\hline
$241$ & $528$ & $10$ & $3$\\
\hline
\end{tabular}\vspace{0.1in}

\begin{tabular}{cccccccc}
\hline
$r$ & $m$ & \ & $v$ & $k$ & $\lambda$ & $\mu$ & $e$\\
\hline
$13$ & $2$ && $26$ & $24$ & $22$ & $24$ & $312$\\
\hline
\end{tabular}
\quad $\Rightarrow$ \quad
\begin{tabular}{cccc}
\hline
$d$ & $n$ & $s$ & $t$\\
\hline
$287$ & $624$ & $10$ & $3$\\
\hline
\end{tabular}
\end{center}\quad
\caption{Complete $r$-partite graph}\vspace{-0.27in}
\end{table}

\noindent
Let us denote $d(m, r) := m^2 r (r-1) / 2 - m r + 1$. For every integer $m$ and $r$ such that $d(m, r) \leqslant 500$, by numerical calculation, we can prove that we have the distance set $A(X) = \{ -1, \: \pm (m r - 1)/(2 d(m, r)), \: \pm (m r - 2)/(2 d(m, r)), \: \pm 2/(2 d(m, r)), \: \pm 1/(2 d(m, r)), \: 0 \}$ (If $r=3$, $A(X)$ does not include $0$.) and that the configuration of the vectors is $(d, n, s, t) = (d(m, r), m^2 r (r-1) / 2, 10, 3)$ (If $r=3$, then $s=9$). Furthermore, for each $x \in X$, we have the following table:
\begin{center}
\begin{tabular}{c|cccc}
$a$ & $-1$ & $\pm (m r - 1)/(2 d(m, r))$ & $\pm (m r - 2)/(2 d(m, r))$\\
\cline{1-4}
$|A_x(X, a)|$ & $1$ & $m (r-2)$ & $2 (m-1)$\\
\multicolumn{2}{c}{\quad}& $\pm 2/(2 d(m, r))$ & $\pm 1/(2 d(m, r))$ & $0$\\
\cline{3-5}
\multicolumn{2}{c}{\quad}& $(m-1)^2$ & $2 m (m-1) (r-2)$ & $m^2 (r-2) (r-3)$\\
\end{tabular}
\end{center}

\newpage

\subsection{Paley graph $P(q)$}\label{sec-pq}
Let $q$ be a prime power such that $q \equiv 1 \pmod{4}$, then the Paley graph $P(q)$ is a strongly regular graph with parameters $(v, k, \lambda, \mu) = (q, (q-1)/2, (q-5)/4, (q-1)/4)$, which is vertex and edge-transitive. Let $V = \mathbb{F}_q$ which is the finite field. Then, we define $E$ as a set of pairs $(x, y)$ of distinct elements of $V$ such that $x-y$ is a square. Now, we can define the Paley graph $P(q)$ as the graph with $V$ and $E$ as a set of vertices and edges, respectively. Moreover, if $q=5$, then the constructed crystal lattice degenerate into the $1$-dimensional standard crystal lattice. Thus, we assume that $q \geqslant 9$.

For the $(d, n, s, t)$-configurations of the vectors of the constructed crystal lattices, we have the following table:

\begin{table}[h]
\begin{center}
\begin{tabular}{ccccccc}
\multicolumn{7}{c}{Paley graph $P(q)$}\\
\hline
$q$ & & $v$ & $k$ & $\lambda$ & $\mu$ & $e$\\
\hline
$9$ && $9$ & $4$ & $1$ & $2$ & $18$\\
$13$ && $13$ & $6$ & $2$ & $3$ & $39$\\
$17$ && $17$ & $8$ & $3$ & $4$ & $68$\\
$25$ && $25$ & $12$ & $5$ & $6$ & $150$\\
$29$ && $29$ & $14$ & $6$ & $7$ & $203$\\
$37$ && $37$ & $18$ & $8$ & $9$ & $333$\\
$41$ && $41$ & $20$ & $9$ & $10$ & $410$\\
$49$ && $49$ & $24$ & $11$ & $12$ & $588$\\
\hline
\end{tabular}
\quad $\Rightarrow$ \quad
\begin{tabular}{cccc}
\multicolumn{4}{c}{Vectors of lattice}\\
\hline
$d$ & $n$ & $s$ & $t$\\
\hline
$10$ & $36$ & $10$ & $3$\\
$27$ & $78$ & $10$ & $3$\\
$52$ & $136$ & $10$ & $3$\\
$126$ & $300$ & $10$ & $3$\\
$175$ & $406$ & $10$ & $3$\\
$297$ & $666$ & $10$ & $3$\\
$370$ & $820$ & $10$ & $3$\\
$540$ & $1176$ & $10$ & $3$\\
\hline
\end{tabular}
\end{center}\quad
\caption{Paley graph $P(q)$}\vspace{-0.27in}
\end{table}

\noindent
For every integer $9 \leqslant q \leqslant 49$, by numerical calculation, we can prove that we have the distance set $A(X) = \{ -1, \: \pm (q-1) / (2 d(q)), \: \pm (q-3) / (2 d(q)), \: \pm 2 / d(q), \: \pm 1 / d(q), \: 0 \}$ and that the configuration of the vectors is $(d, n, s, t) = (d(q), q (q-1) / 2, 10, 3)$, where we denote $d(q) := (q-1) (q-4) / 4$.\\

On the other hand, the complement graph of the Paley graph $P(q)$ is isomorphic to original graph. Moreover, Paley graph $P(9)$ is isomorphic to the square lattice graph $L_2(3)$.

\newpage

\subsection{Other strongly regular graphs}
In this section, we consider the other strongly regular graphs. There are some graphs which have the same parameters $(v, k, \lambda, \mu)$ but are not isomorphic. For such graphs, we refer to the electric data on the website \cite{Se2}. We consider details about them in Section \ref{sec-dif-sr}.

\begin{table}[h]
\begin{center}
\begin{tabular}{lccccc}
\multicolumn{6}{c}{Strongly regular graphs}\\
\hline
& $v$ & $k$ & $\lambda$ & $\mu$ & $e$\\
\hline
Clebsch & $16$ & $5$ & $0$ & $2$ & $40$\\
Shrikhande & $16$ & $6$ & $2$ & $2$ & $48$\\
$\overline{\text{Shrikhande}}$ & $16$ & $9$ & $4$ & $6$ & $72$\\
$\overline{\text{Clebsch}}$ & $16$ & $10$ & $6$ & $6$ & $80$\\
Paulus $[15]$ & $25$ & $12$ & $5$ & $6$ & $150$\\
Paulus $[10]$ & $26$ & $10$ & $3$ & $4$ & $130$\\
$\overline{\text{Paulus}}$ $[10]$ & $26$ & $15$ & $8$ & $9$ & $195$\\
Schl\"{a}fli & $27$ & $10$ & $1$ & $5$ & $135$\\
$\overline{\text{Schl\"{a}fli}}$ & $27$ & $16$ & $10$ & $6$ & $216$\\
Chang $[3]$ & $28$ & $12$ & $6$ & $4$ & $168$\\
$\overline{\text{Chang}}$ $[3]$ & $28$ & $15$ & $6$ & $10$ & $168$\\
$[41]$ & $29$ & $14$ & $6$ & $7$ & $203$\\
$[28]$ & $40$ & $12$ & $2$ & $4$ & $240$\\
$[28]$ & $40$ & $27$ & $18$ & $18$ & $540$\\
{\small Hoffman-Singleton} & $50$ & $7$ & $0$ & $1$ & $175$\\
Gewirtz & $56$ & $10$ & $0$ & $2$ & $280$\\
$M_{22}$ & $77$ & $16$ & $0$ & $4$ & $616$\\
\hline
\end{tabular}
\ $\Rightarrow$ \
\begin{tabular}{ccccc}
\multicolumn{5}{c}{Vectors of lattice}\\
\hline
$d$ & $n$ & $s$ & $t$ & $A(X) \setminus \{ -1, 0\}$\\
\hline
$25$ & $80$ & $8$ & $3$ & $\{ \pm \frac{1}{4}, \pm \frac{1}{10}, \pm \frac{1}{20} \}$\\
$33$ & $96$ & $10$ & $3$ & $\{ \pm \frac{5}{22}, \pm \frac{2}{11}, \pm \frac{1}{11}, \pm \frac{1}{22} \}$\\
$57$ & $144$ & $10$ & $3$ & $\{ \pm \frac{5}{38}, \pm \frac{9}{76}, \pm \frac{1}{38}, \pm \frac{1}{76}\}$\\
$65$ & $160$ & $10$ & $3$ & $\{ \pm \frac{3}{26}, \pm \frac{4}{39}, \pm \frac{1}{39}, \pm \frac{1}{78} \}$\\
$126$ & $300$ & $10$ & $3$ & $\{ \pm \frac{2}{21}, \pm \frac{11}{126}, \pm \frac{1}{63}, \pm \frac{1}{84} \}$\\
$105$ & $260$ & $10$ & $3$ & $\{ \pm \frac{5}{42}, \pm \frac{3}{28}, \pm \frac{1}{42}, \pm \frac{1}{126} \}$\\
$170$ & $390$ & $10$ & $3$ & $\{ \pm \frac{5}{68}, \pm \frac{7}{102}, \pm \frac{1}{102}, \pm \frac{1}{204} \}$\\
$109$ & $270$ & $10$ & $3$ & $\{ \pm \frac{13}{109}, \pm \frac{12}{109}, \pm \frac{2}{109}, \pm \frac{1}{109} \}$\\
$190$ & $432$ & $10$ & $3$ & $\{ \pm \frac{13}{190}, \pm \frac{6}{95}, \pm \frac{1}{95}, \pm \frac{1}{190} \}$\\
$141$ & $336$ & $10$ & $3$ & $\{ \pm \frac{9}{94}, \pm \frac{4}{47}, \pm \frac{1}{47}, \pm \frac{1}{94} \}$\\
$183$ & $420$ & $10$ & $3$ & $\{ \pm \frac{9}{122}, \pm \frac{17}{244}, \pm \frac{1}{122}, \pm \frac{1}{244} \}$\\
$175$ & $406$ & $10$ & $3$ & $\{ \pm \frac{2}{25}, \pm \frac{13}{175}, \pm \frac{2}{175}, \pm \frac{1}{175} \}$\\
$201$ & $480$ & $10$ & $3$ & $\{ \pm \frac{13}{134}, \pm \frac{6}{67}, \pm \frac{1}{67}, \pm \frac{1}{134} \}$\\
$501$ & $1080$ & $10$ & $3$ & $\{ \pm \frac{13}{334}, \pm \frac{25}{668}, \pm \frac{1}{334}, \pm \frac{1}{668} \}$\\
$126$ & $350$ & $6$ & $3$ & $\{ \pm \frac{1}{6}, \pm \frac{1}{36} \}$\\
$225$ & $560$ & $8$ & $3$ & $\{ \pm \frac{1}{9}, \pm \frac{1}{45}, \pm \frac{1}{90} \}$\\
$540$ & $1232$ & $8$ & $3$ & $\{ \frac{1}{15}, \pm \frac{1}{135}, \pm \frac{1}{270} \}$\\
\hline
\end{tabular}
\end{center}\quad
\caption{Other strongly regular graph}\vspace{-0.27in}
\end{table}

\noindent
Here, we denote the number in `$[ \ ]$' by the number of strongly regular graphs with the same parameters, and denote `$\overline{\text{Name}}$' by the complement of the named graph. For the above cases, the sets of vectors of constructed crystal lattices from the strongly regular graphs with the same parameters are not isomorphic, but has the same distance set and the same $(d, n, s, t)$-configuration.

The Shrikhande graph and the Chang graphs has the same parameters as the square lattice graph $L_2(4)$ and the triangular graph $T(8)$, respectively. However, they are not isomorphic to each other.

Furthermore, the above strongly regular graphs with $25$ and $29$ vertices include the graph which is isomorphic to the Paley graph $P(25)$ and $P(29)$, respectively.

\newpage

\subsection{List of strongly regular graphs with at most $30$ vertices}\quad

\begin{table}[h]
\begin{center}
\begin{tabular}{ccccl}
\hline
$v$ & $k$ & $\lambda$ & $\mu$\\
\hline
$1$ & $0$ & $0$ & $0$ & $K_1$\\
$2$ & $1$ & $0$ & $0$ & $K_2$\\
$3$ & $2$ & $1$ & $0$ & $K_3$\\
$4$ & $2$ & $0$ & $2$ & $C_4$\\
& $3$ & $2$ & $0$ & $K_4$\\
$5$ & $2$ & $0$ & $1$ & $C_5$\\
& $4$ & $3$ & $0$ & $K_5$\\
$6$ & $3$ & $0$ & $3$ & $K_{3, 3}$\\
& $4$ & $2$ & $4$ & $CP(3)$\\
& $5$ & $4$ & $0$ & $K_6$\\
$7$ & $6$ & $5$ & $0$ & $K_7$\\
$8$ & $4$ & $0$ & $4$ & $K_{4, 4}$\\
& $6$ & $4$ & $6$ & $CP(4)$\\
& $7$ & $6$ & $0$ & $K_8$\\
$9$ & $4$ & $1$ & $2$ & $L_2(3)$\\
& $6$ & $3$ & $6$ & $K_{3, 3, 3}$\\
& $8$ & $7$ & $0$ & $K_9$\\
$10$ & $3$ & $0$ & $1$ & Petersen\\
& $5$ & $0$ & $5$ & $K_{5, 5}$\\
& $6$ & $3$ & $4$ & $T(5)$\\
& $8$ & $6$ & $8$ & $CP(5)$\\
& $9$ & $8$ & $0$ & $K_{10}$\\
$11$ & $10$ & $9$ & $0$ & $K_{11}$\\
$12$ & $6$ & $0$ & $6$ & $K_{6, 6}$\\
& $8$ & $4$ & $8$ & $K_{4, 4, 4}$\\
& $9$ & $6$ & $9$ & $K_{3, 3, 3, 3}$\\
& $10$ & $8$ & $10$ & $CP(6)$\\
& $11$ & $10$ & $0$ & $K_{12}$\\
$13$ & $6$ & $2$ & $3$ & $P(13)$\\
& $12$ & $11$ & $0$ & $K_{13}$\\
$14$ & $7$ & $0$ & $7$ & $K_{7, 7}$\\
& $12$ & $10$ & $12$ & $CP(7)$\\
& $13$ & $12$ & $0$ & $K_{14}$\\
$15$ & $6$ & $1$ & $3$ & $\overline{T(6)}$\\
& $8$ & $4$ & $4$ & $T(6)$\\
& $10$ & $5$ & $10$ & $K_{5, 5, 5}$\\
& $12$ & $9$ & $12$ & $K_{3, \ldots , 3}$\\
& $14$ & $13$ & $0$ & $K_{15}$\\
\hline\\
\end{tabular}
\quad
\begin{tabular}{ccccl}
\hline
$v$ & $k$ & $\lambda$ & $\mu$\\
\hline
$16$ & $5$ & $0$ & $2$ & Clebsch\\
& $6$ & $2$ & $2$ & $L_2(4)$\\
&&&& Shrikhande\\
& $8$ & $0$ & $8$ & $K_{8, 8}$\\
& $9$ & $4$ & $6$ & $\overline{L_2(4)}$\\
&&&& $\overline{\text{Shrikhande}}$\\
& $10$ & $6$ & $6$ & $\overline{\text{Clebsch}}$\\
& $12$ & $8$ & $12$ & $K_{4, 4, 4, 4}$\\
& $14$ & $12$ & $14$ & $CP(8)$\\
& $15$ & $14$ & $0$ & $K_{16}$\\
$17$ & $8$ & $3$ & $4$ & $P(17)$\\
& $16$ & $15$ & $0$ & $K_{17}$\\
$18$ & $9$ & $0$ & $9$ & $K_{9, 9}$\\
& $12$ & $6$ & $12$ & $K_{6, 6, 6}$\\
& $15$ & $12$ & $15$ & $K_{3, \ldots , 3}$\\
& $16$ & $14$ & $16$ & $CP(9)$\\
& $17$ & $16$ & $0$ & $K_{18}$\\
$19$ & $18$ & $17$ & $0$ & $K_{19}$\\
$20$ & $10$ & $0$ & $10$ & $K_{10, 10}$\\
& $15$ & $10$ & $15$ & $K_{5, 5, 5, 5}$\\
& $16$ & $12$ & $16$ & $K_{4, \ldots , 4}$\\
& $18$ & $16$ & $18$ & $CP(10)$\\
& $19$ & $18$ & $0$ & $K_{20}$\\
$21$ & $10$ & $3$ & $6$ & $\overline{T(7)}$\\
& $10$ & $5$ & $4$ & $T(7)$\\
& $14$ & $7$ & $14$ & $K_{7, 7, 7}$\\
& $18$ & $15$ & $18$ & $K_{3, \ldots , 3}$\\
& $20$ & $19$ & $0$ & $K_{21}$\\
$22$ & $11$ & $0$ & $11$ & $K_{11, 11}$\\
& $20$ & $18$ & $20$ & $CP(11)$\\
& $21$ & $20$ & $0$ & $K_{22}$\\
$23$ & $22$ & $21$ & $0$ & $K_{23}$\\
$24$ & $12$ & $0$ & $12$ & $K_{12, 12}$\\
& $16$ & $8$ & $16$ & $K_{8, 8, 8}$\\
& $18$ & $12$ & $18$ & $K_{6, 6, 6, 6}$\\
& $20$ & $16$ & $20$ & $K_{4, \ldots , 4}$\\
& $21$ & $18$ & $21$ & $K_{3, \ldots , 3}$\\
& $22$ & $20$ & $22$ & $CP(12)$\\
& $23$ & $22$ & $0$ & $K_{24}$\\
\hline
\end{tabular}
\quad
\begin{tabular}{ccccl}
\hline
$v$ & $k$ & $\lambda$ & $\mu$\\
\hline
$25$ & $8$ & $3$ & $2$ & $L_2(5)$\\
& $12$ & $5$ & $6$ & {\scriptsize Paulus $[15]$}\\
&&&& {\scriptsize iwith $P(25)$j}\\
& $16$ & $9$ & $12$ & $\overline{L_2(5)}$\\
& $20$ & $15$ & $20$ & $K_{5, \ldots , 5}$\\
& $24$ & $23$ & $0$ & $K_{25}$\\
$26$ & $10$ & $3$ & $4$ & {\scriptsize Paulus $[10]$}\\
& $13$ & $0$ & $13$ & $K_{13, 13}$\\
& $15$ & $8$ & $9$ & {\scriptsize $\overline{\text{Paulus}}$ $[15]$}\\
& $24$ & $22$ & $24$ & $CP(13)$\\
& $25$ & $24$ & $0$ & $K_{26}$\\
$27$ & $10$ & $1$ & $5$ & Schl\"{a}fli\\
& $16$ & $10$ & $6$ & $\overline{\text{Schl\"{a}fli}}$\\
& $18$ & $9$ & $18$ & $K_{9, 9, 9}$\\
& $24$ & $21$ & $24$ & $K_{3, \ldots , 3}$\\
& $26$ & $25$ & $0$ & $K_{27}$\\
$28$ & $12$ & $6$ & $4$ & $T(8)$\\
&&&& Chang $[3]$\\
& $14$ & $0$ & $14$ & $K_{14, 14}$\\
& $15$ & $6$ & $10$ & $\overline{T(8)}$\\
&&&& $\overline{\text{Chang}}$ $[3]$\\
& $21$ & $14$ & $21$ & $K_{7, 7, 7, 7}$\\
& $24$ & $20$ & $24$ & $K_{4, \ldots , 4}$\\
& $26$ & $24$ & $26$ & $CP(14)$\\
& $27$ & $26$ & $0$ & $K_{28}$\\
$29$ & $14$ & $6$ & $7$ & $P(29)$\\
&&&& otherwise $[40]$\\
& $28$ & $27$ & $0$ & $K_{29}$\\
$30$ & $15$ & $0$ & $15$ & $K_{15, 15}$\\
& $20$ & $10$ & $20$ & $K_{10, 10, 10}$\\
& $24$ & $18$ & $24$ & $K_{6, \ldots , 6}$\\
& $25$ & $20$ & $25$ & $K_{5, \ldots , 5}$\\
& $27$ & $24$ & $27$ & $K_{3, \ldots , 3}$\\
& $28$ & $26$ & $28$ & $CP(15)$\\
& $29$ & $28$ & $0$ & $K_{30}$\\
\hline\\ \\ \\ \\
\end{tabular}
\end{center}\quad
\caption{Strongly regular graphs with at most $30$ vertices}
\end{table}

\noindent
Here, we denote the number in `$[ \ ]$' by the number of strongly regular graphs with the same parameters.

\newpage

\subsection{On a difference among the vectors of crystal lattices from strongly regular graphs with the same parameters}\label{sec-dif-sr}

In the previous sections, the vectors of crystal lattices from strongly regular graphs with the same parameters $(v, k, \lambda, \mu)$have the same $(d, n, s, t)$-configuration and the same distance set. However, there can exist some difference in details. In this section, we want to introduce such examples.\\

\subsubsection{Square lattice graph $L_2(4)$ and Shrikhande graph}\label{sec-dif-l24}

Both the square lattice graph $L_2(4)$ and the Shrikhande graph have the same parameters $(v, k, \lambda, \mu) = (16, 6, 2, 2)$ and they are vertex and edge-transitive, and the vectors of crystal lattices from each graph has $(33, 96, 10, 3)$-configuration and the distance set $A(X) = \{ -1, \: \pm 5/22, \: \pm 2/11, \: \pm 1/11, \: \pm 1/22, \: 0\}$. However, for each $x \in X$ from each graph, we have the following table:
\begin{center}
\begin{tabular}{lccccccc}
\hline
& $a$ & $-1$ & $\pm 5/22$ & $\pm 2/11$ & $\pm 1/11$ & $\pm 1/22$ & $0$\\
\hline
$L_2(4)$ & $|A_x(X, a)|$ & $1$ & $4$ & $6$ & $3$ & $12$ & $44$\\
Shrikhande & $|A_x(X, a)|$ & $1$ & $4$ & $6$ & $1$ & $20$ & $32$\\
\hline
\end{tabular}
\end{center}

\quad\\

On the other hand, the complements of the two graphs have the same parameters $(v, k, \lambda, \mu) = (16, 9, 4, 6)$, and they are vertex-transitive but only the complement of $L_2(4)$ is edge-transitive. Then, the vectors of crystal lattices from each graph has $(57, 144, 10, 3)$-configuration and the distance set $A(X) = \{ -1, \: \pm 5/38, \: \pm 9/76, \: \pm 1/38, \: \pm 1/76, \: 0\}$. However, for the numbers $|A_x(X, a)|$ for each $x \in X$ from the two graphs, we have the following two types, namely type A and type B:
\begin{center}
\begin{tabular}{cccccccc}
\hline
& $a$ & $-1$ & $\pm 5/38$ & $\pm 9/76$ & $\pm 1/38$ & $\pm 1/76$ & $0$\\
\hline
'` & $|A_x(X, a)|$ & $1$ & $8$ & $8$ & $12$ & $24$ & $38$\\
B & $|A_x(X, a)|$ & $1$ & $8$ & $8$ & $10$ & $32$ & $26$\\
\hline
\end{tabular}
\end{center}
Then, all of the vectors of lattice from the complement of $L_2(4)$ are of type A, and $48$ vectors of lattice from the complement of the Shrikhande graph are of type A and the other $96$ vectors are of type B.\\

\subsubsection{Triangular graph $T(8)$ and Chang graphs}\label{sec-dif-t8}

There are just four graphs which have the same parameters as the triangular graph $T(8)$ (See \cite{Cl}). Let us denote the three Chang graphs by $G_1$, $G_2$, $G_3$ for convenience. The triangular graph $T(8)$ and the three Chang graphs have the same parameters $(v, k, \lambda, \mu) = (28, 12, 6, 4)$, and the vectors of crystal lattices from each graph has $(141, 336, 10, 3)$-configuration and the distance set $A(X) = \{ -1, \: \pm 9/94, \: \pm 4/47, \: \pm 1/47, \: \pm 1/94, \: 0\}$. However, for the numbers $|A_x(X, a)|$ for each $x \in X$ from the four graphs, we have the following three types:
\begin{center}
\begin{tabular}{cccccccc}
\hline
& $a$ & $-1$ & $\pm 9/94$ & $\pm 4/47$ & $\pm 1/47$ & $\pm 1/94$ & $0$\\
\hline
A & $|A_x(X, a)|$ & $1$ & $12$ & $10$ & $5$ & $60$ & $160$\\
B & $|A_x(X, a)|$ & $1$ & $12$ & $10$ & $3$ & $68$ & $148$\\
C & $|A_x(X, a)|$ & $1$ & $12$ & $10$ & $1$ & $76$ & $136$\\
\hline
\end{tabular}
\end{center}
Then, all of the vectors of lattice from $T(8)$ is of type A, and the numbers of vectors of lattices from $G_1$, $G_2$, $G_3$ of each type is $($A, B, C$)$ $=$ $(48, 192, 96)$, $(36, 120, 180)$, $(0, 192, 144)$, respectively.\\

On the other hand, the complements of the four graphs have the same parameters $(v, k, \lambda, \mu) = (28, 15, 6, 10)$, and the vectors of crystal lattices from each graph has $(183, 420, 10, 3)$-configuration and the distance set $A(X) = \{ -1, \: \pm 9/122, \: \pm 17/244, \: \pm 1/122, \: \pm 1/244, \: 0\}$. However, for the numbers $|A_x(X, a)|$ for each $x \in X$ from the four graphs, we have the following three types:
\begin{center}
\begin{tabular}{cccccccc}
\hline
& $a$ & $-1$ & $\pm 9/122$ & $\pm 17/244$ & $\pm 1/122$ & $\pm 1/244$ & $0$\\
\hline
'` & $|A_x(X, a)|$ & $1$ & $12$ & $16$ & $48$ & $80$ & $106$\\
B & $|A_x(X, a)|$ & $1$ & $12$ & $16$ & $46$ & $88$ & $94$\\
C & $|A_x(X, a)|$ & $1$ & $12$ & $16$ & $44$ & $96$ & $82$\\
\hline
\end{tabular}
\end{center}
Then, all of the vectors of lattice from the complement of $T(8)$ is of type A, and the numbers of vectors of lattices from the complements of $G_1$, $G_2$, $G_3$ of each type is $($A, B, C$)$ $=$ $(132, 192, 96)$, $(0, 360, 60)$, $(36, 288, 96)$, respectively.

\newpage

\subsubsection{Paulus graphs with $25$ vertices}

A. J. L. Paulus classified the strongly regular graphs with the parameters $(v, k, \lambda, \mu) = (25, 12, 5, 6)$, then we have $15$ graphs up to isomorphic (See \cite{P}). One of the Paulus graphs with $25$ vertices is isomorphic to the Paley graph $P(25)$, and we can separate the other graphs into $7$ pairs of graphs which are complement to each other. Let us choose one of each $7$ pairs and denote by $G_1$, $G_2$, \ldots, $G_7$ for convenience. The vectors of crystal lattices from each graph has $(126, 300, 10, 3)$-configuration and the distance set $A(X) = \{ -1, \: \pm 2/21, \: \pm 11/126, \: \pm 1/63, \: \pm 1/126, \: 0\}$. However, for the numbers $|A_x(X, a)|$ for each $x \in X$ from the graphs, we have the following three types:
\begin{center}
\begin{tabular}{cccccccc}
\hline
& $a$ & $-1$ & $\pm 2/21$ & $\pm 11/126$ & $\pm 1/63$ & $\pm 1/126$ & $0$\\
\hline
A & $|A_x(X, a)|$ & $1$ & $10$ & $12$ & $18$ & $60$ & $98$\\
B & $|A_x(X, a)|$ & $1$ & $10$ & $12$ & $16$ & $68$ & $86$\\
C & $|A_x(X, a)|$ & $1$ & $10$ & $12$ & $14$ & $76$ & $74$\\
\hline
\end{tabular}
\end{center}
Then, for the numbers of vectors of lattices from the $15$ graphs $P(25)$, $G_1$, $\overline{G_1}$, \ldots, $G_7$, $\overline{G_7}$ of each type, we have the following table:
\begin{center}
\begin{tabular}{cccccccccc}
\hline
& $P(25)$ & $G_1, \overline{G_1}$ & $G_2, \overline{G_2}$ & $G_3, \overline{G_3}$ & $G_4, \overline{G_4}$ & $G_5, \overline{G_5}$ & $G_6, \overline{G_6}$ & \ $G_7$ \ & \ $\overline{G_7}$ \ \\
\hline
'` & -- & $18$ & $114$ & $154$ & $170$ & $172$ & $180$ & $48$ & $120$\\
B & $300$ & $240$ & $168$ & $136$ & $128$ & $124$ & $120$ & $252$ & $108$\\
C & -- & $42$ & $18$ & $10$ & $2$ & $4$ & -- & -- & $72$\\
\hline
\end{tabular}
\end{center}

\quad\\

\subsubsection{Paulus graphs with $26$ vertices}

A. J. L. Paulus also classified the strongly regular graphs with the parameters $(v, k, \lambda, \mu) = (26, 10, 3, 4)$, then we have $10$ graphs up to isomorphic (See \cite{P}). Let us denote the $10$ graphs by $G_1$, $G_2$, \ldots, $G_{10}$ for convenience. The vectors of crystal lattices from each graph has $(105, 260, 10, 3)$-configuration and the distance set $A(X) = \{ -1, \: \pm 5/42, \: \pm 3/28, \: \pm 1/42, \: \pm 1/84, \: 0\}$. However, for the numbers $|A_x(X, a)|$ for each $x \in X$ from the graphs, we have the following two types:
\begin{center}
\begin{tabular}{cccccccc}
\hline
& $a$ & $-1$ & $\pm 5/42$ & $\pm 3/28$ & $\pm 1/42$ & $\pm 1/84$ & $0$\\
\hline
A & $|A_x(X, a)|$ & $1$ & $6$ & $12$ & $14$ & $52$ & $90$\\
B & $|A_x(X, a)|$ & $1$ & $6$ & $12$ & $12$ & $60$ & $78$\\
\hline
\end{tabular}
\end{center}

On the other hand, the complements of the Paulus graphs with $26$ vertices have the same parameters $(v, k, \lambda, \mu) = (26, 15, 8, 9)$, and the vectors of crystal lattices from each graph has $(170, 390, 10, 3)$-configuration and the distance set $A(X) = \{ -1, \: \pm 5/68, \: \pm 7/102, \: \pm 1/102, \: \pm 1/204, \: 0\}$. However, for the numbers $|A_x(X, a)|$ for each $x \in X$ from the graphs, we have the following four types:
\begin{center}
\begin{tabular}{cccccccc}
\hline
& $a$ & $-1$ & $\pm 5/68$ & $\pm 7/102$ & $\pm 1/102$ & $\pm 1/204$ & $0$\\
\hline
C & $|A_x(X, a)|$ & $1$ & $16$ & $12$ & $24$ & $72$ & $140$\\
D & $|A_x(X, a)|$ & $1$ & $16$ & $12$ & $22$ & $80$ & $128$\\
E & $|A_x(X, a)|$ & $1$ & $16$ & $12$ & $20$ & $88$ & $116$\\
F & $|A_x(X, a)|$ & $1$ & $16$ & $12$ & $18$ & $96$ & $104$\\
\hline
\end{tabular}
\end{center}

Now, for the numbers of vectors of lattices from the $10$ graphs and their complements $G_1$, $\overline{G_1}$, \ldots, $G_{10}$, $\overline{G_{10}}$ of each type, we have the following table:
\begin{center}
\begin{tabular}{ccccccccccc}
\hline
& $G_1$ & $G_2$ & $G_3$ & $G_4$ & $G_5$ & $G_6$ & $G_7$ & $G_8$ & $G_9$ & $G_{10}$\\
\hline
'` & -- & -- & $12$ & $24$ & $24$ & $96$ & $120$ & $120$ & $144$ & $156$\\
B & $260$ & $260$ & $248$ & $236$ & $236$ & $164$ & $140$ & $140$ & $116$ & $104$\\
\hline
& $\overline{G_1}$ & $\overline{G_2}$ & $\overline{G_3}$ & $\overline{G_4}$ & $\overline{G_5}$ & $\overline{G_6}$ & $\overline{G_7}$ & $\overline{G_8}$ & $\overline{G_9}$ & $\overline{G_{10}}$\\
\hline
C & -- & -- & -- & -- & -- & -- & $60$ & $20$ & $20$ & --\\
D & $72$ & $76$ & $84$ & $90$ & $90$ & $96$ & -- & $96$ & $114$ & $156$\\
E & $246$ & $238$ & $234$ & $234$ & $234$ & $294$ & $330$ & $258$ & $246$ & $234$\\
F & $72$ & $76$ & $72$ & $66$ & $66$ & -- & -- & $16$ & $10$ & --\\
\hline
\end{tabular}
\end{center}

\newpage

\subsubsection{Strongly regular graphs with the parameters $(v, k, \lambda, \mu) = (29, 14, 6, 7)$}

We have $41$ graphs as the strongly regular graphs with the parameters $(v, k, \lambda, \mu) = (29, 14, 6, 7)$ up to isomorphic (classified by F. C. Bussemaker and E. Spence independently by computer search). One of these graphs is isomorphic to the Paley graph $P(29)$, and we can separate the other graphs into $20$ pairs of graphs which are complement to each other. Let us choose one of each $20$ pairs and denote by $G_1$, $G_2$, \ldots, $G_{20}$ for convenience. The vectors of crystal lattices from each graph has $(175, 406, 10, 3)$-configuration and the distance set $A(X) = \{ -1, \: \pm 2/25, \: \pm 13/175, \: \pm 2/175, \: \pm 1/175, \: 0\}$. However, for the numbers $|A_x(X, a)|$ for each $x \in X$ from the graphs, we have the following four types:
\begin{center}
\begin{tabular}{cccccccc}
\hline
& $a$ & $-1$ & $\pm 2/25$ & $\pm 13/175$ & $\pm 2/175$ & $\pm 1/175$ & $0$\\
\hline
A & $|A_x(X, a)|$ & $1$ & $12$ & $14$ & $26$ & $78$ & $144$\\
B & $|A_x(X, a)|$ & $1$ & $12$ & $14$ & $24$ & $86$ & $132$\\
C & $|A_x(X, a)|$ & $1$ & $12$ & $14$ & $22$ & $94$ & $120$\\
D & $|A_x(X, a)|$ & $1$ & $12$ & $14$ & $20$ & $102$ & $108$\\
\hline
\end{tabular}
\end{center}
Then, for the numbers of vectors of lattices from the $41$ graphs $P(29)$, $G_1$, $\overline{G_1}$, \ldots, $G_{20}$, $\overline{G_{20}}$ of each type, we have the following table:
{\small \begin{center}
\begin{tabular}{cccccccccccccccc}
\hline
& $P(29)$ & $G_1$ & $\overline{G_1}$ & $G_2$ & $\overline{G_2}$ & $G_3$ & $\overline{G_3}$ & $G_4$ & $\overline{G_4}$ & $G_5$ & $\overline{G_5}$ & $G_6$ & $\overline{G_6}$ & $G_7$ & $\overline{G_7}$\\
\hline
'` & -- & -- & -- & -- & -- & -- & $2$ & -- & $2$ & -- & $6$ & -- & $6$ & -- & $6$\\
B & $406$ & $208$ & $202$ & $208$ & $194$ & $208$ & $200$ & $200$ & $212$ & $202$ & $190$ & $196$ & $200$ & $190$ & $190$\\
C & -- & $180$ & $192$ & $180$ & $208$ & $180$ & $190$ & $196$ & $166$ & $192$ & $198$ & $204$ & $178$ & $216$ & $198$\\
D & -- & $18$ & $12$ & $18$ & $4$ & $18$ & $14$ & $10$ & $26$ & $12$ & $12$ & $6$ & $22$ & -- & $12$\\
\hline \vspace{-0.08in} \\
\hline
&& $G_8$ & $\overline{G_8}$ & $G_9$ & $\overline{G_9}$ & $G_{10}$ & $\overline{G_{10}}$ & $G_{11}$ & $\overline{G_{11}}$ & $G_{12}$ & $\overline{G_{12}}$ & $G_{13}$ & $\overline{G_{13}}$ & $G_{14}$ & $\overline{G_{14}}$\\
\hline
'` && -- & $8$ & -- & $12$ & $2$ & $4$ & $2$ & $6$ & $2$ & $6$ & $2$ & $12$ & $2$ & $14$\\
B && $204$ & $182$ & $232$ & $102$ & $206$ & $190$ & $198$ & $208$ & $198$ & $196$ & $196$ & $188$ & $228$ & $212$\\
C && $188$ & $208$ & $168$ & $192$ & $178$ & $204$ & $194$ & $162$ & $194$ & $186$ & $198$ & $184$ & $170$ & $166$\\
D && $14$ & $8$ & $6$ & -- & $20$ & $8$ & $12$ & $30$ & $12$ & $18$ & $10$ & $22$ & $6$ & $14$\\
\hline \vspace{-0.08in} \\
\cline{1-14}
&& $G_{15}$ & $\overline{G_{15}}$ & $G_{16}$ & $\overline{G_{16}}$ & $G_{17}$ & $\overline{G_{17}}$ & $G_{18}$ & $\overline{G_{18}}$ & $G_{19}$ & $\overline{G_{19}}$ & $G_{20}$ & $\overline{G_{20}}$\\
\cline{1-14}
'` && $4$ & $4$ & $4$ & $8$ & $6$ & $6$ & $6$ & $12$ & $8$ & $10$ & $8$ & $12$\\
B && $196$ & $194$ & $204$ & $186$ & $202$ & $184$ & $202$ & $184$ & $186$ & $194$ & $194$ & $180$\\
C && $192$ & $196$ & $176$ & $200$ & $174$ & $210$ & $174$ & $192$ & $200$ & $178$ & $184$ & $200$\\
D && $14$ & $12$ & $22$ & $12$ & $24$ & $6$ & $24$ & $18$ & $12$ & $24$ & $20$ & $14$\\
\cline{1-14}
\end{tabular}
\end{center}}

\quad

\subsubsection{Strongly regular graphs with the parameters $(v, k, \lambda, \mu) = (40, 12, 2, 4)$}

We have $28$ graphs as the strongly regular graphs with the parameters $(v, k, \lambda, \mu) = (40, 12, 2, 4)$ up to isomorphic (classified by E. Spence, see \cite{Se1}). The vectors of crystal lattices from each graph has $(201, 480, 10, 3)$-configuration and the distance set $A(X) = \{ -1, \: \pm 13/134, \: \pm 6/67, \: \pm 1/67, \: \pm 1/134, \: 0 \}$. However, for the numbers $|A_x(X, a)|$ for each $x \in X$ from the graphs, we have the following two types:
\begin{center}
\begin{tabular}{cccccccc}
\hline
& $a$ & $-1$ & $\pm 13/134$ & $\pm 6/67$ & $\pm 1/67$ & $\pm 1/134$ & $0$\\
\hline
A & $|A_x(X, a)|$ & $1$ & $4$ & $18$ & $27$ & $108$ & $164$\\
B & $|A_x(X, a)|$ & $1$ & $4$ & $18$ & $25$ & $116$ & $152$\\
\hline
\end{tabular}
\end{center}

On the other hand, the complements of these $28$ graphs have the same parameters $(v, k, \lambda, \mu) = (40, 27, 18, 18)$, and the vectors of crystal lattices from each graph has $(501, 1080, 10, 3)$-configuration and the distance set $A(X) = \{ -1, \: \pm 13/334, \: \pm 25/668, \: \pm 1/334, \: \pm 1/668, \: 0\}$. However, for the numbers $|A_x(X, a)|$ for each $x \in X$ from the graphs, we have the following four types:
\begin{center}
\begin{tabular}{cccccccc}
\hline
& $a$ & $-1$ & $\pm 13/334$ & $\pm 25/668$ & $\pm 1/334$ & $\pm 1/668$ & $0$\\
\hline
C & $|A_x(X, a)|$ & $1$ & $36$ & $16$ & $40$ & $240$ & $414$\\
D & $|A_x(X, a)|$ & $1$ & $36$ & $16$ & $38$ & $248$ & $402$\\
\hline
\end{tabular}
\end{center}

Now, we have eleven combinations of the numbers of the vectors of each type from each graphs and its complement, then we have the following table:
\begin{center}
\begin{tabular}{cccccccccccc}
\hline
& $(1)$ & $(2)$ & $(3)$ & $(4)$ & $(5)$ & $(6)$ & $(7)$ & $(8)$ & $(9)$ & $(10)$ & $(11)$\\
\hline
'` & -- & $48$ & $84$ & $96$ & $120$ & $144$ & $156$ & $192$ & $264$ & $288$ & $480$\\
B & $480$ & $432$ & $396$ & $384$ & $360$ & $336$ & $324$ & $288$ & $216$ & $192$ & --\\
\hline
\hline
C & $600$ & $648$ & $684$ & $696$ & $720$ & $744$ & $756$ & $792$ & $864$ & $888$ & $1080$\\
D & $480$ & $432$ & $396$ & $384$ & $360$ & $336$ & $324$ & $288$ & $216$ & $192$ & --\\
\hline
\hline
number of graphs & $1$ & $3$ & $5$ & $3$ & $1$ & $4$ & $2$ & $4$ & $2$ & $1$ & $2$\\
\hline
\end{tabular}
\end{center}

\newpage

\section{Graphs from association schemes}\label{sec-as}

We can consider graphs from every adjacency matrices of association schemes (See Definition \ref{def-as}).

For example, we have {\it distance regular graphs} (See \cite{BCN}). Let $X_0 = (V, E)$ be a regular graph of diameter $d$, and define $s_i := \{ (x, y) \in V \times V \ ; \ d(x, y) = i \}$ for every $0 \leqslant i \leqslant d$. If $(V, \{ s_i \}_{0 \leqslant i \leqslant d})$ is an association scheme, then the graph $X_0$ is a distance regular graph. In particular, if $d \leqslant 2$, then the graph $X_0$ is a strongly regular graph.

In this section, we consider graphs from various association schemes, where we do not limit to distance regular graph.

\subsection{Hamming graph $H(m, q)$}\label{sec-ham}

Let $\Omega$ be a finite set of order $q \geqslant 2$, and we denote $V := \Omega^m$. We define a distance on $V$ by $d(x, y) := | \{ i \ ; \ x_i \ne y_i \} |$ for every $x = (x_1, \ldots, x_m), y = (y_1, \ldots, y_m) \in V$, we also define the relation $s_i := \{ (x, y) \in V \times V \ ; \ d(x, y) = i \}$ for every $0 \leqslant i \leqslant m$. Then, we have an association scheme $(V, \{ s_i \}_{0 \leqslant i \leqslant m})$, which is called {\it Hamming scheme}. Also, $(V, s_1)$ is a distance regular graph, that is {\it Hamming graph $H(m, q)$}. In this paper, as a generalization of Hamming graph, we call the graph $(V, s_k)$ the {\it Hamming graph of order $k$} for every $0 \leqslant k \leqslant m$. Note that every Hamming graph of order $k$ is vertex and edge-transitive.

In particular, if $m = 2$, then the Hamming graph $H(2, q)$ is the square lattice graph $L_2(q)$ (See Section \ref{sec-l2}). On the other hand, if $q = 2$, then the Hamming graph $H(m, 2)$ is isometric to the graph from $m$-dimensional regular $2 m$-cell polytope (See Section \ref{sec-cube}).

The Hamming graphs of order $k \geqslant 2$ are not always connected. In this section, we consider the crystal lattice from each connected subgraph, every number $P$ in the following table are the number of connected subgraphs. Here, since the Hamming graphs are vertex-transitive, all of the connected subgraphs are isomorphic to each other. We also consider whether distance regular (`DR') or not, moreover, consider whether strongly regular (`SR') or not (`T'rue or `F'aulse).

If $q=2$, for the $(d, n, s, t)$-configurations of the vectors of the constructed crystal lattices, we have the following table:\vspace{-0.05in}

\begin{table}[h]
\begin{center}
\begin{tabular}{cccccccc}
\multicolumn{8}{c}{Hamming graph $H(m, q)$}\\
\hline
$(m, q)$ & $k$ & $v$ & $e$ & P & {\small DR} & {\small SR}\\
\hline
$(3, 2)$ & $1$ & $8$ & $12$ & $1$ & T & F\\
& $2$ && $12$ & $2$ & T & T & $2 \cdot K_4$\\
& $3$ && $4$ & $4$ & T & T & $4 \cdot K_2$\\
$(4, 2)$ & $1$ & $16$ & $32$ & $1$ & T & F\\
& $2$ && $48$ & $2$ & T & T & $2 \cdot CP(4)$\\
& $3$ && $32$ & $1$ & T & F\\
& $4$ && $8$ & $8$ & T & T & $8 \cdot K_2$\\
$(5, 2)$ & $1$ & $32$ & $80$ & $1$ & T & F\\
& $2$ && $160$ & $2$ & T & T & $2 \cdot$ $\overline{\text{Clebsch}}$\\
& $3$ && $160$ & $1$ & T & F\\
& $4$ && $80$ & $2$ & T & T & $2 \cdot$ Clebsch\\
& $5$ && $16$ & $16$ & T & T & $16 \cdot K_2$\\
$(6, 2)$ & $1$ & $64$ & $192$ & $1$ & T & F\\
& $2$ && $480$ & $2$ & T & F\\
& $3$ && $640$ & $1$ & F\\
& $4$ && $480$ & $2$ & T & F\\
& $5$ && $192$ & $1$ & T & F\\
& $6$ && $32$ & $32$ & T & T & $32 \cdot K_2$\\
$(7, 2)$ & $1$ & $128$ & $448$ & $1$ & T & F\\
& $2$ && $1344$ & $2$ & T & F\\
& $3$ && $2240$ & $1$ & F\\
& $4$ && $2240$ & $2$ & T & T & $(64, 35, 18, 20)$\\
& $5$ && $1344$ & $1$ & F\\
& $6$ && $448$ & $2$ & T & F\\
& $7$ && $64$ & $64$ & T & T & $64 \cdot K_2$\\
\hline
\end{tabular}
\qquad $\Rightarrow$ \qquad
\begin{tabular}{cccc}
\multicolumn{4}{c}{Vectors of lattice}\\
\hline
$d$ & $n$ & $s$ & $t$\\
\hline
$5$ & $24$ & $9$ & $3$\\
$3$ & $12$ & $4$ & $3$\\
$1$ & $2$ & $1$ & --\\
$17$ & $64$ & $11$ & $3$\\
$17$ & $48$ & $10$ & $3$\\
$17$ & $64$ & $11$ & $3$\\
$1$ & $2$ & $1$ & --\\
$49$ & $160$ & $17$ & $3$\\
$65$ & $160$ & $10$ & $3$\\
$129$ & $320$ & $9$ & $3$\\
$25$ & $80$ & $8$ & $3$\\
$1$ & $2$ & $1$ & --\\
$129$ & $384$ & $19$ & $3$\\
$209$ & $480$ & $16$ & $3$\\
$577$ & $1280$ & $13$ & $3$\\
$209$ & $480$ & $12$ & $3$\\
$129$ & $320$ & $9$ & $3$\\
$1$ & $2$ & $1$ & --\\
$321$ & $896$ & $25$ & $3$\\
$609$ & $1344$ & $14$ & $3$\\
$2113$\\
$1057$\\
$1207$\\
$161$ & $448$ & $12$ & $3$\\
$1$ & $2$ & $1$ & --\\
\hline
\end{tabular}
\end{center}\quad
\caption{Hamming graph $H(m, 2)$ of order $k \leqslant m$}\vspace{-0.27in}
\end{table}
Here, Hamming graph $H(4, 2)$ of order $1$ and that of order $3$ is isomorphic to each other.

\newpage

If $q \geqslant 3$, for every case which we consider in the following table, the Hamming graph is connected. For the $(d, n, s, t)$-configurations of the vectors of the constructed crystal lattices, we have the following table:

\begin{table}[h]
\begin{center}
\begin{tabular}{ccccccc}
\multicolumn{7}{c}{Hamming graph $H(m, q)$}\\
\hline
$(m, q)$ & $k$ & $v$ & $e$ & {\small DR} & {\small SR}\\
\hline
$(3, 3)$ & $1$ & $27$ & $81$ & T & F\\
& $2$ && $162$ & F\\
& $3$ && $108$ & F\\
$(4, 3)$ & $1$ & $81$ & $324$ & T & F\\
& $2$ && $972$ & T & T & $(81, 24, 9, 6)$\\
& $3$ && $1296$ & T & T & $(81, 32, 13, 12)$\\
& $4$ && $648$ & F\\
$(5, 3)$ & $1$ & $243$ & $1215$ & T & F\\
& $2$ && $4860$ & F\\
& $3$ && $9720$ & F\\
& $4$ && $9720$ & F\\
& $5$ && $3888$ & F\\
\hline
\end{tabular}
\qquad $\Rightarrow$ \qquad
\begin{tabular}{cccc}
\multicolumn{4}{c}{Vectors of lattice}\\
\hline
$d$ & $n$ & $s$ & $t$\\
\hline
$55$ & $162$ & $16$ & $3$\\
$136$ & $324$ & $22$ & $3$\\
$82$ & $216$ & $24$ & $3$\\
$244$ & $648$ & $16$ & $3$\\
$892$\\
$1216$\\
$568$ & $1296$ & $40$ & $3$\\
$973$\\
$4618$\\
$9478$\\
$9478$\\
$3646$\\
\hline
\end{tabular}\vspace{0.1in}

\begin{tabular}{cccccccc}
\hline
$(m, q)$ & $k$ & $v$ & $e$ & {\small DR} & {\small SR}\\
\hline
$(3, 4)$ & $1$ & $64$ & $288$ & T & F\\
& $2$ && $864$ & T & T & $(64, 27, 10, 12)$\\
& $3$ && $864$ & F\\
$(4, 4)$ & $1$ & $256$ & $1536$ & T & F\\
& $2$ && $6912$ & F\\
& $3$ && $13824$ & F\\
& $4$ && $10368$ & F\\
\hline
\end{tabular}
\qquad $\Rightarrow$ \qquad
\begin{tabular}{cccc}
\hline
$d$ & $n$ & $s$ & $t$\\
\hline
$225$ & $576$ & $14$ & $3$\\
$801$ & $1728$ & $10$ & $3$\\
$801$ & $1728$ & $26$ & $3$\\
$1281$\\
$6657$\\
$13569$\\
$10113$\\
\hline
\end{tabular}\vspace{0.1in}

\begin{tabular}{cccccccc}
\hline
$(m, q)$ & $k$ & $v$ & $e$ & {\small DR} & {\small SR} & \hspace{0.8in} \quad \\
\hline
$(3, 5)$ & $1$ & $125$ & $750$ &  T & F\\
& $2$ && $3000$ &  F\\
& $3$ && $4000$ &  F\\
\hline
\end{tabular}
\qquad $\Rightarrow$ \qquad
\begin{tabular}{cccc}
\hline
$d$ & $n$ & $s$ & $t$\\
\hline
$626$ & $1500$ & $16$ & $3$\\
$2876$\\
$3876$\\
\hline
\end{tabular}
\end{center}\quad
\caption{Hamming graph $H(m, q)$ of order $k \leqslant m$}
\end{table}

\newpage

\subsection{Johnson graph $J(q, m)$}\label{sec-john}

Let $\Omega$ be a finite set of order $q \geqslant 2$, and we denote $V$ be a family of all the $m$-point subsets of $\Omega$ for a positive integer $m \leqslant q/2$. We define a distance on $V$ by $d(x, y) := k - |x \cap y|$ for every $x, y \in V$, we also define the relation $s_i := \{ (x, y) \in V \times V \ ; \ d(x, y) = i \}$ for every $0 \leqslant i \leqslant m$. Then, we have an association scheme $(V, \{ s_i \}_{0 \leqslant i \leqslant m})$, which is called {\it Johnson scheme}. Also, $(V, s_1)$ is a distance regular graph, that is {\it Johnson graph $J(q, m)$}. In this paper, as a generalization of Johnson graph, we call the graph $(V, s_k)$ the {\it Johnson graph of order $k$} for every $0 \leqslant k \leqslant m$. Note that every Johnson graph of order $k$ is vertex and edge-transitive.

In particular, if $m = 2$, then the Johnson graph $J(q, 2)$ is the triangular graph $T(q)$ (See Section \ref{sec-tri}). Moreover, Johnson graph $J(q, m)$ of order $m$ is called the {\it Kneser graph}.

Some Johnson graphs of order $k \geqslant 2$ are not connected, then we consider the crystal lattice from each connected subgraph. Here, since the Johnson graphs are vertex-transitive, all of the connected subgraphs are isomorphic to each other. We also consider whether distance regular (`DR') (moreover, strongly regular (`SR')) or not.

If $m \geqslant 3$, for the $(d, n, s, t)$-configurations of the vectors of the constructed crystal lattices, we have the following table:\vspace{-0.05in}

\begin{table}[h]
\begin{center}
\begin{tabular}{ccccccc}
\multicolumn{7}{c}{Johnson graph $J(q, m)$}\\
\hline
$(q, m)$ & $k$ & $v$ & $e$ & {\small DR} & {\small SR}\\
\hline
$(6, 3)$ & $1$ & $20$ & $90$ & T & F\\
& $2$ && $90$ & T & F\\
& $3$ && $10$ & T & T & $10 \cdot K_2$ (non-connected)\\
$(7, 3)$ & $1$ & $35$ & $210$ & T & F\\
& $2$ && $315$ & T & T & $(35, 18, 9, 9)$\\
& $3$ && $70$ & T & F\\
$(8, 3)$ & $1$ & $56$ & $420$ & T & F\\
& $2$ && $840$ & F\\
& $3$ && $280$ & F\\
$(9, 3)$ & $1$ & $84$ & $756$ & T & F\\
& $2$ && $1890$ & F\\
& $3$ && $840$ & F\\
\hline
\end{tabular}
\qquad $\Rightarrow$ \qquad
\begin{tabular}{cccc}
\multicolumn{4}{c}{Vectors of lattice}\\
\hline
$d$ & $n$ & $s$ & $t$\\
\hline
$71$ & $180$ & $16$ & $3$\\
$71$ & $180$ & $16$ & $3$\\
$1$ & $2$ & $1$ & --\\
$176$ & $420$ & $16$ & $3$\\
$281$ & $630$ & $10$ & $3$\\
$36$ & $140$ & $9$ & $3$\\
$365$ & $840$ & $14$ & $3$\\
$785$ & $1680$ & $26$ & $3$\\
$225$ & $560$ & $18$ & $3$\\
$673$ & $1512$ & $16$ & $3$\\
$1807$\\
$757$ & $1680$ & $22$ & $3$\\
\hline
\end{tabular}\vspace{0.1in}

\begin{tabular}{ccccccc}
\hline
$(q, m)$ & $k$ & $v$ & $e$ & {\small DR} & {\small SR}\\
\hline
$(8, 4)$ & $1$ & $70$ & $560$ & T & F\\
& $2$ && $1260$ & F\\
& $3$ && $560$ & F\\
& $3$ && $35$ & T & T & $35 \cdot K_2$ (non-connected)\\
$(9, 4)$ & $1$ & $126$ & $1260$ & T & F\\
& $2$ && $3780$ & F\\
& $3$ && $2520$ & F\\
& $4$ && $315$ & T & F\\
\hline
\end{tabular}
\qquad $\Rightarrow$ \qquad
\begin{tabular}{cccc}
\hline
$d$ & $n$ & $s$ & $t$\\
\hline
$491$ & $1120$ & $22$ & $3$\\
$1191$\\
$491$ & $1120$ & $30$ & $3$\\
$1$ & $2$ & $1$ & --\\
$1135$\\
$3655$\\
$2395$\\
$190$ & $630$ & $12$ & $3$\\
\hline
\end{tabular}
\end{center}\quad
\caption{Johnson graph $J(q, m)$ of order $k \leqslant m$}\vspace{-0.27in}
\end{table}
Here, Johnson graph $J(6, 3)$ of order $1$ and that of order $2$ is isomorphic to each other.

\newpage

\subsection{Cyclotomic graph $Cyc(q, m)$}\label{sec-cyc}

Let $q$ be a prime power, $m$ be a divisor of $q-1$ less than $q-1$, and $\omega$ be a primitive element of the finite field $\mathbb{F}_q$. We define $C_0 = \{ 0 \}, C_1, \ldots , C_d$ a partition of $\mathbb{F}_q$ as follows:
\begin{equation*}
C_i := \{ \omega^{i + a d} \: ; \: a = 0, 1, \ldots , (q-1)/m \} \quad \text{for every} \ 1 \leqslant i \leqslant m.
\end{equation*}
Now, we define $s_i := \{ (x, y) \: ; \: x - y \in C_i \}$ the relation on $\mathbb{F}_q$ for every $0 \leqslant i \leqslant m$. Then, we have an association scheme $(\mathbb{F}_q, \{ s_i \}_{0 \leqslant i \leqslant m})$, which is called {\it cyclotomic scheme}. Also, if the scheme $(\mathbb{F}_q, \{ s_i \}_{0 \leqslant i \leqslant m})$ is symmetric, we call $(\mathbb{F}_q, s_1)$ the {\it cyclotomic graph $Cyc(q, m)$}. Note that every graph $(\mathbb{F}_q, s_1)$ is vertex and edge-transitive but not always distance regular, and also note that all the graphs $(\mathbb{F}_q, s_i)$ for $1 \leqslant i \leqslant m$ are isomorphic to each other.

In particular, if $m = 2$, then the scheme is symmetric if only if $q \equiv 1 \pmod{4}$, and then the cyclotomic graph $cyc(q, 2)$ is the Paley graph $P(q)$ (See Section \ref{sec-pq}).

All the cyclotomic graphs $Cyc(q, m)$ for $q \leqslant 49$ are in the following table:\vspace{-0.05in}

\begin{table}[h]
{\small \begin{center}
\begin{tabular}{cccl}
\hline
$q$ & $m$ & $e$\\
\hline
$2$ & $1$ & $1$ & $K_2$\\
$3$ & $1$ & $3$ & $K_3$\\
$4$ & $1$ & $6$ & $K_4$\\
$5$ & $1$ & $10$ & $K_5$\\
 & $2$ & $5$ & $P(5)$\\
$7$ & $1$ & $21$ & $K_7$\\
 & $3$ & $7$ & $C_7$\\
$8$ & $1$ & $28$ & $K_8$\\
$9$ & $1$ & $36$ & $K_9$\\
 & $2$ & $18$ & $P(9)$\\
 & $4$ & $9$ & $3 \cdot K_3$\\
$11$ & $1$ & $55$ & $K_{11}$\\
 & $5$ & $11$ & $C_{11}$\\
$13$ & $1$ & $78$ & $K_{13}$\\
 & $2$ & $39$ & $P(13)$\\
 & $3$ & $26$\\
 & $6$ & $13$ & $C_{13}$\\
$16$ & $1$ & $120$ & $K_{16}$\\
\hline
\end{tabular}
\
\begin{tabular}{cccl}
\hline
$q$ & $m$ & $e$\\
\hline
$16$ & $3$ & $40$ & Clebsch\\
& $5$ & $24$ & $4 \cdot K_4$\\
$17$ & $1$ & $136$ & $K_{17}$\\
& $2$ & $68$ & $P(17)$\\
& $4$ & $34$\\
& $8$ & $17$ & $C_{17}$\\
$19$ & $1$ & $171$ & $K_{19}$\\
& $3$ & $57$\\
& $9$ & $19$ & $C_{19}$\\
$23$ & $1$ & $253$ & $K_{23}$\\
& $11$ & $23$ & $C_{23}$\\
$25$ & $1$ & $300$ & $K_{25}$\\
& $2$ & $150$ & $P(25)$\\
& $3$ & $100$ & $L_2(5)$\\
& $4$ & $75$\\
& $6$ & $50$ & $5 \cdot K_5$\\
& $12$ & $25$ & $5 \cdot C_5$\\
$27$ & $1$ & $351$ & $K_{27}$\\
\hline
\end{tabular}
\
\begin{tabular}{cccl}
\hline
$q$ & $m$ & $e$\\
\hline
$27$ & $13$ & $27$ & $9 \cdot K_3$\\
$29$ & $1$ & $406$ & $K_{29}$\\
 & $2$ & $203$ & $P(29)$\\
 & $7$ & $58$\\
 & $14$ & $29$ & $C_{29}$\\
$31$ & $1$ & $465$ & $K_{31}$\\
 & $3$ & $155$\\
 & $5$ & $93$\\
 & $15$ & $31$ & $C_{31}$\\
$32$ & $1$ & $496$ & $K_{32}$\\
$37$ & $1$ & $666$ & $K_{37}$\\
 & $2$ & $333$ & $P(37)$\\
 & $3$ & $222$\\
 & $6$ & $111$\\
 & $9$ & $74$\\
 & $18$ & $37$ & $C_{37}$\\
$41$ & $1$ & $840$ & $K_{41}$\\
 & $2$ & $420$ & $P(41)$\\
\hline
\end{tabular}
\
\begin{tabular}{cccl}
\hline
$q$ & $m$ & $e$\\
\hline
$41$ & $4$ & $210$\\
 & $5$ & $168$\\
 & $10$ & $84$\\
 & $20$ & $41$ & $C_{41}$\\
$43$ & $1$ & $903$ & $K_{43}$\\
 & $3$ & $301$\\
 & $7$ & $129$\\
 & $21$ & $43$ & $C_{43}$\\
$47$ & $1$ & $1081$ & $K_{47}$\\
 & $23$ & $47$ & $C_{47}$\\
$49$ & $1$ & $1176$ & $K_{49}$\\
& $2$ & $588$ & $P(49)$\\
& $3$ & $392$\\
& $4$ & $294$ & $L_2(7)$\\
& $6$ & $196$\\
& $8$ & $147$ & $7 \cdot K_7$\\
& $12$ & $98$\\
& $24$ & $49$ & $7 \cdot C_7$\\
\hline
\end{tabular}
\end{center}}
\caption{Cyclotomic graph $Cyc(q, m)$ for $q \leqslant 49$}\vspace{-0.27in}
\end{table}

Here, in the above table, graphs which are distance regular and not strongly regular are only {\it cycle graphs $C_p$}. Furthermore, for the $(d, n, s, t)$-configurations of the vectors of the constructed crystal lattices from graphs which are not distance regular, we have the following table:\vspace{-0.05in}

\begin{table}[h]
\begin{center}
\begin{tabular}{cccccccc}
\hline
$q$ & $m$ & $e$ & \ & $d$ & $n$ & $s$ & $t$\\
\hline
$13$ & $3$ & $26$ && $14$ & $52$ & $16$ & $3$\\
$17$ & $4$ & $34$ && $18$ & $68$ & $24$ & $3$\\
$19$ & $3$ & $57$ && $37$ & $114$ & $35$ & $3$\\
$25$ & $4$ & $75$ && $51$ & $150$ & $29$ & $3$\\
$29$ & $7$ & $58$ && $30$ & $116$ & $44$ & $3$\\
$31$ & $3$ & $155$ && $125$ & $310$ & $24$ & $3$\\
 & $5$ & $93$ && $63$ & $186$ & $51$ & $3$\\
$37$ & $3$ & $222$ && $186$ & $444$ & $24$ & $3$\\
 & $6$ & $111$ && $75$ & $222$ & $69$ & $3$\\
\hline
\end{tabular}
\qquad
\begin{tabular}{cccccccc}
\hline
$q$ & $m$ & $e$ & \ & $d$ & $n$ & $s$ & $t$\\
\hline
$37$ & $9$ & $74$ && $38$ & $148$ & $56$ & $3$\\
$41$ & $4$ & $210$ && $165$ & $410$ & $46$ & $3$\\
 & $5$ & $168$ && $124$ & $328$ & $62$ & $3$\\
 & $10$ & $84$ && $42$ & $164$ & $62$ & $3$\\
$43$ & $3$ & $301$ && $259$ & $602$ & $24$ & $3$\\
 & $7$ & $129$ && $87$ & $258$ & $79$ & $3$\\
$49$ & $3$ & $392$ && $344$ & $784$ & $26$ & $3$\\
& $6$ & $196$ && $148$ & $392$ & $70$ & $3$\\
& $12$ & $98$ && $50$ & $196$ & $34$ & $3$\\
\hline
\end{tabular}
\end{center}\quad
\caption{Cyclotomic graph which are not distance regular}\vspace{-0.27in}
\end{table}

Furthermore, in the above table, the graphs which are not distance regular and have prime number vertices are {\it circulant graph ${Ci}_n(i_1, \ldots , i_k)$} as follows:

\begin{center}
\begin{tabular}{ccl}
\hline
$q$ & $m$\\
\hline
$13$ & $3$ & ${Ci}_{13} (1, 5)$\\
$17$ & $4$ & ${Ci}_{17} (1, 4)$\\
$19$ & $3$ & ${Ci}_{19} (1, 7, 8)$\\
$29$ & $7$ & ${Ci}_{29} (1, 10)$\\
\hline\\
\end{tabular}
\qquad
\begin{tabular}{ccl}
\hline
$q$ & $m$\\
\hline
$31$ & $3$ & ${Ci}_{31} (1, 2, 4, 8, 15)$\\
& $5$ & ${Ci}_{31} (1, 5, 6)$\\
$37$ & $3$ & ${Ci}_{37} (1, 6, 8, 10, 11, 14)$\\
& $6$ & ${Ci}_{37} (1, 10, 11)$\\
& $9$ & ${Ci}_{37} (1, 6)$\\
\hline
\end{tabular}
\qquad
\begin{tabular}{ccl}
\hline
$q$ & $m$\\
\hline
$41$ & $4$ & ${Ci}_{41} (1, 4, 10, 16, 18)$\\
& $5$ & ${Ci}_{41} (1, 3, 9, 14)$\\
& $10$ & ${Ci}_{41} (1, 9)$\\
$43$ & $3$ & ${Ci}_{43} (1, 2, 4, 8, 11, 16, 21)$\\
& $7$ & ${Ci}_{43} (1, 6, 7)$\\
\hline
\end{tabular}
\end{center}

\newpage

On the other hand, the remaining graphs which are not distance regular are the following:
\begin{figure}[h]
\begin{center}
$Cyc(25, 4)$\includegraphics[width=1.5in]{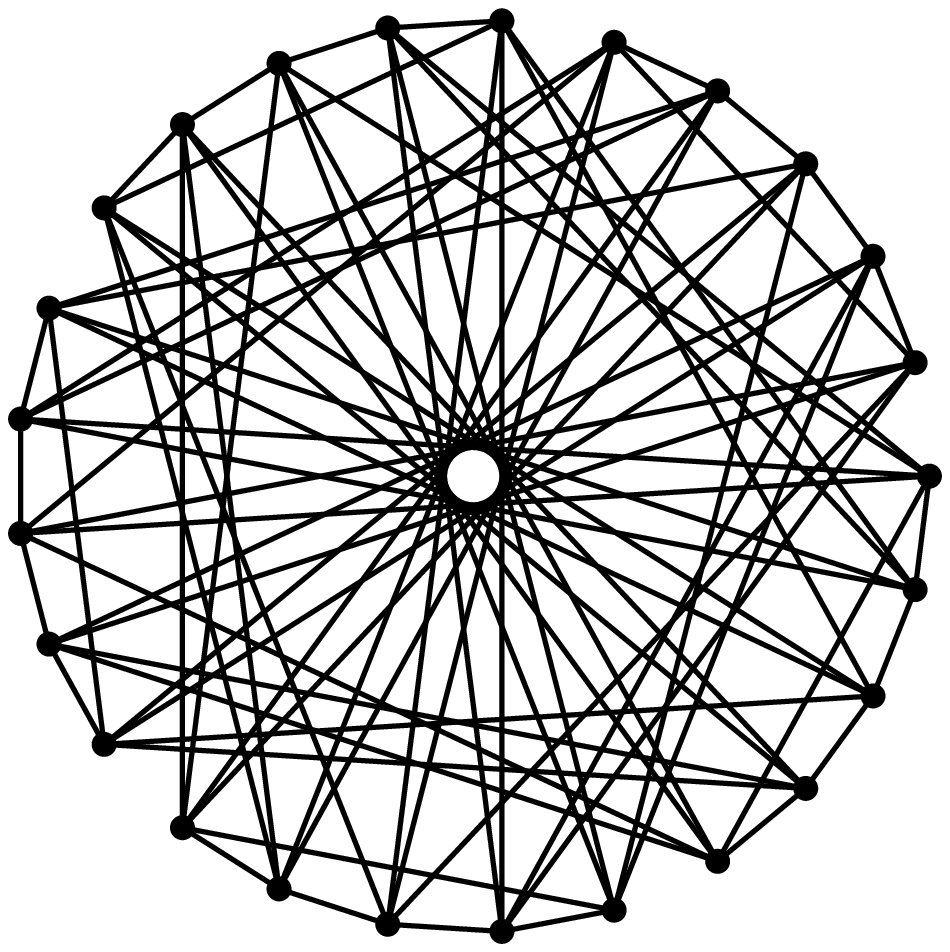} \quad
$Cyc(49, 3)$\includegraphics[width=1.5in]{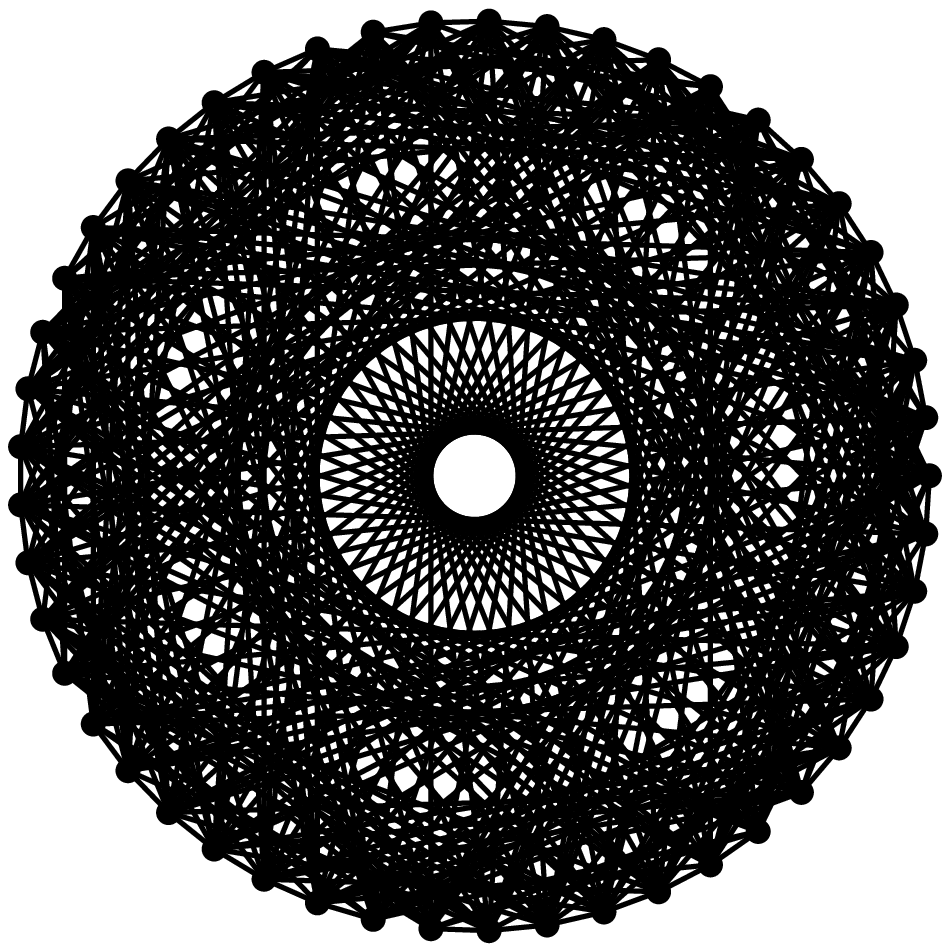}
\quad\\
$Cyc(49, 6)$\includegraphics[width=1.5in]{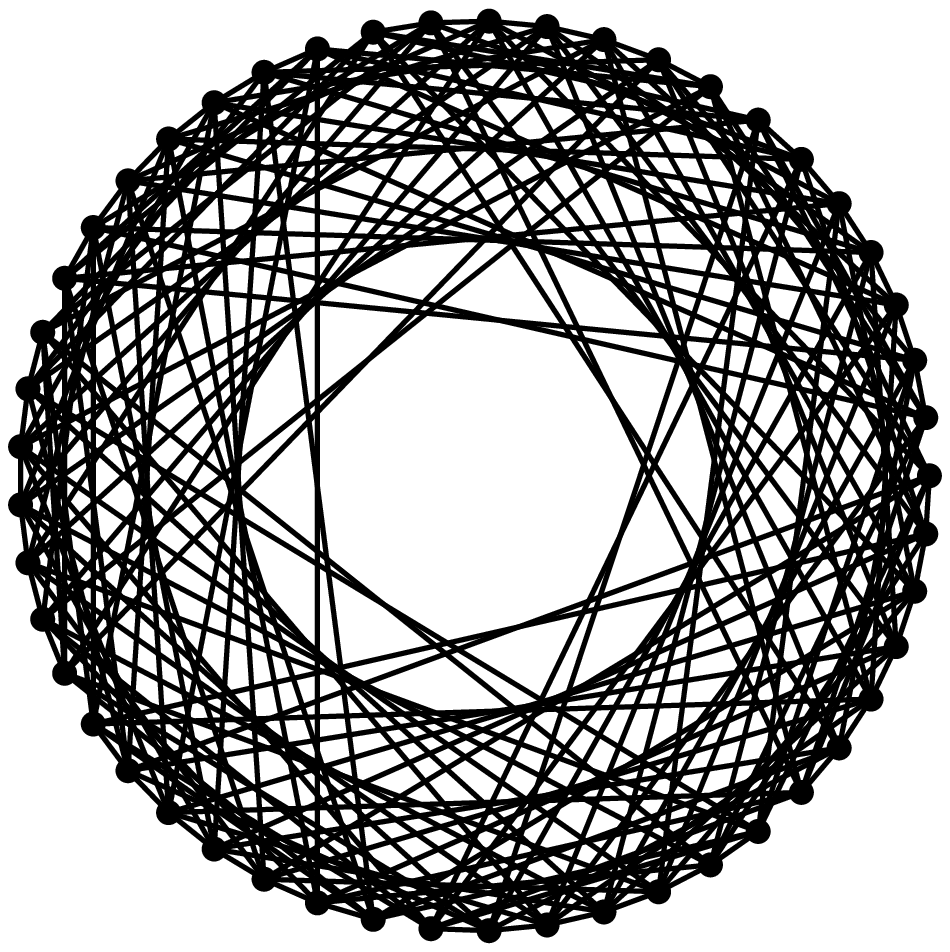} \qquad
$Cyc(49, 12)$\includegraphics[width=1.5in]{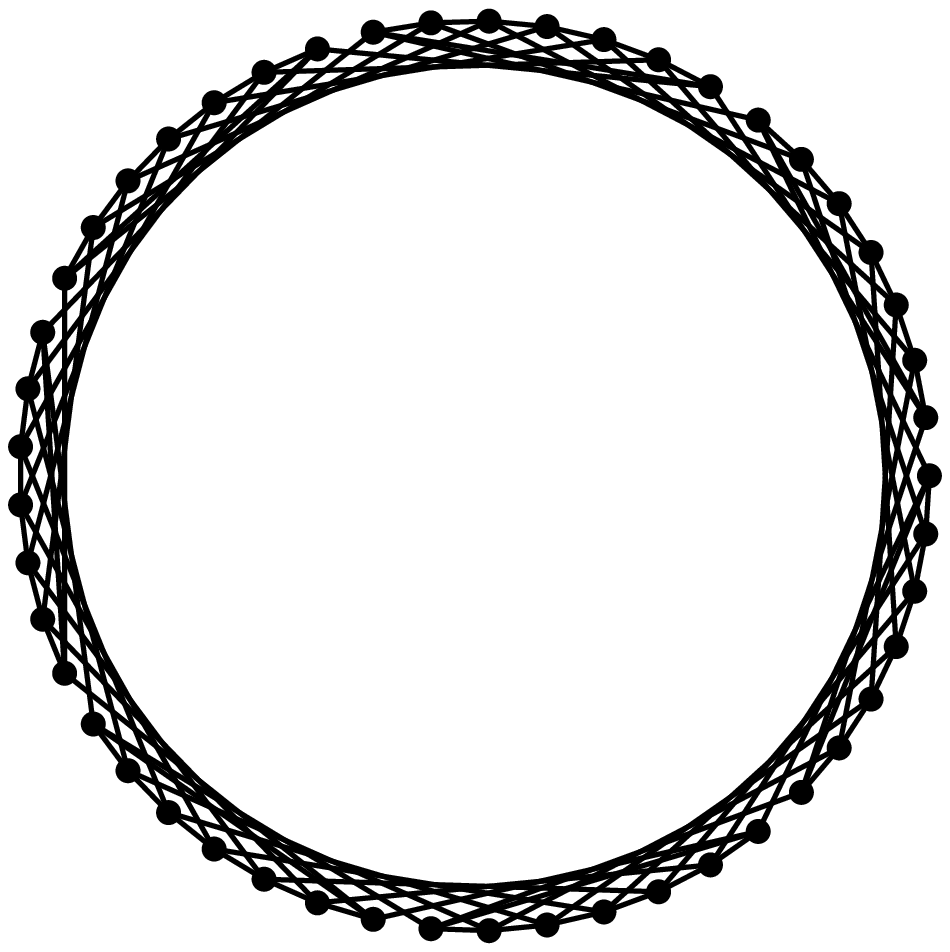}
\end{center}
\caption{certain cyclotomic graph}\label{fig-net}
\end{figure}

\subsection{Graphs from the other association schemes}\label{sec-oth-as}

By \cite{HM1,HM2,HM3,HM4,Hm,HS,N,Sy,SS}, association schemes with small vertices are classified. The following table is the number of association schemes up to isomorphic (`AS') and non-Schurian association schemes (`nSS') with $|X|$ vertices (See \cite{HM1}).\vspace{-0.05in}

\begin{table}[h]
\begin{center}
\begin{tabular}{rrr}
\hline
$|X|$ & AS & nSS\\
\hline
$1$ & $1$ & $0$\\
$2$ & $1$ & $0$\\
$3$ & $2$ & $0$\\
$4$ & $4$ & $0$\\
$5$ & $5$ & $0$\\
$6$ & $8$ & $0$\\
$7$ & $4$ & $0$\\
$8$ & $21$ & $0$\\
$9$ & $12$ & $0$\\
$10$ & $13$ & $0$\\
\hline
\end{tabular}
\quad
\begin{tabular}{rrr}
\hline
$|X|$ & AS & nSS\\
\hline
$11$ & $4$ & $0$\\
$12$ & $59$ & $0$\\
$13$ & $6$ & $0$\\
$14$ & $16$ & $0$\\
$15$ & $25$ & $1$\\
$16$ & $222$ & $16$\\
$17$ & $5$ & $0$\\
$18$ & $95$ & $2$\\
$19$ & $7$ & $1$\\
$20$ & $95$ & $0$\\
\hline
\end{tabular}
\quad
\begin{tabular}{rrr}
\hline
$|X|$ & AS & nSS\\
\hline
$21$ & $32$ & $0$\\
$22$ & $16$ & $0$\\
$23$ & $22$ & $18$\\
$24$ & $750$ & $81$\\
$25$ & $45$ & $13$\\
$26$ & $34$ & $10$\\
$27$ & $502$ & $380$\\
$28$ & $185$ & $61$\\
$29$ & $26$ & $20$\\
$30$ & $243$ & $15$\\
\hline
\end{tabular}
\quad
\begin{tabular}{rrr}
\hline
$|X|$ & AS & nSS\\
\hline
$31$ & $?$\\
$32$ & $18210$ & $13949$\\
$33$ & $27$ & $0$\\
$34$ & $20$ & $0$\\
$35$ & $?$\\
$36$ & $?$\\
$37$ & $?$\\
$38$ & $33$ & $11$\\
$39$ & $?$\\
$40$ & $?$\\
\hline
\end{tabular}
\end{center}\quad
\caption{Classification of association schemes}\vspace{-0.27in}
\end{table}

Now, we want to consider Conjecture \ref{conj-as}. By Theorem \ref{th-et} and Proposition \ref{prop-ss-et}, we can prove the conjecture for the graphs from Schurian schemes. Thus, the remaining part is for the graphs from non-Schurian association schemes. We calculated all of the graphs from non-Schurian schemes with $v$ vertices for $1 \leqslant v \leqslant 30$ and $v = 38$, where the vectors of every constructed crystal lattice have the same norm. (I used the data in \cite{HM1}. For the case of $v=32$, the number of association schemes is too large for my calculation.) Thus, we can prove Conjecture \ref{conj-as} for the graphs with  $v$ vertices for $1 \leqslant v \leqslant 30$, and $v = 33, 34, 38$, by numerical calculation.

\newpage

\section{List of edge-transitive graphs}\label{sec-et}
In this section, we consider connected simple graphs with at most $15$ vertices. The following table is the list of them. For vertex-transitive graphs, we also consider whether distance regular (`DR') (moreover, strongly regular (`SR')) or not (`T'rue or `F'aulse).

\begin{table}[h]
\begin{flushleft}``Vertex-transitive graph''\end{flushleft}
{\small
\begin{center}
\begin{tabular}{lcccc}
\hline
& $v$ & $e$ & {\small DR} & {\small SR}\\
\hline
$K_2$ & $2$ & $1$ & T & T\\
$K_3$ & $3$ & $3$ & T & T\\
$C_4$ & $4$ & $4$ & T & T\\
$K_4$ & $4$ & $6$ & T & T\\
$C_5$ & $5$ & $5$ & T & T\\
$K_5$ & $5$ & $10$ & T & T\\
$C_6$ & $6$ & $6$ & T & {\bf F}\\
$K_{3, 3}$ & $6$ & $9$ & T & T\\
$CP(3)$ & $6$ & $12$ & T & T\\
$K_6$ & $6$ & $15$ & T & T\\
$C_7$ & $7$ & $7$ & T & {\bf F}\\
$K_7$ & $7$ & $21$ & T & T\\
$C_8$ & $8$ & $8$ & T & {\bf F}\\
$H(3, 2)$ & $8$ & $12$ & T & {\bf F}\\
$K_{4, 4}$ & $8$ & $16$ & T & T\\
$CP(4)$ & $8$ & $24$ & T & T\\
$K_8$ & $8$ & $28$ & T & T\\
$C_9$ & $9$ & $9$ & T & {\bf F}\\
$L_2(3)$ & $9$ & $18$ & T & T\\
$K_{3, 3, 3}$ & $9$ & $27$ & T & T\\
$K_9$ & $9$ & $36$ & T & T\\
\hline
\end{tabular}
\quad
\begin{tabular}{lcccc}
\hline
& $v$ & $e$ & {\small DR} & {\small SR}\\
\hline
$C_{10}$ & $10$ & $10$ & T & {\bf F}\\
Petersen & $10$ & $15$ & T & T\\
${Ci}_{10}(1, 3)$ & $10$ & $20$ & T & {\bf F}\\
${Ci}_{10}(1, 4)$ & $10$ & $20$ & {\bf F}\\
$K_{5, 5}$ & $10$ & $25$ & T & T\\
$T(5)$ & $10$ & $30$ & T & T\\
$CP(5)$ & $10$ & $40$ & T & T\\
$K_{10}$ & $10$ & $45$ & T & T\\
$C_{11}$ & $11$ & $11$ & T & {\bf F}\\
$K_{11}$ & $11$ & $55$ & T & T\\
$C_{12}$ & $12$ & $12$ & T & {\bf F}\\
${Ci}_{12}(1, 5)$ & $12$ & $24$ & T & {\bf F}\\
$rg(12, 4)$ & $12$ & $24$ & {\bf F}\\
$rg(12, 5)$ & $12$ & $30$ & T & {\bf F}\\
{\footnotesize ${Ci}_{12}(1, 2, 5)$} & $12$ & $36$ & {\bf F}\\
$K_{6, 6}$ & $12$ & $36$ & T & T\\
$K_{4, 4, 4}$ & $12$ & $48$ & T & T\\
$K_{3, 3, 3, 3}$ & $12$ & $54$ & T & T\\
$CP(6)$ & $12$ & $60$ & T & T\\
$K_{12}$ & $12$ & $66$ & T & T\\
$C_{13}$ & $13$ & $13$ & T & {\bf F}\\
\hline
\end{tabular}
\quad
\begin{tabular}{lcccc}
\hline
& $v$ & $e$ & {\small DR} & {\small SR}\\
\hline
${Ci}_{13}(1, 5)$ & $13$ & $26$ & {\bf F}\\
$P(13)$ & $13$ & $39$ & T & T\\
$K_{13}$ & $13$ & $78$ & T & T\\
$C_{14}$ & $14$ & $14$ & T & {\bf F}\\
$rg(14, 3)$ & $14$ & $21$ & T & {\bf F}\\
${Ci}_{14}(1, 6)$ & $14$ & $28$ & {\bf F}\\
$rg(14, 4)$ & $14$ & $28$ & T & {\bf F}\\
{\scriptsize ${Ci}_{14}(1, 3, 5)$} & $14$ & $42$ & T & {\bf F}\\
$K_{7, 7}$ & $14$ & $49$ & T & T\\
$CP(7)$ & $14$ & $84$ & T & T\\
$K_{14}$ & $14$ & $91$ & T & T\\
$C_{15}$ & $15$ & $15$ & T & {\bf F}\\
${Ci}_{15}(1, 4)$ & $15$ & $30$ & T & {\bf F}\\
$rg(15, 4)$ & $15$ & $30$ & T & {\bf F}\\
{\footnotesize ${Ci}_{15}(1, 4, 6)$} & $15$ & $45$ & {\bf F}\\
$\overline{T(6)}$ & $15$ & $45$ & T & T\\
{\tiny ${Ci}_{15}(1, 2, 4, 7)$} & $15$ & $60$ & {\bf F}\\
$T(6)$ & $15$ & $60$ & T & T\\
$K_{5, 5, 5}$ & $15$ & $75$ & T & T\\
$K_{3, \ldots , 3}$ & $15$ & $90$ & T & T\\
$K_{15}$ & $15$ & $105$ & T & T\\
\hline
\end{tabular}
\end{center}}\quad\\ \quad\\
\begin{flushleft}``Bipartite graph (other than the above graphs)''\end{flushleft}
{\small
\begin{center}
\begin{tabular}{lcccc}
\hline
& $m_1$ & $m_2$ & $v$ & $e$\\
\hline
$K_{1, 2}$ & $1$ & $2$ & $3$ & $2$\\
$K_{1, 3}$ & $1$ & $3$ & $4$ & $3$\\
$K_{1, 4}$ & $1$ & $4$ & $5$ & $4$\\
$K_{2, 3}$ & $2$ & $3$ & $5$ & $6$\\
$K_{1, 5}$ & $1$ & $5$ & $6$ & $5$\\
$K_{2, 4}$ & $2$ & $4$ & $6$ & $8$\\
$K_{1, 6}$ & $1$ & $6$ & $7$ & $6$\\
$K_{2, 5}$ & $2$ & $5$ & $7$ & $10$\\
$K_{3, 4}$ & $3$ & $4$ & $7$ & $12$\\
$K_{1, 7}$ & $1$ & $7$ & $8$ & $7$\\
$K_{2, 6}$ & $2$ & $6$ & $8$ & $12$\\
$K_{3, 5}$ & $3$ & $5$ & $8$ & $15$\\
$K_{1, 8}$ & $1$ & $8$ & $9$ & $8$\\
$K_{2, 7}$ & $2$ & $7$ & $9$ & $14$\\
$K_{3, 6}$ & $3$ & $6$ & $9$ & $18$\\
$bp(3, 6)$ & $3$ & $6$ & $9$ & $12$\\
$K_{4, 5}$ & $4$ & $5$ & $9$ & $20$\\
$K_{1, 9}$ & $1$ & $9$ & $10$ & $9$\\
$K_{2, 8}$ & $2$ & $8$ & $10$ & $16$\\
$K_{3, 7}$ & $3$ & $7$ & $10$ & $21$\\
$K_{4, 6}$ & $4$ & $6$ & $10$ & $24$\\
$bp(4, 6)$ & $4$ & $6$ & $10$ & $12$\\
\hline
\end{tabular}
\qquad
\begin{tabular}{lcccc}
\hline
& $m_1$ & $m_2$ & $v$ & $e$\\
\hline
$K_{1, 10}$ & $1$ & $10$ & $11$ & $10$\\
$K_{2, 9}$ & $2$ & $9$ & $11$ & $18$\\
$K_{3, 8}$ & $3$ & $8$ & $11$ & $24$\\
$K_{4, 7}$ & $4$ & $7$ & $11$ & $28$\\
$K_{5, 6}$ & $5$ & $6$ & $11$ & $30$\\
$K_{1, 11}$ & $1$ & $11$ & $12$ & $11$\\
$K_{2, 10}$ & $2$ & $10$ & $12$ & $20$\\
$K_{3, 9}$ & $3$ & $9$ & $12$ & $27$\\
$bp(3, 9)$ & $3$ & $9$ & $12$ & $18$\\
$K_{4, 8}$ & $4$ & $8$ & $12$ & $32$\\
$bp(4, 8)_1$ & $4$ & $8$ & $12$ & $24$\\
$bp(4, 8)_2$ & $4$ & $8$ & $12$ & $16$\\
$K_{5, 7}$ & $5$ & $7$ & $12$ & $35$\\
$K_{1, 12}$ & $1$ & $12$ & $13$ & $12$\\
$K_{2, 11}$ & $2$ & $11$ & $13$ & $22$\\
$K_{3, 10}$ & $3$ & $10$ & $13$ & $30$\\
$K_{4, 9}$ & $4$ & $9$ & $13$ & $36$\\
$K_{5, 8}$ & $5$ & $8$ & $13$ & $40$\\
$K_{6, 7}$ & $6$ & $7$ & $13$ & $42$\\
$K_{1, 13}$ & $1$ & $13$ & $14$ & $13$\\
$K_{2, 12}$ & $2$ & $12$ & $14$ & $24$\\
$K_{3, 11}$ & $3$ & $11$ & $14$ & $33$\\
\hline
\end{tabular}
\qquad
\begin{tabular}{lcccc}
\hline
& $m_1$ & $m_2$ & $v$ & $e$\\
\hline
$K_{4, 10}$ & $4$ & $10$ & $14$ & $40$\\
$K_{5, 9}$ & $5$ & $9$ & $14$ & $45$\\
$K_{6, 8}$ & $6$ & $8$ & $14$ & $48$\\
$bp(6, 8)_1$ & $6$ & $8$ & $14$ & $24$\\
$bp(6, 8)_2$ & $6$ & $8$ & $14$ & $24$\\
$K_{1, 14}$ & $1$ & $14$ & $15$ & $14$\\
$K_{2, 13}$ & $2$ & $13$ & $15$ & $26$\\
$K_{3, 12}$ & $3$ & $12$ & $15$ & $36$\\
$bp(3, 12)$ & $3$ & $12$ & $15$ & $24$\\
$K_{4, 11}$ & $4$ & $11$ & $15$ & $44$\\
$K_{5, 10}$ & $5$ & $10$ & $15$ & $50$\\
$bp(5,10)_1$ & $5$ & $10$ & $15$ & $40$\\
$bp(5,10)_2$ & $5$ & $10$ & $15$ & $30$\\
$bp(5,10)_3$ & $5$ & $10$ & $15$ & $20$\\
$bp(5,10)_4$ & $5$ & $10$ & $15$ & $20$\\
$K_{6, 9}$ & $6$ & $9$ & $15$ & $54$\\
$bp(6, 9)_1$ & $6$ & $9$ & $15$ & $36$\\
$bp(6, 9)_2$ & $6$ & $9$ & $15$ & $36$\\
$bp(6, 9)_3$ & $6$ & $9$ & $15$ & $18$\\
$K_{7, 8}$ & $7$ & $8$ & $15$ & $56$\\
\hline\\ \\
\end{tabular}
\end{center}}
\caption{Edge-transitive graph with at most $15$ vertices}
\end{table}

\clearpage
\begin{figure}[h]
\begin{center}
$rg(12, 4)$\includegraphics[width=1.3in]{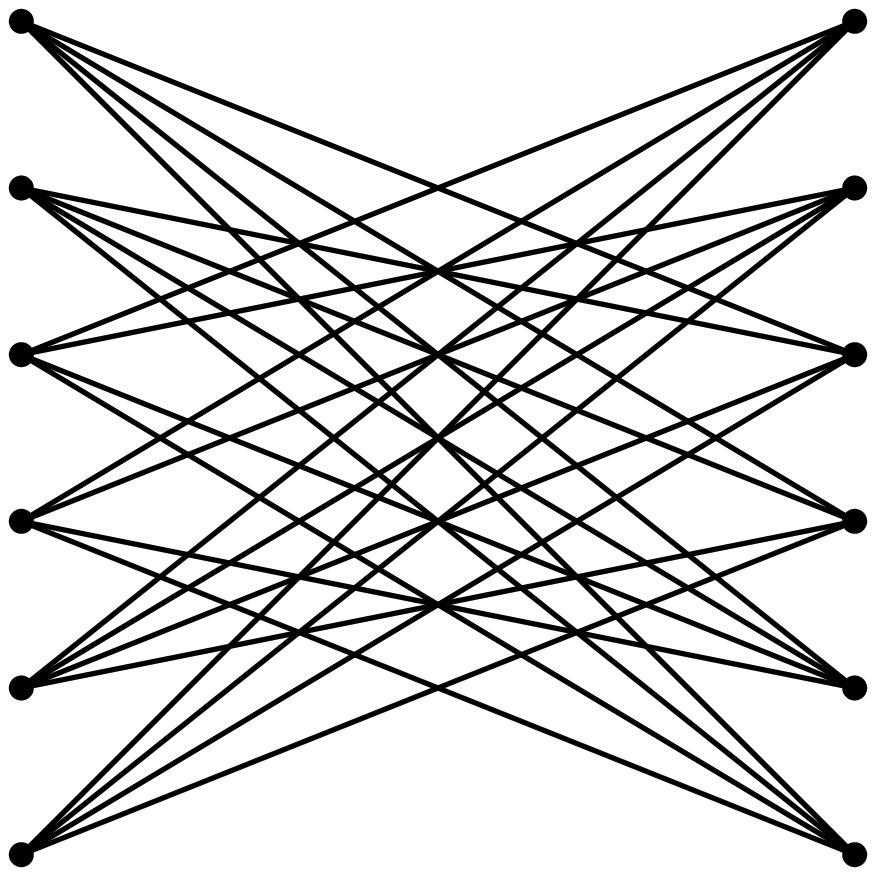} \quad
$rg(12, 5)$\includegraphics[width=1.3in]{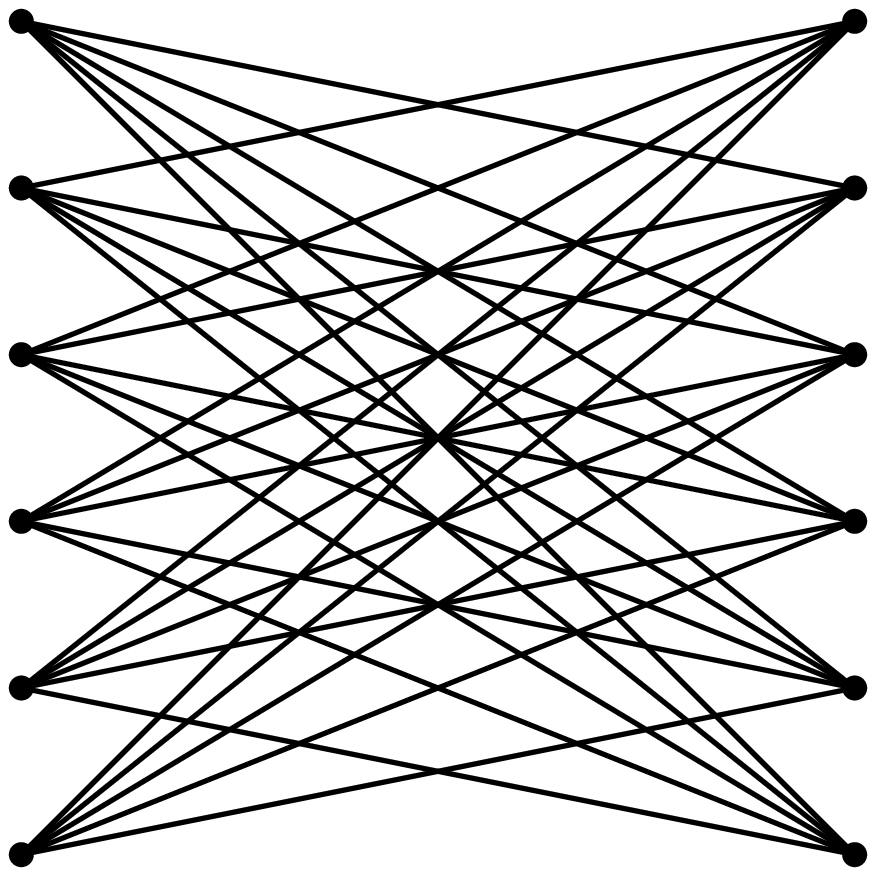} \\
$rg(14, 3)$\includegraphics[width=1.3in]{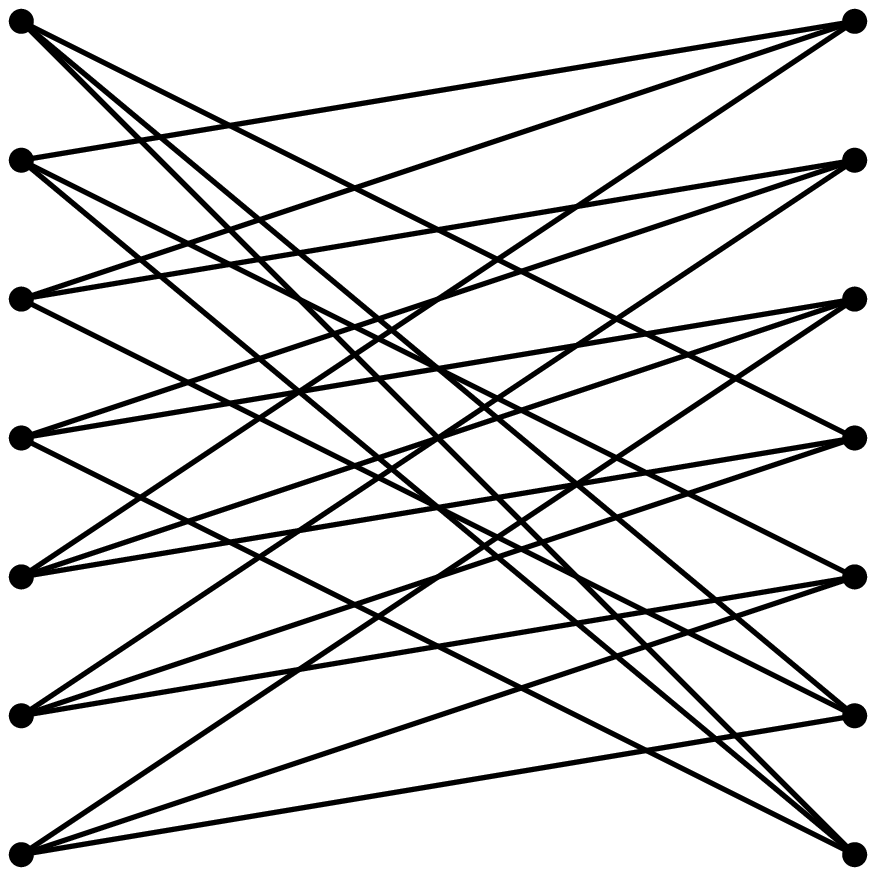} \quad
$rg(14, 4)$\includegraphics[width=1.3in]{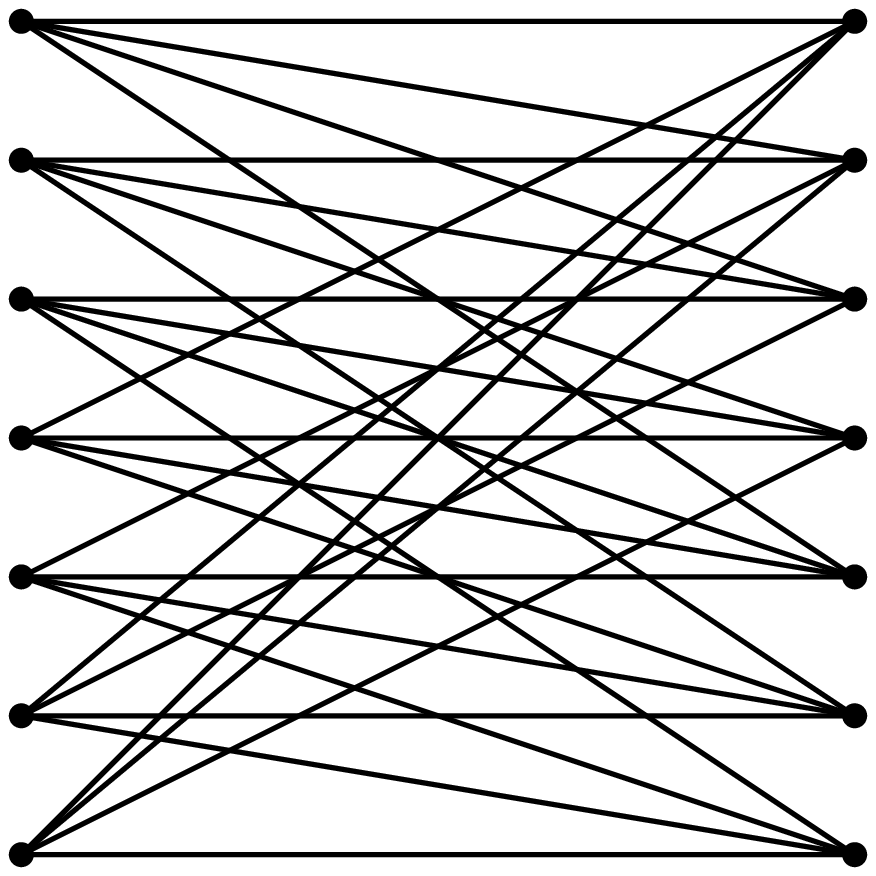} \quad
$rg(15, 4)$\includegraphics[width=1.3in]{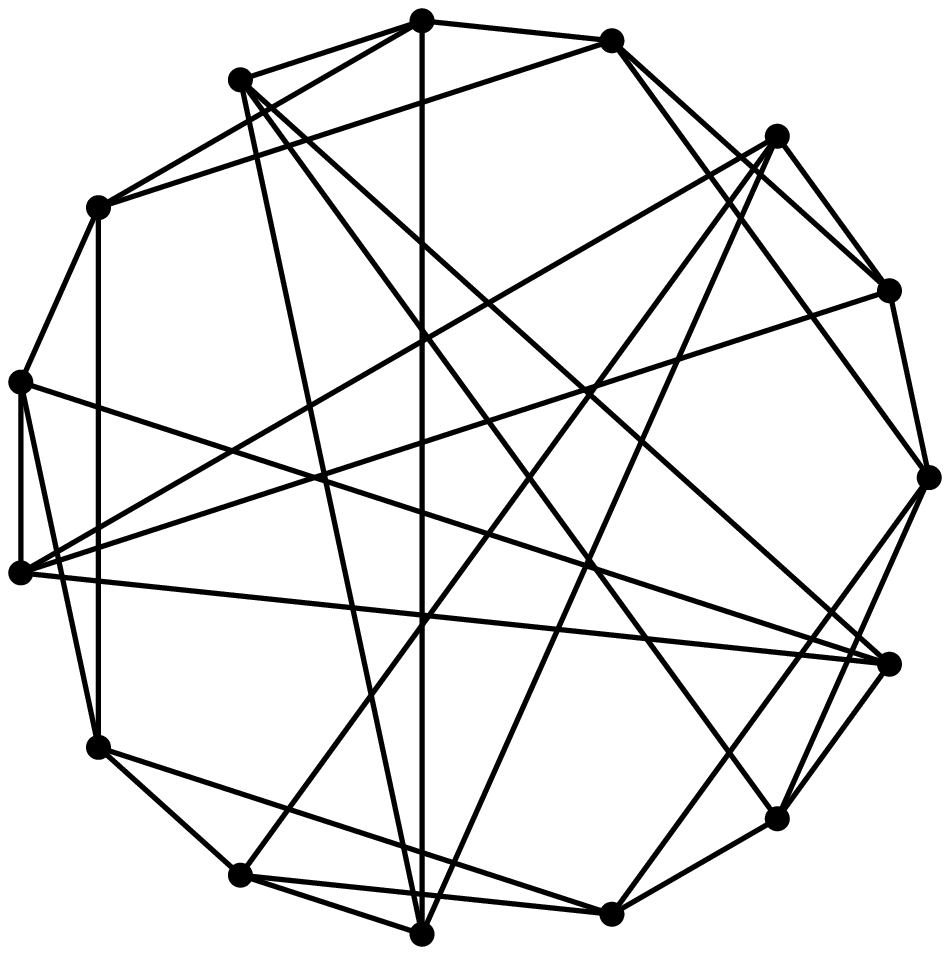} \vspace{0.25in} \\
$bp(3, 6)$\includegraphics[width=0.8in]{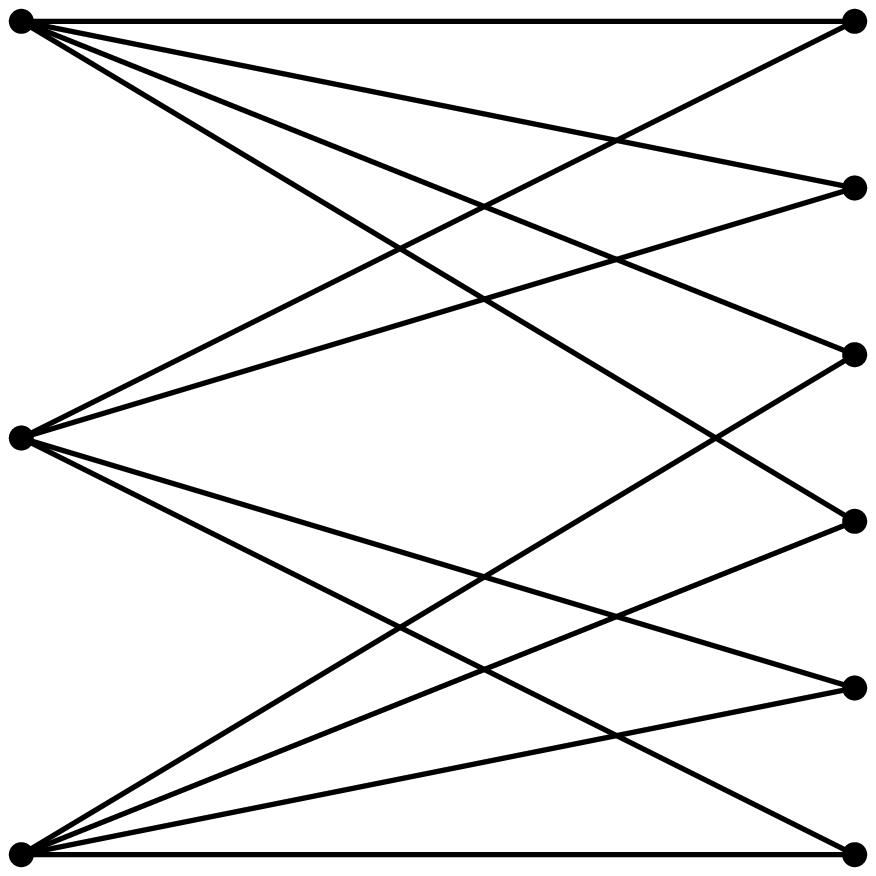} \
$bp(4, 6)$\includegraphics[width=0.8in]{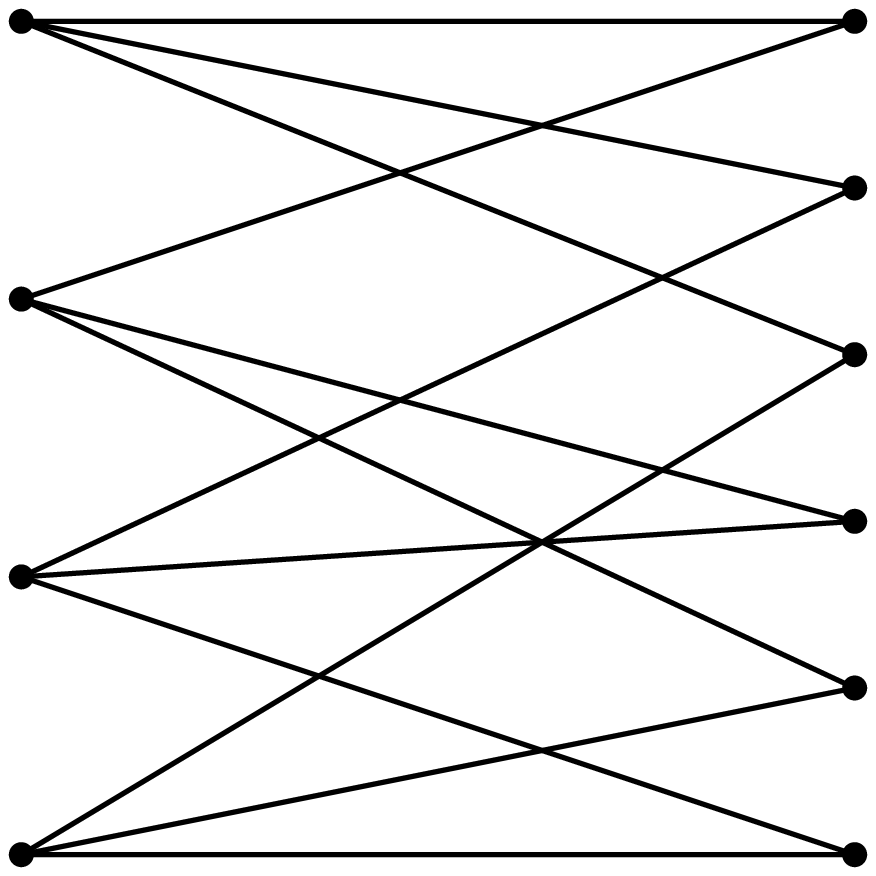} \
$bp(3, 9)$\includegraphics[width=0.8in]{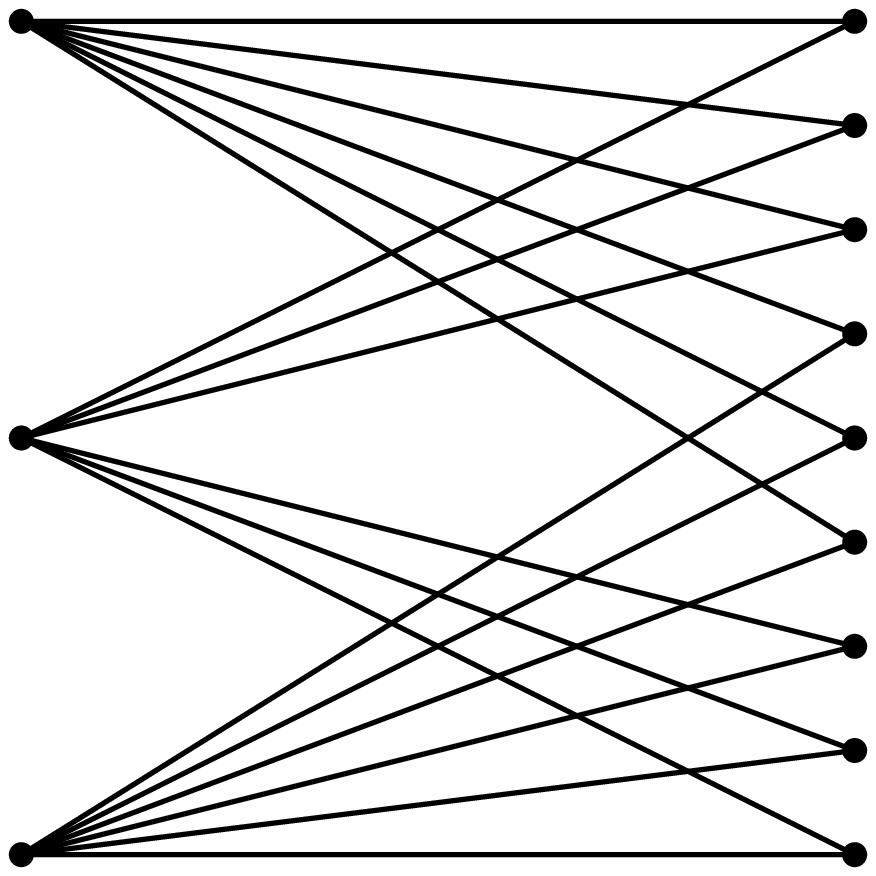} \
$bp(4, 8)_1$\includegraphics[width=0.8in]{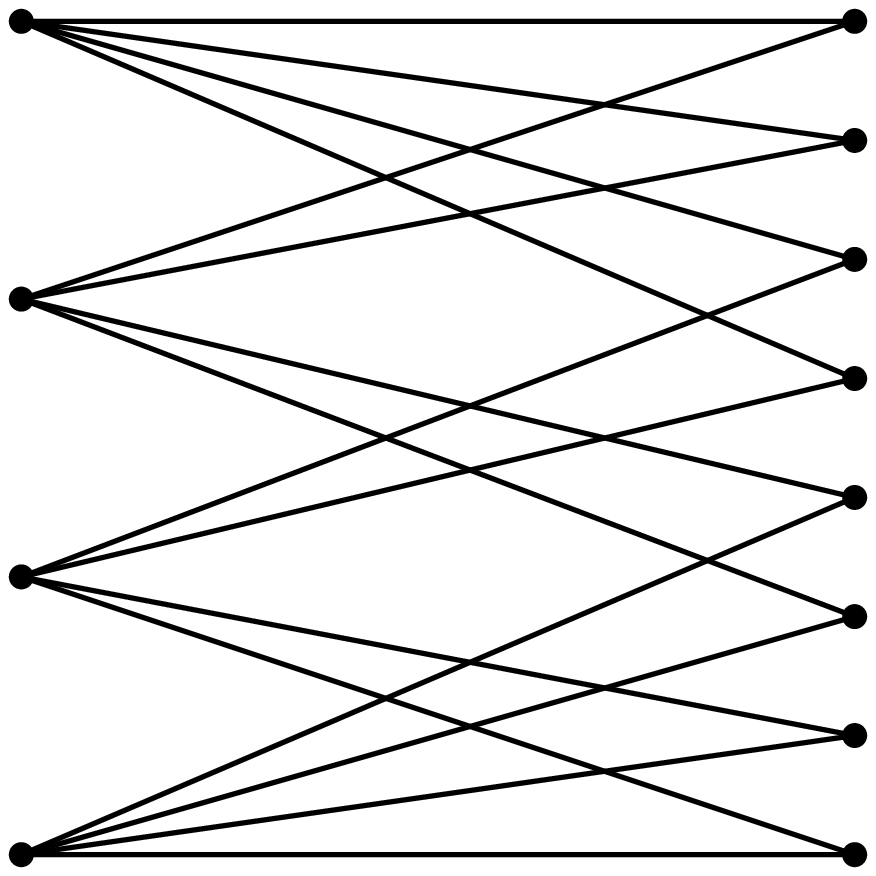} \\
$bp(4, 8)_2$\includegraphics[width=1.1in]{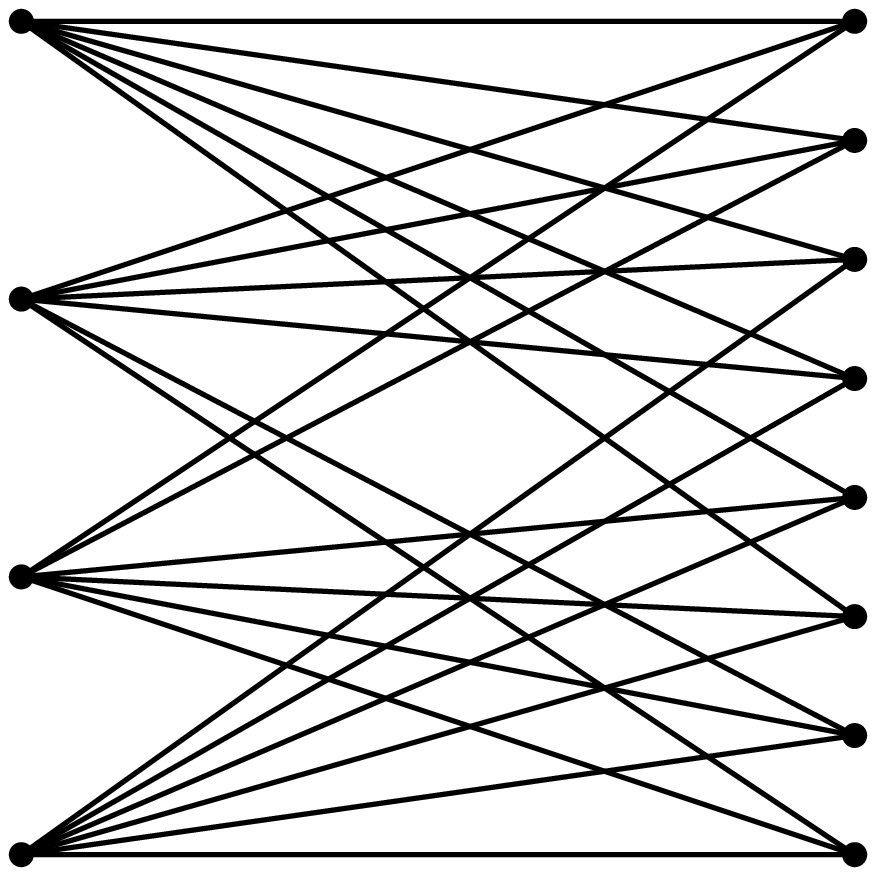} \qquad
$bp(6, 8)_1$\includegraphics[width=1.1in]{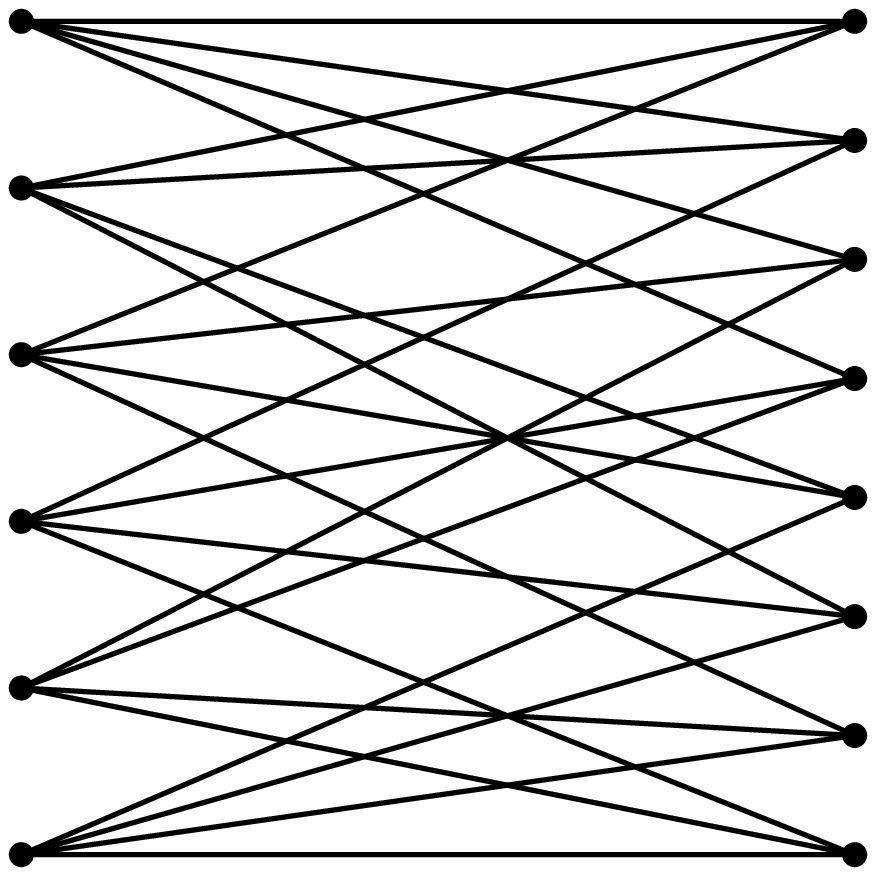} \qquad
$bp(6, 8)_2$\includegraphics[width=1.1in]{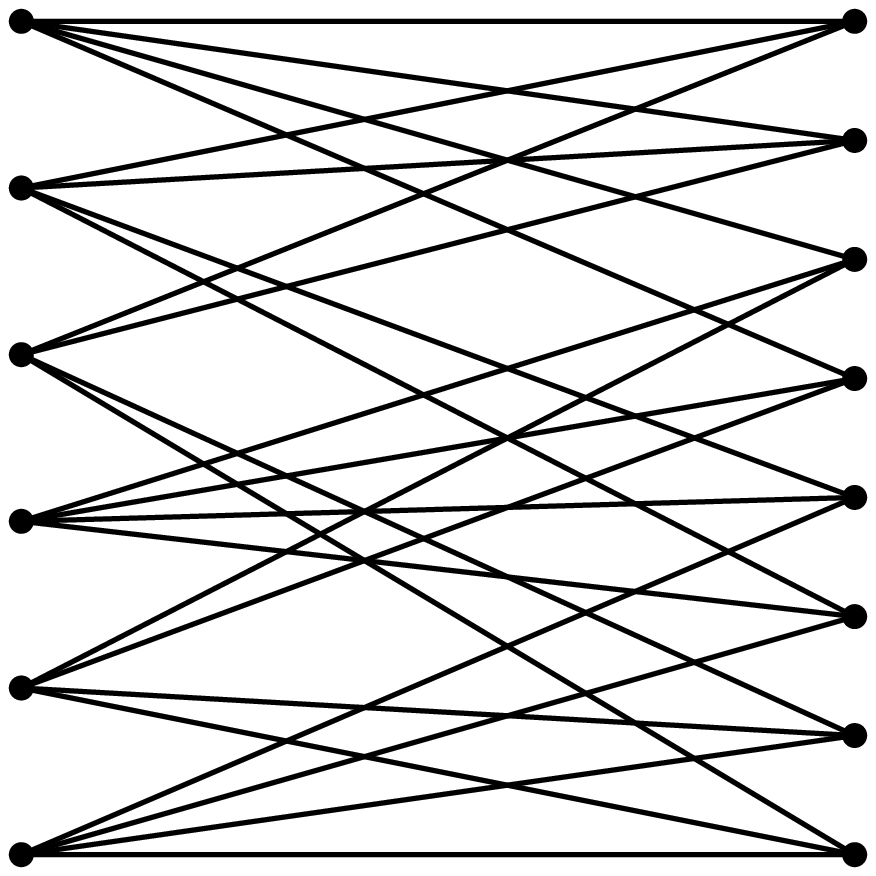} \\
$bp(3, 12)$\includegraphics[width=1.1in]{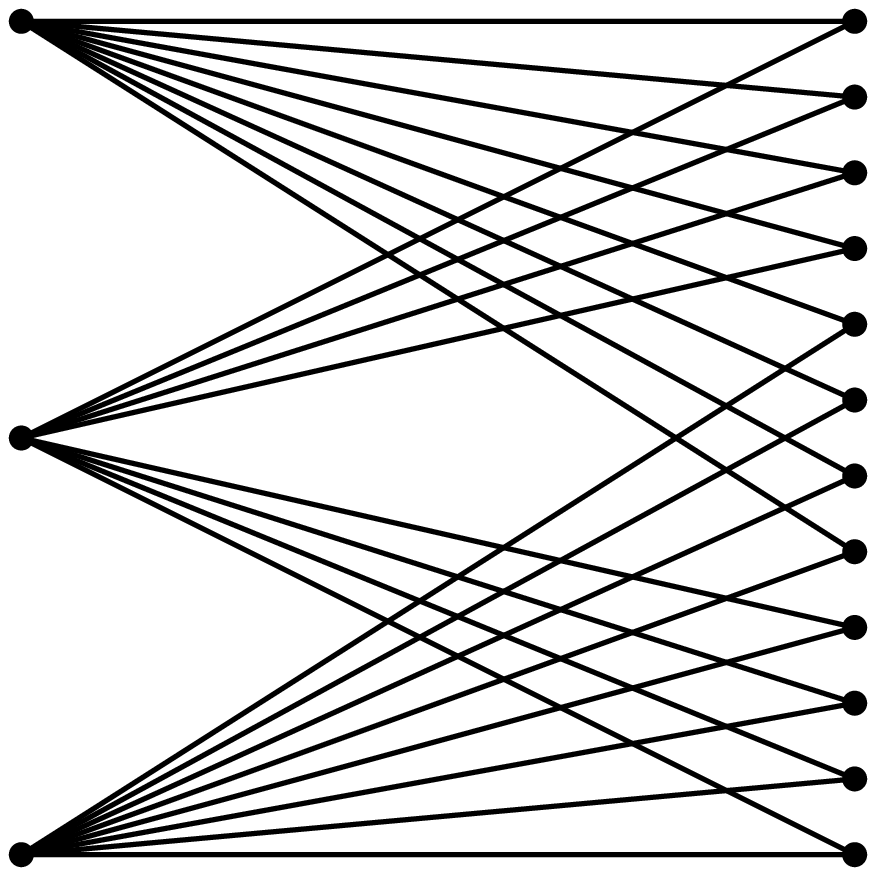} \qquad
$bp(5, 10)_1$\includegraphics[width=1.1in]{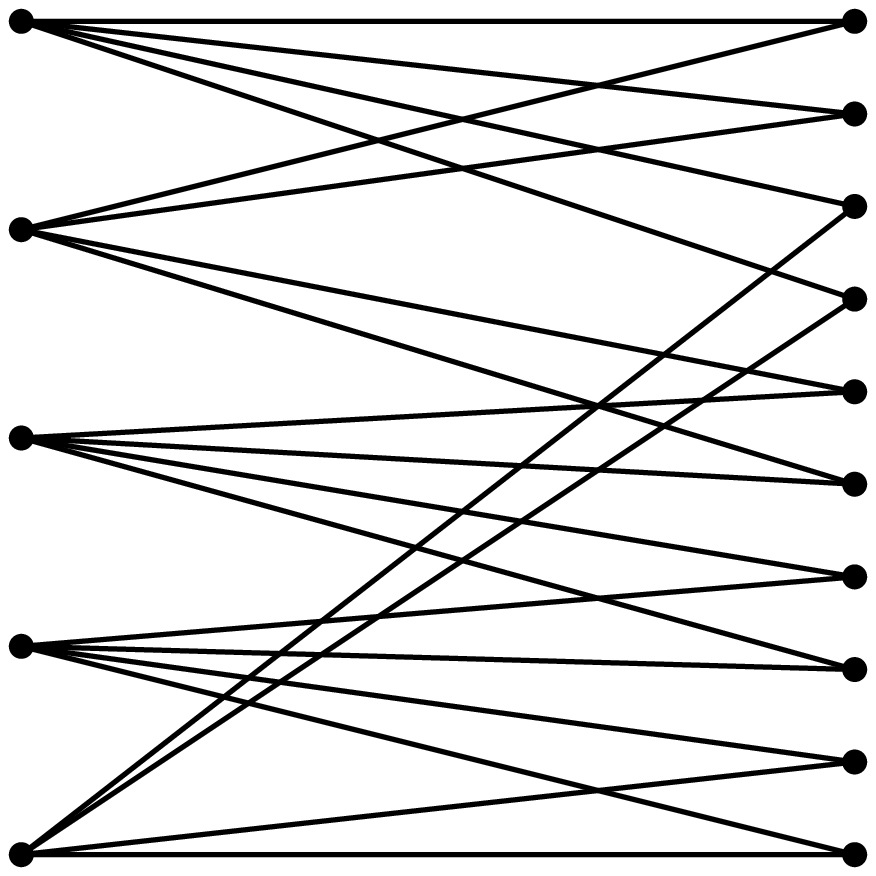} \qquad
$bp(5, 10)_2$\includegraphics[width=1.1in]{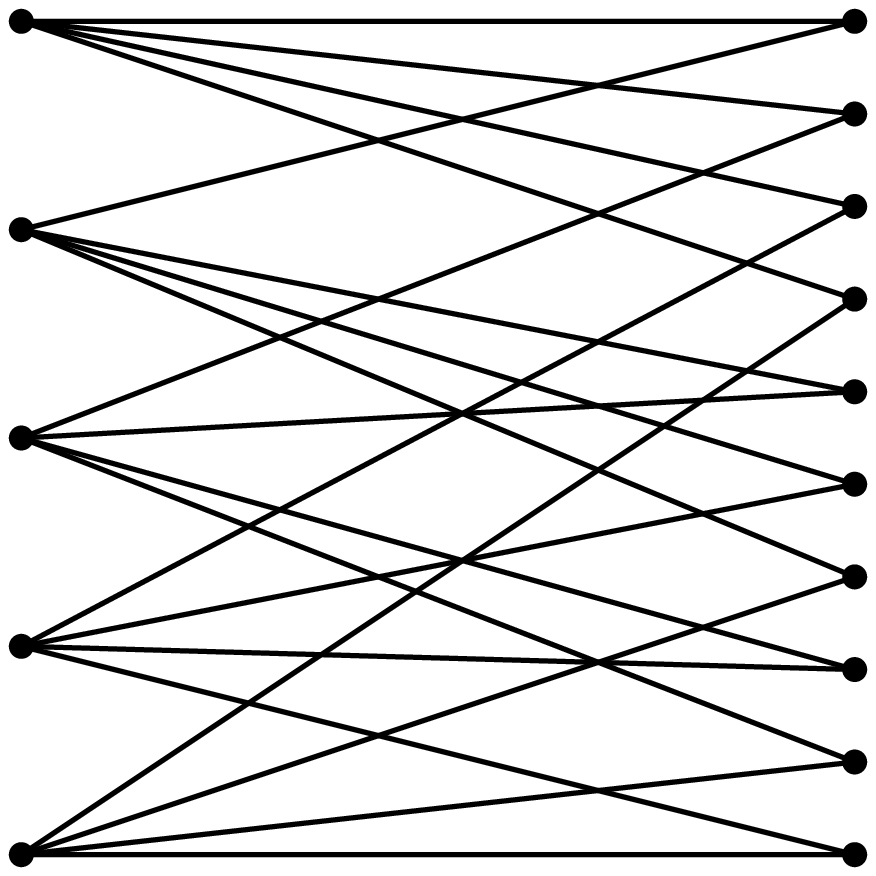} \\
$bp(5, 10)_3$\includegraphics[width=1.3in]{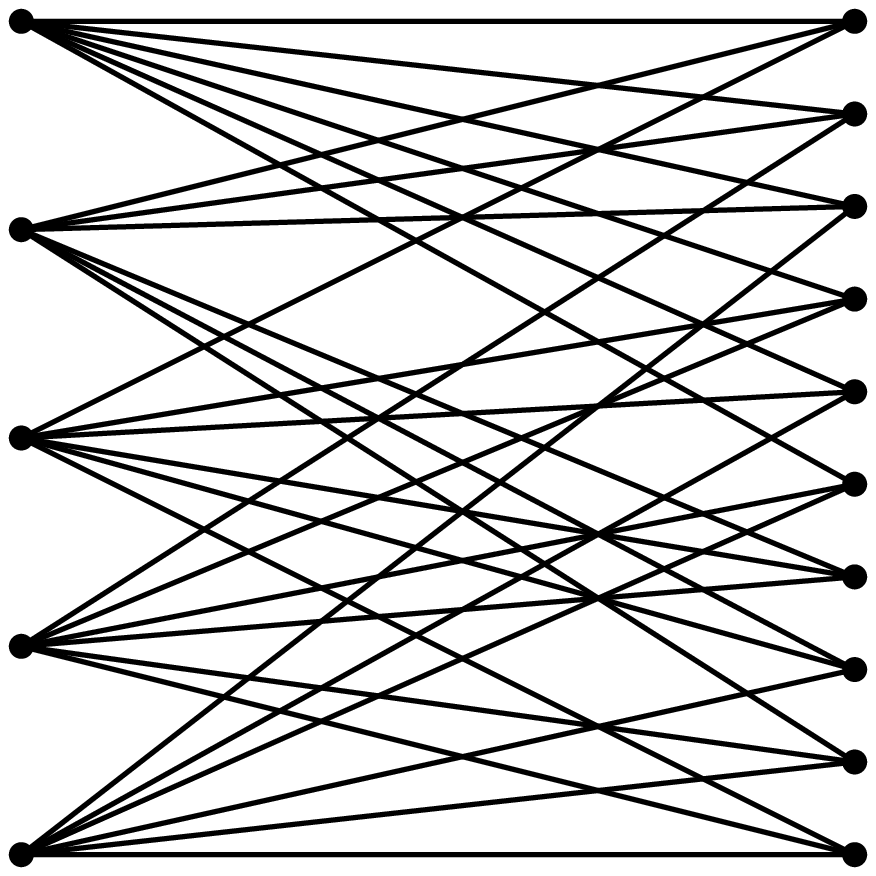} \quad
$bp(5, 10)_4$\includegraphics[width=1.3in]{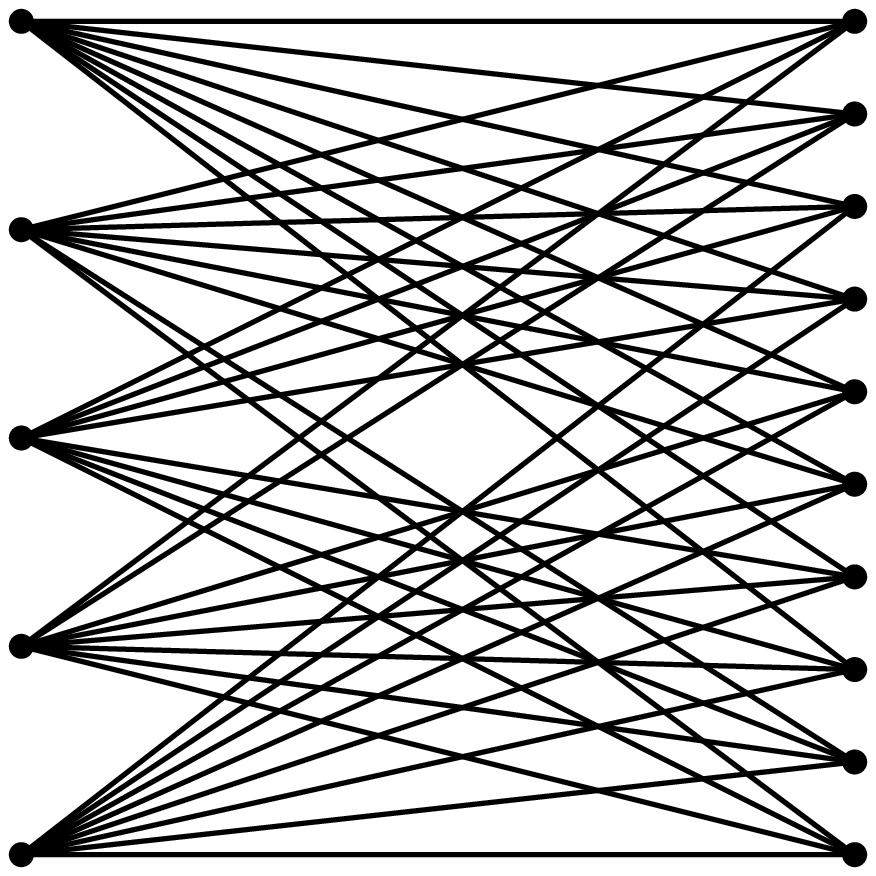} \\
$bp(6, 9)_1$\includegraphics[width=1.3in]{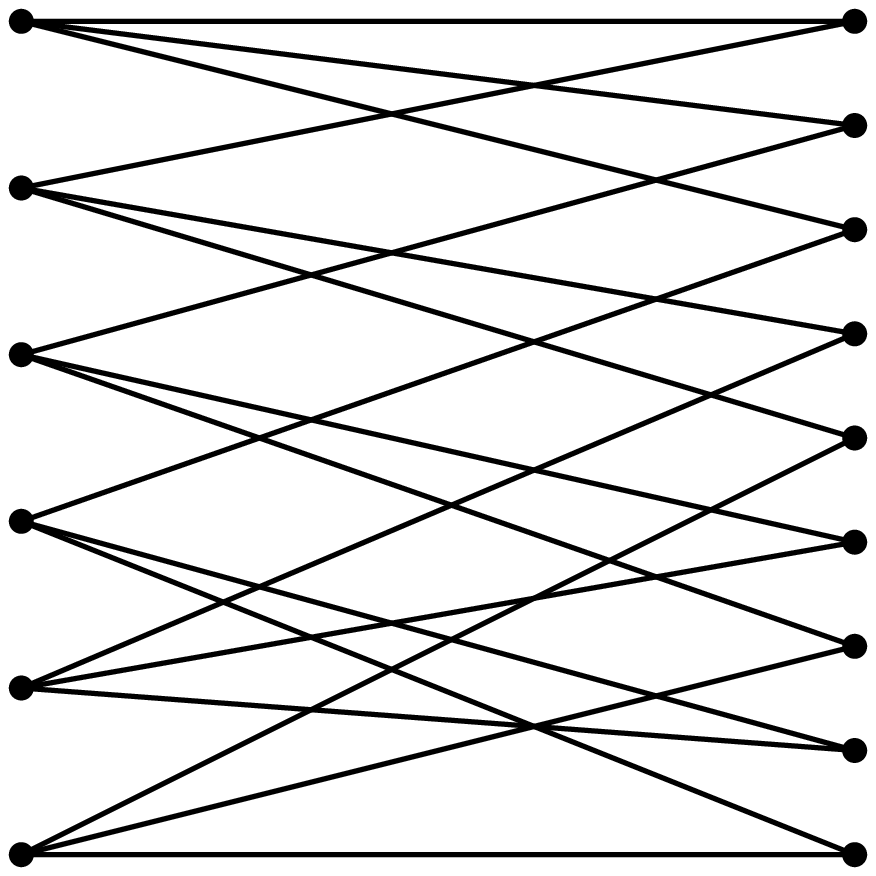} \quad
$bp(6, 9)_2$\includegraphics[width=1.3in]{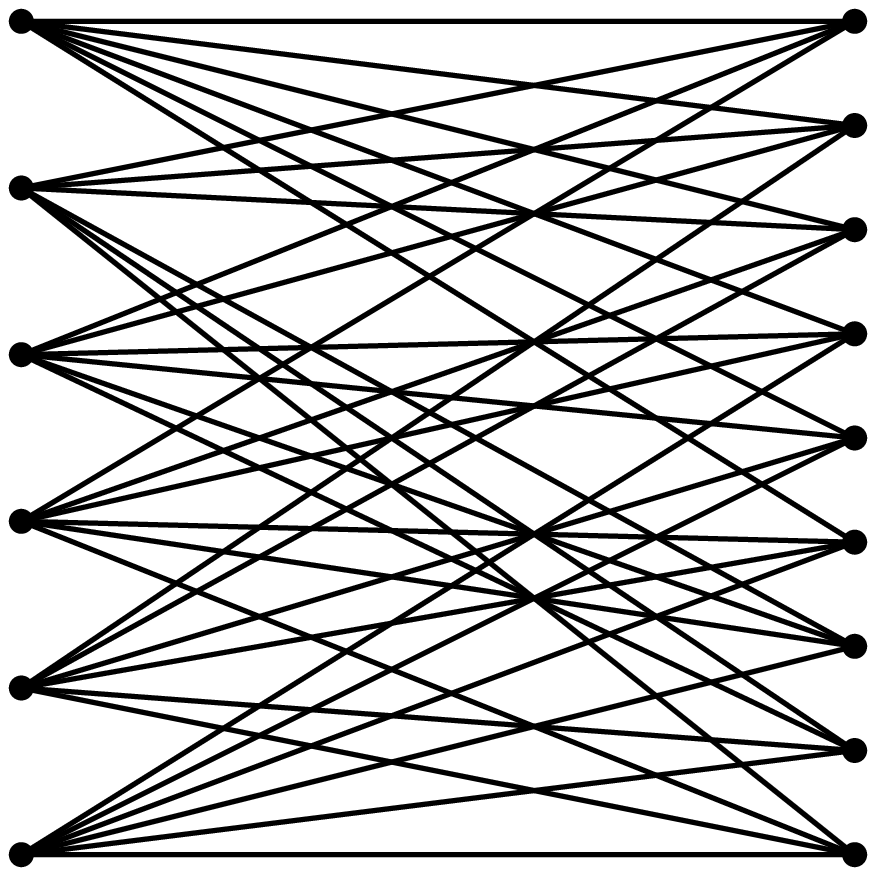} \quad
$bp(6, 9)_3$\includegraphics[width=1.3in]{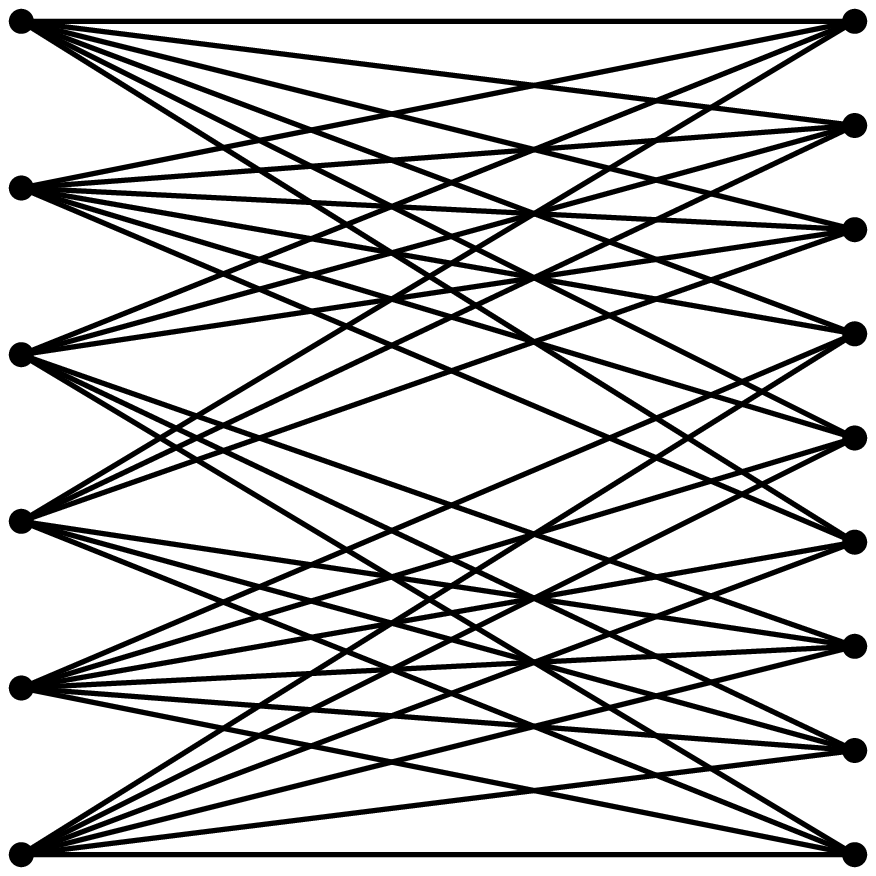} 
\end{center}
\caption{Some edge-transitive graph}\label{fig-etg}
\end{figure}
\clearpage

Now, we consider graphs from which the vectors of crystal lattice are nondegenerate when we construct vectors with the method in Section \ref{sec-sr-lattice}. For the $(d, n, s, t)$-configurations of the vectors of the constructed crystal lattices from such graphs, we have the following table:

\begin{table}[h]
\begin{center}
\begin{tabular}{lccccccc}
\hline
& $v$ & $e$ & \ &  $d$ & $n$ & $s$ & $t$\\
\hline
$K_4$ & $4$ & $6$ && $3$ & $12$ & $4$ & $3$\\
$K_5$ & $5$ & $10$ && $6$ & $20$ & $4$ & $3$\\
$K_{3, 3}$ & $6$ & $9$ && $4$ & $18$ & $5$ & $3$\\
$CP(3)$ & $6$ & $12$ && $7$ & $24$ & $9$ & $3$\\
$K_6$ & $6$ & $15$ && $10$ & $30$ & $4$ & $3$\\
$K_{3, 4}$ & $7$ & $12$ && $6$ & $24$ & $7$ & $3$\\
$K_7$ & $7$ & $21$ && $15$ & $42$ & $4$ & $3$\\
$H(3, 2)$ & $8$ & $12$ && $5$ & $24$ & $9$ & $3$\\
$K_{3, 5}$ & $8$ & $15$ && $8$ & $30$ & $7$ & $3$\\
$K_{4, 4}$ & $8$ & $16$ && $9$ & $32$ & $5$ & $3$\\
$CP(4)$ & $8$ & $24$ && $17$ & $48$ & $10$ & $3$\\
$K_8$ & $8$ & $28$ && $21$ & $56$ & $4$ & $3$\\
$L_2(3)$ & $9$ & $18$ && $10$ & $36$ & $10$ & $3$\\
$K_{3, 6}$ & $9$ & $18$ && $10$ & $36$ & $7$ & $3$\\
$K_{4, 5}$ & $9$ & $20$ && $12$ & $40$ & $7$ & $3$\\
$K_{3, 3, 3}$ & $9$ & $27$ && $19$ & $54$ & $9$ & $3$\\
$K_9$ & $9$ & $36$ && $28$ & $72$ & $4$ & $3$\\
Petersen & $10$ & $15$ && $6$ & $30$ & $6$ & $3$\\
${Ci}_{10}(1, 3)$ & $10$ & $20$ && $11$ & $40$ & $9$ & $3$\\
${Ci}_{10}(1, 4)$ & $10$ & $20$ && $11$ & $40$ & $9$ & $3$\\
$K_{3, 7}$ & $10$ & $21$ && $12$ & $42$ & $7$ & $3$\\
$K_{4, 6}$ & $10$ & $24$ && $15$ & $48$ & $7$ & $3$\\
$K_{5, 5}$ & $10$ & $25$ && $16$ & $50$ & $5$ & $3$\\
$T(5)$ & $10$ & $30$ && $21$ & $60$ & $10$ & $3$\\
$CP(5)$ & $10$ & $40$ && $31$ & $80$ & $10$ & $3$\\
$K_{10}$ & $10$ & $45$ && $36$ & $90$ & $4$ & $3$\\
$K_{3, 8}$ & $11$ & $24$ && $14$ & $48$ & $7$ & $3$\\
$K_{4, 7}$ & $11$ & $28$ && $18$ & $56$ & $7$ & $3$\\
$K_{5, 6}$ & $11$ & $30$ && $20$ & $60$ & $7$ & $3$\\
$K_{11}$ & $11$ & $55$ && $45$ & $110$ & $4$ & $3$\\
${Ci}_{12}(1, 5)$ & $12$ & $24$ && $13$ & $48$ & $9$ & $3$\\
$rg(12, 4)$ & $12$ & $24$ && $13$ & $48$ & $16$ & $3$\\
$bp(4, 8)_2$ & $12$ & $24$ && $13$ & $48$ & $15$ & $3$\\
$K_{3, 9}$ & $12$ & $27$ && $16$ & $54$ & $7$ & $3$\\
$rg(12, 5)$ & $12$ & $30$ && $19$ & $60$ & $9$ & $3$\\
$K_{4, 8}$ & $12$ & $32$ && $21$ & $64$ & $7$ & $3$\\
$K_{5, 7}$ & $12$ & $35$ && $24$ & $70$ & $7$ & $3$\\
{\footnotesize ${Ci}_{12}(1, 2, 5)$} & $12$ & $36$ && $25$ & $72$ & $21$ & $3$\\
$K_{6, 6}$ & $12$ & $36$ && $25$ & $72$ & $5$ & $3$\\
$K_{4, 4, 4}$ & $12$ & $48$ && $37$ & $96$ & $9$ & $3$\\
$K_{3, 3, 3, 3}$ & $12$ & $54$ && $43$ & $108$ & $10$ & $3$\\
$CP(6)$ & $12$ & $60$ && $49$ & $120$ & $10$ & $3$\\
$K_{12}$ & $12$ & $66$ && $55$ & $132$ & $4$ & $3$\\
\hline
\end{tabular}
\qquad
\begin{tabular}{lccccccc}
\hline
& $v$ & $e$ & \ &  $d$ & $n$ & $s$ & $t$\\
\hline
${Ci}_{13}(1, 5)$ & $13$ & $26$ && $14$ & $52$ & $16$ & $3$\\
$K_{3, 10}$ & $13$ & $30$ && $18$ & $60$ & $7$ & $3$\\
$K_{4, 9}$ & $13$ & $36$ && $24$ & $72$ & $7$ & $3$\\
$P(13)$ & $13$ & $39$ && $27$ & $78$ & $10$ & $3$\\
$K_{5, 8}$ & $13$ & $40$ && $28$ & $80$ & $7$ & $3$\\
$K_{6, 7}$ & $13$ & $42$ && $30$ & $84$ & $7$ & $3$\\
$K_{13}$ & $13$ & $78$ && $66$ & $156$ & $4$ & $3$\\
$rg(14, 3)$ & $14$ & $21$ && $8$ & $42$ & $7$ & $3$\\
$bp(6, 8)_1$ & $14$ & $24$ && $11$ & $48$ & $12$ & $3$\\
$bp(6, 8)_2$ & $14$ & $24$ && $11$ & $48$ & $16$ & $3$\\
${Ci}_{14}(1, 6)$ & $14$ & $28$ && $15$ & $56$ & $9$ & $3$\\
$rg(14, 4)$ & $14$ & $28$ && $15$ & $56$ & $7$ & $3$\\
$K_{3, 11}$ &$14$ & $33$ && $20$ & $66$ & $7$ & $3$\\
$K_{4, 10}$ & $14$ & $40$ && $27$ & $80$ & $7$ & $3$\\
{\scriptsize ${Ci}_{14}(1, 3, 5)$} & $14$ & $42$ && $29$ & $84$ & $9$ & $3$\\
$K_{5, 9}$ & $14$ & $45$ && $32$ & $90$ & $7$ & $3$\\
$K_{6, 8}$ & $14$ & $48$ && $35$ & $96$ & $7$ & $3$\\
$K_{7, 7}$ & $14$ & $49$ && $36$ & $98$ & $5$ & $3$\\
$CP(7)$ & $14$ & $84$ && $71$ & $168$ & $10$ & $3$\\
$K_{14}$ & $14$ & $91$ && $78$ & $182$ & $4$ & $3$\\
${Ci}_{15}(1, 4)$ & $15$ & $30$ && $16$ & $60$ & $21$ & $3$\\
$rg(15, 4)$  & $15$ & $30$ && $16$ & $60$ & $14$ & $3$\\
$bp(5, 10)_3$ & $15$ & $30$ && $16$ & $60$ & $11$ & $3$\\
$K_{3, 12}$ & $15$ & $36$ && $22$ & $72$ & $7$ & $3$\\
$bp(6, 9)_2$ & $15$ & $36$ && $22$ & $72$ & $13$ & $3$\\
$bp(6, 9)_3$ & $15$ & $36$ && $22$ & $72$ & $25$ & $3$\\
$bp(5, 10)_4$ & $15$ & $40$ && $26$ & $80$ & $15$ & $3$\\
$K_{4, 11}$ & $15$ & $44$ && $30$ & $88$ & $7$ & $3$\\
{\footnotesize ${Ci}_{15}(1, 4, 6)$} & $15$ & $45$ && $31$ & $90$ & $9$ & $3$\\
$\overline{T(6)}$ & $15$ & $45$ && $31$ & $90$ & $10$ & $3$\\
$K_{5, 10}$ & $15$ & $50$ && $36$ & $100$ & $7$ & $3$\\
$K_{6, 9}$ & $15$ & $54$ && $40$ & $108$ & $7$ & $3$\\
$K_{7, 8}$ & $15$ & $56$ && $42$ & $112$ & $7$ & $3$\\
{\tiny ${Ci}_{15}(1, 2, 4, 7)$} & $15$ & $60$ && $46$ & $120$ & $21$ & $3$\\
$T(6)$ & $15$ & $60$ && $46$ & $120$ & $10$ & $3$\\
$K_{5, 5, 5}$ & $15$ & $75$ && $61$ & $150$ & $9$ & $3$\\
$K_{3, \ldots , 3}$ & $15$ & $90$  && $76$ & $180$ & $10$ & $3$\\
$K_{15}$ & $15$ & $105$ && $91$ & $120$ & $4$ & $3$\\
\hline\\ \\ \\ \\ \\
\end{tabular}
\end{center}\quad
\caption{Edge-transitive graph with at most $15$ vertices}
\end{table}

\newpage

\section{Not edge-transitive graphs}\label{sec-net}

In this section, we consider graphs which are not regular and not edge-transitive, from which we construct vectors of crystal lattices of the same norm. The following graphs are all of such graphs with at most $9$ vertices which are irreducible and nondegenerate when we construct vectors with the method in Section \ref{sec-sr-lattice}.

\begin{figure}[htbp]
\begin{center}
$netg(6, 10)$\includegraphics[width=1in]{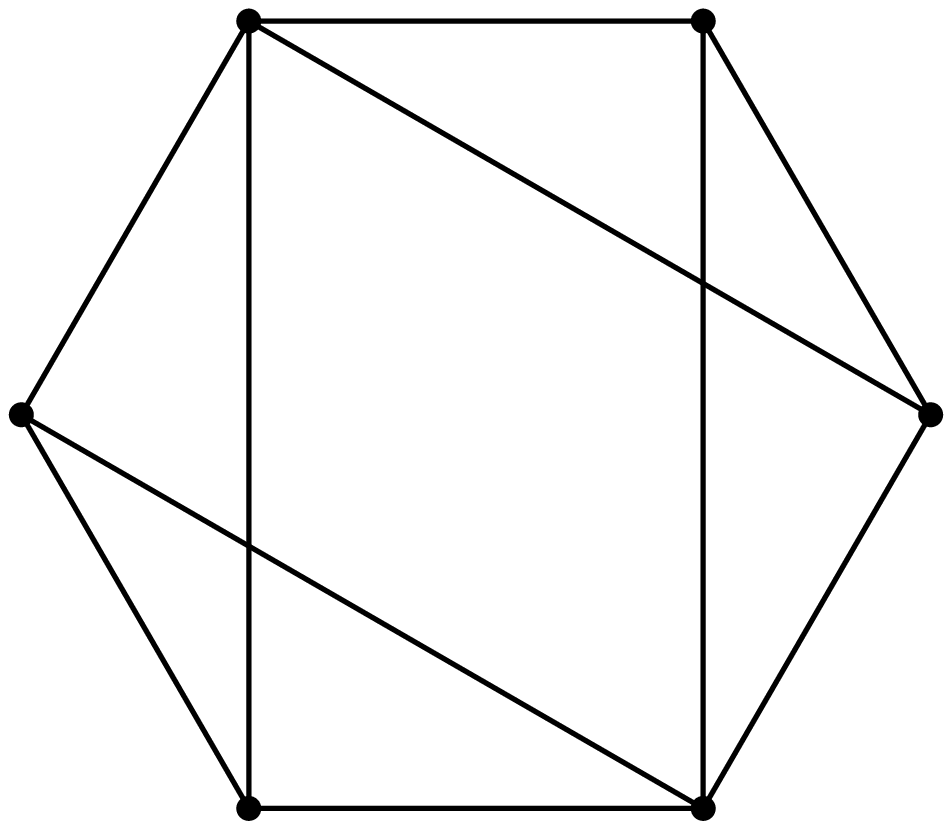} \quad
$netg(8, 14)_1$\includegraphics[width=1.1in]{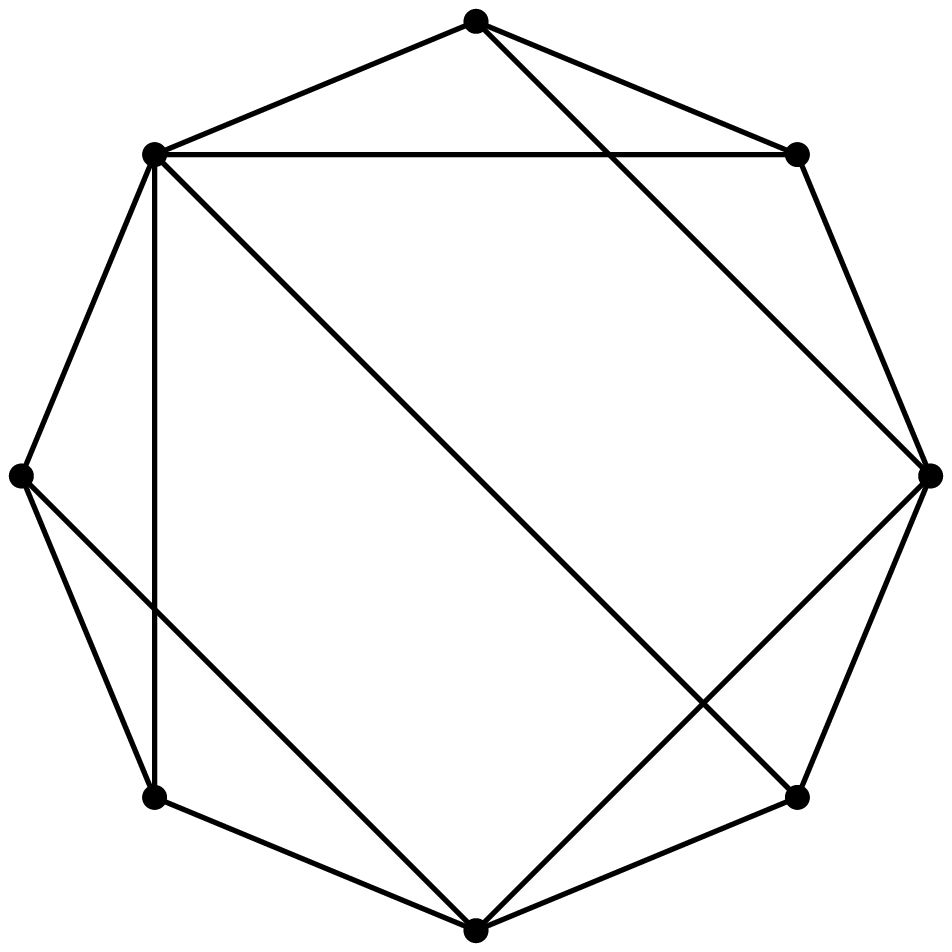} \quad
$netg(8, 14)_2$\includegraphics[width=1.1in]{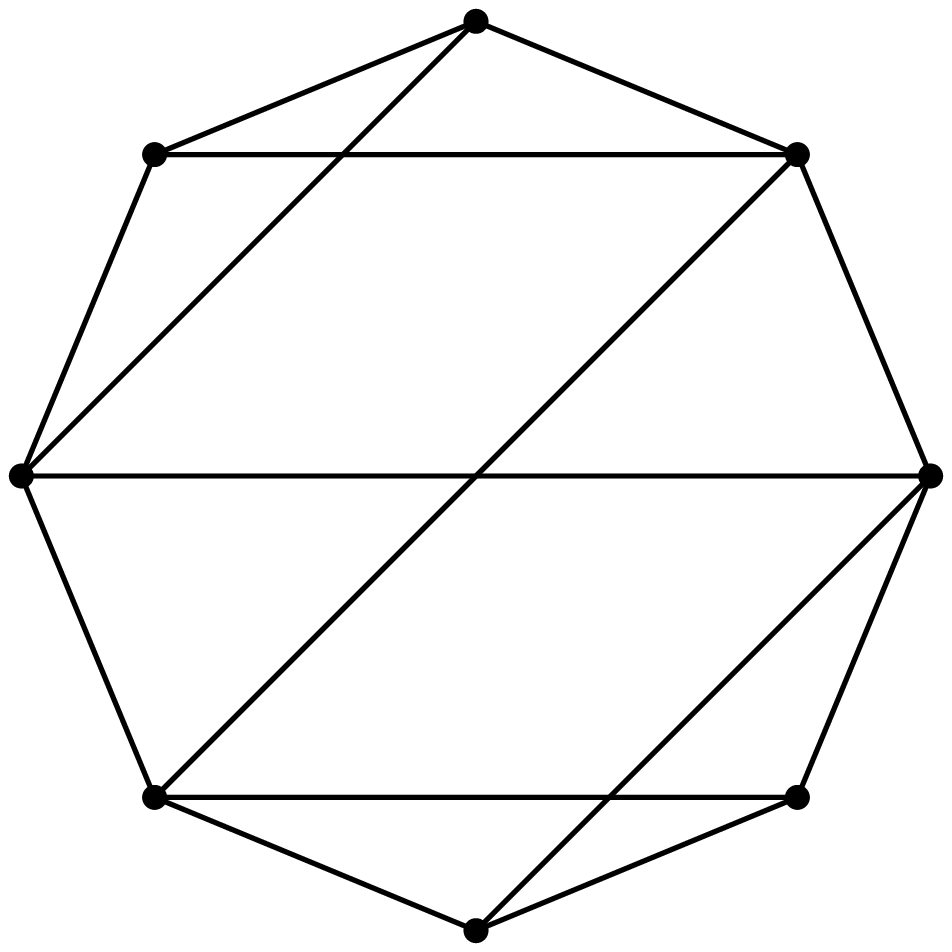}\\
\quad\\
$netg(9, 16)$\includegraphics[width=1.3in]{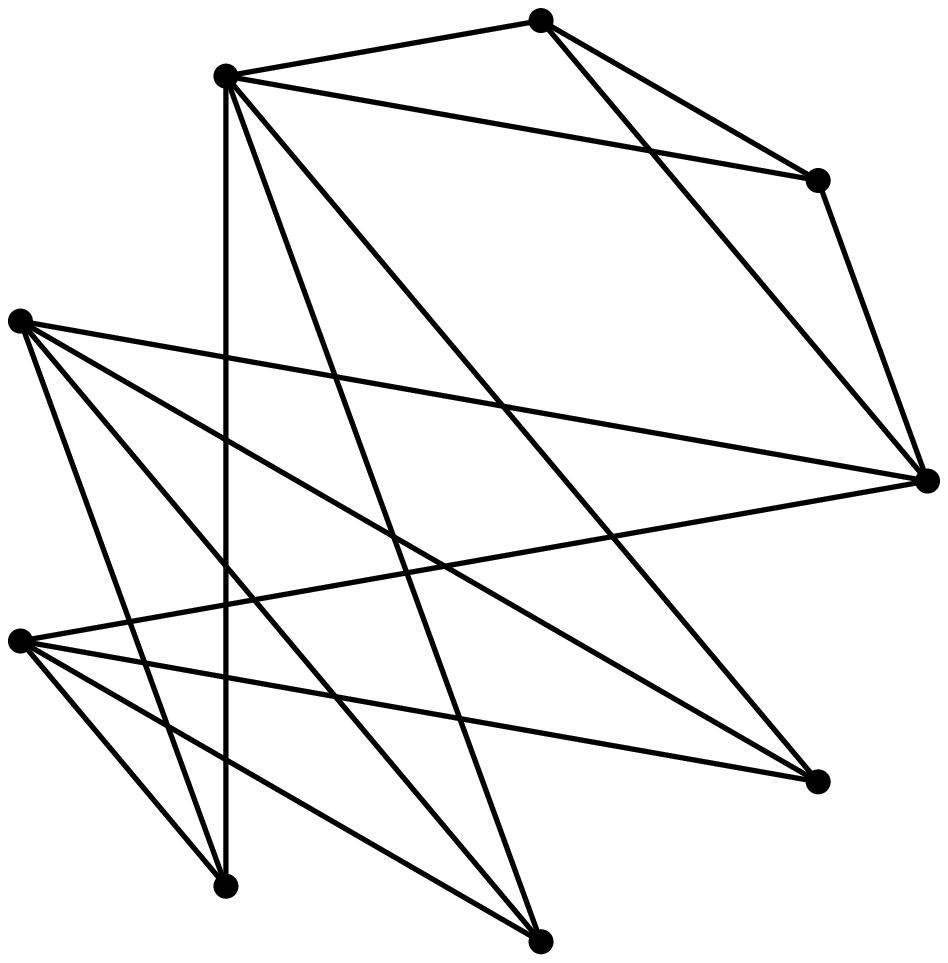} \qquad
$netg(9, 24)$\includegraphics[width=1.3in]{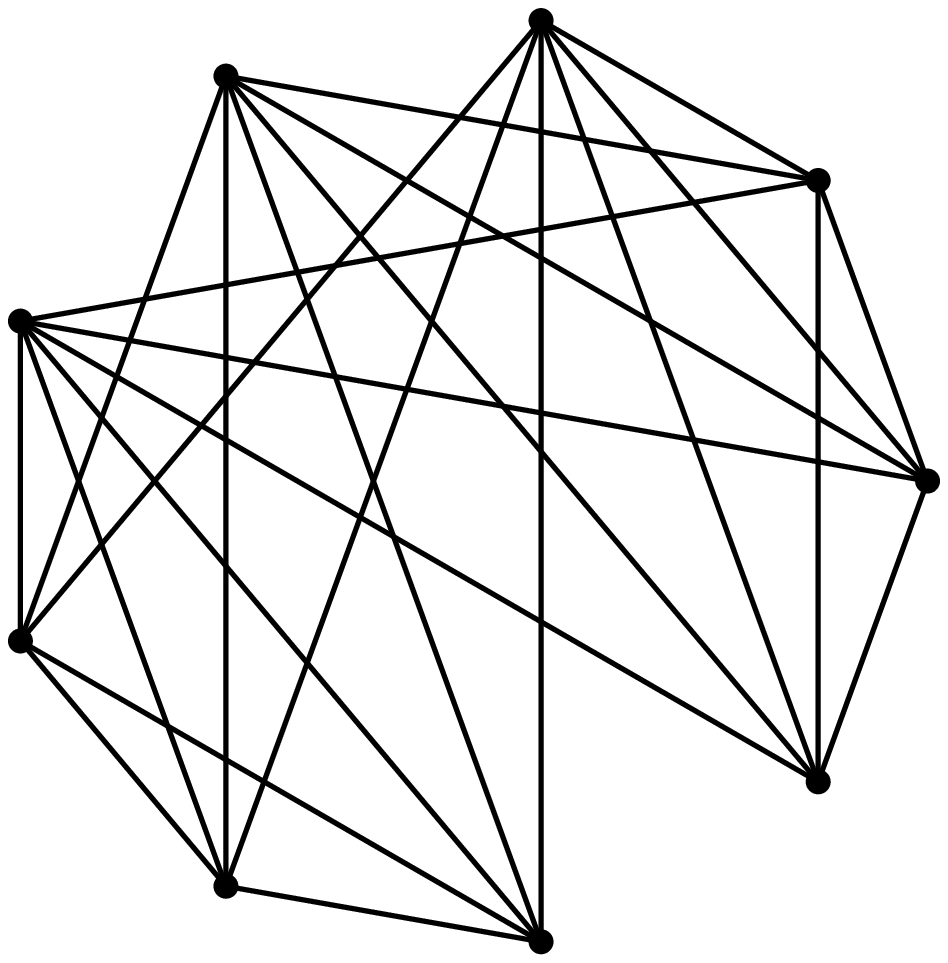}
\end{center}
\caption{Not edge-transitive graphs with at most $9$ vertices}\label{fig-net}
\end{figure}

For the $(d, n, s, t)$-configurations of the vectors of the constructed crystal lattices from the above graphs, we have the following table:

\begin{table}[h]
\begin{center}
\begin{tabular}{lccccccc}
\hline
& $v$ & $e$ & \ &  $d$ & $n$ & $s$ & $t$\\
\hline
$netg(6, 10)$ & $6$ & $10$ && $5$ & $20$ & $6$ & $3$\\
$netg(8, 14)_1$ & $8$ & $14$ && $7$ & $28$ & $8$ & $3$\\
$netg(8, 14)_2$ & $8$ & $14$ && $7$ & $28$ & $6$ & $3$\\
$netg(9, 16)$ & $9$ & $16$ && $8$ & $32$ & $12$ & $3$\\
$netg(9, 24)$ & $9$ & $24$ && $16$ & $48$ & $8$ & $3$\\
\hline
\end{tabular}
\end{center}\quad
\caption{Not edge-transitive graphs with at most $9$ vertices}
\end{table}

\newpage

\section{Graphs with multiple edges}\label{sec-mul-g}
We consider new graphs from finite graphs to change simple edges to $m$-fold multiple edges. Most simple example is the diamond crystal lattice from the graph with two vertices and multiple edge.

\subsection{Graphs with multiple edges from regular polygons}\label{sec-mul-2_reg_poly}
If we construct crystal lattice from $k$-gon, then the clattice degenerate into the $1$-dimensional standard crystal lattice. However, if we consider $k$-gon with multiple edges, then the constructed crystal lattice does not degenerate. For the $(d, n, s, t)$-configurations of the vectors of the constructed crystal lattices, we have the following table:\\

\paragraph{Graphs with double edges}\quad\vspace{-0.05in}

\begin{table}[h]
\begin{center}
\begin{tabular}{ccccc}
\hline
$k$ & $d$ & $n$ & $s$ & $t$\\
\hline
$3$ & $4$ & $12$ & $5$ & $3$\\
$4$ & $5$ & $16$ & $5$ & $3$\\
$5$ & $6$ & $20$ & $5$ & $3$\\
$6$ & $7$ & $24$ & $5$ & $3$\\
$7$ & $8$ & $28$ & $5$ & $3$\\
$8$ & $9$ & $32$ & $5$ & $3$\\
$9$ & $10$ & $36$ & $5$ & $3$\\
$10$ & $11$ & $40$ & $5$ & $3$\\
$11$ & $12$ & $44$ & $5$ & $3$\\
$12$ & $13$ & $48$ & $5$ & $3$\\
\hline
\end{tabular}
\quad
\begin{tabular}{ccccc}
\hline
$k$ & $d$ & $n$ & $s$ & $t$\\
\hline
$13$ & $14$ & $52$ & $5$ & $3$\\
$14$ & $15$ & $56$ & $5$ & $3$\\
$15$ & $16$ & $60$ & $5$ & $3$\\
$16$ & $17$ & $64$ & $5$ & $3$\\
$17$ & $18$ & $68$ & $5$ & $3$\\
$18$ & $19$ & $72$ & $5$ & $3$\\
$19$ & $20$ & $76$ & $5$ & $3$\\
$20$ & $21$ & $80$ & $5$ & $3$\\
$21$ & $22$ & $84$ & $5$ & $3$\\
$22$ & $23$ & $88$ & $5$ & $3$\\
\hline
\end{tabular}
\quad
\begin{tabular}{ccccc}
\hline
$k$ & $d$ & $n$ & $s$ & $t$\\
\hline
$23$ & $24$ & $92$ & $5$ & $3$\\
$24$ & $25$ & $96$ & $5$ & $3$\\
$25$ & $26$ & $100$ & $5$ & $3$\\
$26$ & $27$ & $104$ & $5$ & $3$\\
$27$ & $28$ & $108$ & $5$ & $3$\\
$28$ & $29$ & $112$ & $5$ & $3$\\
$29$ & $30$ & $116$ & $5$ & $3$\\
$30$ & $31$ & $120$ & $5$ & $3$\\
$31$ & $32$ & $124$ & $5$ & $3$\\
$32$ & $33$ & $128$ & $5$ & $3$\\
\hline
\end{tabular}
\quad
\begin{tabular}{ccccc}
\hline
$k$ & $d$ & $n$ & $s$ & $t$\\
\hline
$33$ & $34$ & $132$ & $5$ & $3$\\
$34$ & $35$ & $136$ & $5$ & $3$\\
$35$ & $36$ & $140$ & $5$ & $3$\\
$36$ & $37$ & $144$ & $5$ & $3$\\
$37$ & $38$ & $148$ & $5$ & $3$\\
$38$ & $39$ & $152$ & $5$ & $3$\\
$39$ & $40$ & $156$ & $5$ & $3$\\
$40$ & $41$ & $160$ & $5$ & $3$\\
$41$ & $42$ & $164$ & $5$ & $3$\\
$42$ & $43$ & $168$ & $5$ & $3$\\
\hline
\end{tabular}
\end{center}\quad
\caption{Graph with double edges from regular polygon}\vspace{-0.27in}
\end{table}

For every integer $3 \leqslant k \leqslant 100$, by numerical calculation, we can prove that we have the distance set $A(X) = \{ -1, \: \pm (k-1) / (k+1), \: \pm 1 / (k+1) \}$ and that the configuration of the vectors is $(d, n, s, t) = (k + 1, 4 k, 5, 3)$. Furthermore, for each $x \in X$, we have the following table:
\begin{center}
\begin{tabular}{c|ccc}
$a$ & $-1$ & $\pm (k-1) / (k+1)$ & $\pm 1 / (k+1)$\\
\hline
$|A_x(X, a)|$ & $1$ & $1$ & $2 (k-1)$\\
\end{tabular}
\end{center}\quad

\paragraph{Graphs with triple edges}\quad\vspace{-0.05in}

\begin{table}[h]
\begin{center}
\begin{tabular}{ccccc}
\hline
$k$ & $d$ & $n$ & $s$ & $t$\\
\hline
$3$ & $7$ & $18$ & $5$ & $3$\\
$4$ & $9$ & $24$ & $5$ & $3$\\
$5$ & $11$ & $30$ & $5$ & $3$\\
$6$ & $13$ & $36$ & $5$ & $3$\\
$7$ & $15$ & $42$ & $5$ & $3$\\
\hline
\end{tabular}
\quad
\begin{tabular}{ccccc}
\hline
$k$ & $d$ & $n$ & $s$ & $t$\\
\hline
$8$ & $17$ & $48$ & $5$ & $3$\\
$9$ & $19$ & $54$ & $5$ & $3$\\
$10$ & $21$ & $60$ & $5$ & $3$\\
$11$ & $23$ & $66$ & $5$ & $3$\\
$12$ & $25$ & $72$ & $5$ & $3$\\
\hline
\end{tabular}
\quad
\begin{tabular}{ccccc}
\hline
$k$ & $d$ & $n$ & $s$ & $t$\\
\hline
$13$ & $27$ & $78$ & $5$ & $3$\\
$14$ & $29$ & $84$ & $5$ & $3$\\
$15$ & $31$ & $90$ & $5$ & $3$\\
$16$ & $33$ & $96$ & $5$ & $3$\\
$17$ & $35$ & $102$ & $5$ & $3$\\
\hline
\end{tabular}
\quad
\begin{tabular}{ccccc}
\hline
$k$ & $d$ & $n$ & $s$ & $t$\\
\hline
$18$ & $37$ & $108$ & $5$ & $3$\\
$19$ & $39$ & $114$ & $5$ & $3$\\
$20$ & $41$ & $120$ & $5$ & $3$\\
$21$ & $43$ & $126$ & $5$ & $3$\\
$22$ & $45$ & $132$ & $5$ & $3$\\
\hline
\end{tabular}
\end{center}\quad
\caption{Graph with triple edges from regular polygon}\vspace{-0.27in}
\end{table}

For every integer $3 \leqslant k \leqslant 50$, by numerical calculation, we can prove that we have the distance set $A(X) = \{ -1, \: \pm (k-1) / (2k+1), \: \pm 1 / (2k+1) \}$ and that the configuration of the vectors is $(d, n, s, t) = (2 k + 1, 6 k, 5, 3)$. Furthermore, for each $x \in X$, we have the following table:
\begin{center}
\begin{tabular}{c|ccc}
$a$ & $-1$ & $\pm (k-1) / (2k+1)$ & $\pm 1 / (2k+1)$\\
\hline
$|A_x(X, a)|$ & $1$ & $2$ & $3 (k-1)$\\
\end{tabular}
\end{center}\quad

\paragraph{Graphs with quadruple edges}\quad\vspace{-0.05in}

\begin{table}[h]
\begin{center}
\begin{tabular}{ccccc}
\hline
$k$ & $d$ & $n$ & $s$ & $t$\\
\hline
$3$ & $10$ & $24$ & $5$ & $3$\\
$4$ & $13$ & $32$ & $5$ & $3$\\
$5$ & $16$ & $40$ & $5$ & $3$\\
$6$ & $19$ & $48$ & $5$ & $3$\\
$7$ & $22$ & $56$ & $5$ & $3$\\
\hline
\end{tabular}
\quad
\begin{tabular}{ccccc}
\hline
$k$ & $d$ & $n$ & $s$ & $t$\\
\hline
$8$ & $25$ & $64$ & $5$ & $3$\\
$9$ & $28$ & $72$ & $5$ & $3$\\
$10$ & $31$ & $80$ & $5$ & $3$\\
$11$ & $34$ & $88$ & $5$ & $3$\\
$12$ & $37$ & $96$ & $5$ & $3$\\
\hline
\end{tabular}
\quad
\begin{tabular}{ccccc}
\hline
$k$ & $d$ & $n$ & $s$ & $t$\\
\hline
$13$ & $40$ & $104$ & $5$ & $3$\\
$14$ & $43$ & $112$ & $5$ & $3$\\
$15$ & $46$ & $120$ & $5$ & $3$\\
$16$ & $49$ & $128$ & $5$ & $3$\\
$17$ & $52$ & $136$ & $5$ & $3$\\
\hline
\end{tabular}
\quad
\begin{tabular}{ccccc}
\hline
$k$ & $d$ & $n$ & $s$ & $t$\\
\hline
$18$ & $55$ & $144$ & $5$ & $3$\\
$19$ & $58$ & $152$ & $5$ & $3$\\
$20$ & $61$ & $160$ & $5$ & $3$\\
$21$ & $64$ & $168$ & $5$ & $3$\\
$22$ & $67$ & $176$ & $5$ & $3$\\
\hline
\end{tabular}
\end{center}\quad
\caption{Graph with quadruple edges from regular polygon}\vspace{-0.27in}
\end{table}

\newpage

For every integer $3 \leqslant k \leqslant 50$, by numerical calculation, we can prove that we have the distance set $A(X) = \{ -1, \: \pm (k-1) / (3k+1), \: \pm 1 / (3k+1) \}$ and that the configuration of the vectors is $(d, n, s, t) = (3 k + 1, 8 k, 5, 3)$. Furthermore, for each $x \in X$, we have the following table:
\begin{center}
\begin{tabular}{c|ccc}
$a$ & $-1$ & $\pm (k-1) / (3k+1)$ & $\pm 1 / (3k+1)$\\
\hline
$|A_x(X, a)|$ & $1$ & $3$ & $4 (k-1)$\\
\end{tabular}
\end{center}\quad

\paragraph{Graphs with $m$-fold multiple edges}

As a conjecture from above results by numerical calculations, we may have the distance set $A(X) = \{ -1, \: \pm (k-1) / (k (m-1) + 1), \: \pm 1 / (k (m-1) + 1) \}$ and that the configuration of the vectors may be $(d, n, s, t) = (k (m-1) + 1, 2 k m, 5, 3)$. Furthermore, for each $x \in X$, we may have the following table:
\begin{center}
\begin{tabular}{c|ccc}
$a$ & $-1$ & $\pm (k-1) / (k (m-1) + 1)$ & $\pm 1 / (k (m-1) + 1)$\\
\hline
$|A_x(X, a)|$ & $1$ & $m-1$ & $m (k-1)$\\
\end{tabular}
\end{center}\quad

\subsection{Graphs with multiple edges from regular polyhedrons}\label{sec-mul-3_reg_poly}
In this section, we consider graphs with multiple edges from regular polyhedrons (cf. Section \ref{sec-3_reg_poly}). For the $(d, n, s, t)$-configurations of the vectors of the constructed crystal lattices, we have the following table:\vspace{-0.05in}

\begin{table}[h]
\begin{flushleft}
\qquad
\begin{tabular}{ccc}
& \qquad\\
&&\\
Tetrahedron && $\cdots$\\
Hexahedron && $\cdots$\\
Octahedron && $\cdots$\\
Dodecahedron && $\cdots$\\
Icosahedron && $\cdots$
\end{tabular}
\quad
\begin{tabular}{cccc}
\multicolumn{4}{c}{\small Graph with simple edge}\\
\hline
$d$ & $n$ & $s$ & $t$ \\
\hline
$3$ & $12$ & $4$ & $3$\\
$5$ & $24$ & $9$ & $3$\\
$7$ & $24$ & $9$ & $3$\\
$11$ & $60$ & $12$ & $3$\\
$19$ \ & \ $60$ \ & \ $12$ \ & \ $3$\\
\hline
\end{tabular}
\quad $\Leftrightarrow$
\end{flushleft}
\quad\\
\begin{center}
\begin{tabular}{cccc}
\multicolumn{4}{c}{Double edge}\\
\hline
$d$ & $n$ & $s$ & $t$\\
\hline
$9$ & $24$ & $6$ & $3$\\
$17$ & $48$ & $11$ & $3$\\
$19$ & $48$ & $11$ & $3$\\
$41$ & $120$ & $14$ & $3$\\
$49$ & $120$ & $14$ & $3$\\
\hline
\end{tabular} \quad
\begin{tabular}{cccc}
\multicolumn{4}{c}{Triple edge}\\
\hline
$d$ & $n$ & $s$ & $t$\\
\hline
$15$ & $36$ & $6$ & $3$\\
$29$ & $72$ & $11$ & $3$\\
$31$ & $72$ & $11$ & $3$\\
$71$ & $180$ & $14$ & $3$\\
$79$ & $180$ & $14$ & $3$\\
\hline
\end{tabular} \quad
\begin{tabular}{cccc}
\multicolumn{4}{c}{Quadruple edge}\\
\hline
$d$ & $n$ & $s$ & $t$\\
\hline
$21$ & $48$ & $6$ & $3$\\
$41$ & $96$ & $11$ & $3$\\
$43$ & $96$ & $11$ & $3$\\
$101$ & $240$ & $14$ & $3$\\
$109$ & $240$ & $14$ & $3$\\
\hline
\end{tabular} \quad
\begin{tabular}{cccc}
\multicolumn{4}{c}{\small $m$-fold multiple edges (conjecture)}\\
\hline
$d$ & $n$ & $s$ & $t$\\
\hline
$6m-3$ & $12m$ & $6$ & $3$\\
$12m-7$ & $24m$ & $11$ & $3$\\
$12m-5$ & $24m$ & $11$ & $3$\\
$30m-19$ & $60m$ & $14$ & $3$\\
$30m-11$ \ & \ $60m$ \ & \ $14$ \ & \ $3$\\
\hline
\end{tabular}
\end{center}\quad
\caption{Graph with multiple edges from regular polyhedron}
\end{table}

\subsection{Graphs with multiple edges from semi-regular polyhedrons}\label{sec-mul-3_semireg_poly}
In this section, we consider graphs with multiple edges from semi-regular polyhedrons which are edge-transitive, cuboctahedron and icosidodecahedron (cf. Section \ref{sec-3_semireg_poly}). For the $(d, n, s, t)$-configurations of the vectors of the constructed crystal lattices, we have the following table:\vspace{-0.05in}

\begin{table}[h]
\begin{flushleft}
\qquad
\begin{tabular}{ccc}
& \qquad\\
&&\\
Cuboctahedron && $\cdots$\\
Icosidodecahedron && $\cdots$
\end{tabular}
\quad
\begin{tabular}{cccc}
\multicolumn{4}{c}{\small Graph with simple edge}\\
\hline
$d$ & $n$ & $s$ & $t$\\
\hline
$13$ & $48$ & $16$ & $3$\\
$31$ \ & \ $120$ \ & \ $33$ \ & \ $3$\\
\hline
\end{tabular}
\quad $\Leftrightarrow$
\end{flushleft}
\quad\\
\begin{center}
\begin{tabular}{cccc}
\multicolumn{4}{c}{Double edge}\\
\hline
$d$ & $n$ & $s$ & $t$\\
\hline
$37$ & $96$ & $18$ & $3$\\
$91$ & $240$ & $35$ & $3$\\
\hline
\end{tabular} \quad
\begin{tabular}{cccc}
\multicolumn{4}{c}{Triple edge}\\
\hline
$d$ & $n$ & $s$ & $t$\\
\hline
$61$ & $144$ & $18$ & $3$\\
$151$ & $360$ & $35$ & $3$\\
\hline
\end{tabular} \quad
\begin{tabular}{cccc}
\multicolumn{4}{c}{Quadruple edge}\\
\hline
$d$ & $n$ & $s$ & $t$\\
\hline
$85$ & $192$ & $18$ & $3$\\
$211$ & $480$ & $35$ & $3$\\
\hline
\end{tabular} \quad
\begin{tabular}{cccc}
\multicolumn{4}{c}{\small $m$-fold multiple edges (conjecture)}\\
\hline
$d$ & $n$ & $s$ & $t$\\
\hline
$24m-11$ & $48m$ & $18$ & $3$\\
$60m-29$ \ & \ $120m$ \ & \ $35$ \ & \ $3$\\
\hline
\end{tabular}
\end{center}\quad
\caption{Graph with multiple edges from semi-regular polyhedron}
\end{table}

\end{document}